\documentclass[reqno,11pt]{amsart}
\usepackage{amssymb,amsmath,amsthm,amsfonts,color}

\usepackage{mathrsfs,dsfont,comment,mathscinet}

\usepackage{enumerate,esint}
\usepackage[left=2.8cm,right=2.8cm,top=2.8cm,bottom=2.8cm]{geometry}
\usepackage{mathtools}
\usepackage{graphicx}

\usepackage{epstopdf}
\usepackage[colorlinks=true, pdfstartview=FitV, linkcolor=blue, 
            citecolor=blue, urlcolor=blue]{hyperref}

\usepackage[T1]{fontenc}
\usepackage{lmodern}
\usepackage{microtype}
\usepackage{amsthm}
\theoremstyle{plain}
\newtheorem{theorem}{Theorem}[section]
\newtheorem{proposition}[theorem]{Proposition}
\newtheorem{lemma}[theorem]{Lemma}
\newtheorem{coro}[theorem]{Corollary}

\theoremstyle{definition}
\newtheorem{definition}[theorem]{Definition}
\newtheorem{remark}[theorem]{Remark}
\newtheorem*{remark*}{Remark}
\newtheorem{hypo}[theorem]{Hypothesis}
\newtheorem*{hypo*}{Hypothesis}
\newtheorem{example}[theorem]{Example}

\theoremstyle{plain}
\newtheorem{result}{Metatheorem}[section]

\usepackage{amsmath,amssymb,mathrsfs,mathtools,braket,latexsym,esint}
\numberwithin{equation}{section}

\def\nat{\mathbb{N}}
\def\inte{\mathbb{Z}}

\def\real{\mathbb{R}}

\def\io{\int_\Omega}
\def\dd{\,\mathrm{d}}

\def\cH{\mathcal{H}}
\def\cL{\mathcal{L}}
\def\cW{\mathcal{W}}

\DeclareMathOperator{\PV}{P.V.}

\usepackage{enumerate}
\usepackage[shortlabels]{enumitem}

\def\red#1{\textcolor{red}{#1}}

\title[Existence results for singular nonlocal Cahn-Hilliard systems]{Existence results for Cahn-Hilliard-type systems\\
driven by nonlocal integrodifferential operators\\
with singular kernels}
\date{\today}
\author[E. Davoli]{Elisa Davoli}
\author[C. Gavioli]{Chiara Gavioli}
\author[L. Lombardini]{Luca Lombardini}	

\AtEndDocument{
	\bigskip
	\footnotesize{
        \noindent
		{\sc E. Davoli}, Institute of Analysis and Scientific Computing, TU Wien, Wiedner Hauptstra\ss e 8-10, 1040 Vienna, Austria.
		E-mail: \href{mailto:elisa.davoli@tuwien.ac.at}{\tt elisa.davoli@tuwien.ac.at}.

        \medskip
  
 	\noindent
		{\sc C. Gavioli}, Faculty of Civil Engineering, Czech Technical University, Th\'akurova 7, 16629 Prague 6, Czech Republic.	
		E-mail: \href{mailto:chiara.gavioli@cvut.cz}{\tt chiara.gavioli@cvut.cz}.

        \medskip

        \noindent
		{\sc L. Lombardini}, Institute of Analysis and Scientific Computing, TU Wien, Wiedner Hauptstra\ss e 8-10, 1040 Vienna, Austria.
		E-mail: \href{mailto:luca.lombardini@tuwien.ac.at}{\tt luca.lombardini@tuwien.ac.at}.
	}
}

\begin{document}

	
\begin{abstract}

    We introduce a fractional variant of the Cahn-Hilliard equation settled in a bounded domain and with a possibly singular potential. We first focus on the case of homogeneous Dirichlet boundary conditions, and show how to prove the existence and uniqueness of a weak solution. The proof relies on the variational method known as \textit{minimizing movements scheme}, which fits naturally with the gradient-flow structure of the equation. The interest of the proposed method lies in its extreme generality and flexibility. In particular, relying on the variational structure of the equation, we prove the existence of a solution for a general class of integrodifferential operators, not necessarily linear or symmetric, which include fractional versions of the $q$-Laplacian.
    
    In the second part of the paper, we adapt the argument in order to prove the existence of solutions in the case of regional fractional operators.    
    As a byproduct, this yields an existence result in the interesting cases of homogeneous fractional Neumann boundary conditions or periodic boundary conditions. 
  
		\medskip
  
		

\end{abstract}
	
	\maketitle
	
	{\parskip=0em \tableofcontents}

\section{Introduction}

Since its introduction in \cite{cahn1958free}, the Cahn-Hilliard equation has acquired a pivotal role in the modeling of spinodal decomposition phenomena in binary alloys, as well as in the diffuse-interface modeling of evolutionary problems ranging from physics and biology to imaging and materials science. Indeed, as observed in the nice Introduction of \cite{akagi2016fractional}, Cahn-Hilliard systems find application, for example, in the description of phase transition and separation phenomena, viscoelasticity, damaging, complex fluids, tumor growth, and whenever diffuse interfaces appear. In its classical formulation, on a bounded domain $\Omega\subset\real^d$ with Lipschitz boundary, and on a time-interval $(0,T)$, the Cahn-Hilliard equation takes the form
\begin{align}
\partial_t u - \Delta w = 0 &\qquad \text{in } \Omega \times (0,T),\label{eq:CH1local} \\
w = -\Delta u + F'(u) &\qquad \text{in } \Omega \times (0,T),\label{eq:CH2local} 
\end{align}
complemented by suitable initial data and boundary conditions, where $u:\Omega\to \real$ is a concentration (or order) parameter, $F$ is a double-well potential, $F'$ is its differential, and the chemical potential $w$ is computed as the differential of the corresponding energy functional 
\begin{equation}\label{eq:intro_local_free_energy}
E(u)\coloneqq \int_\Omega |\nabla u(x)|^2\dd x+\int_\Omega F(u(x))\dd x.
\end{equation}
For an overview of the mathematical aspects, we refer to the survey \cite{miranville2017cahn}.

In this paper, we consider the problem
\begin{align}
\partial_t u + \mathfrak{L} w = 0 &\qquad \text{in } \Omega \times (0,T), \label{eq:CH1intro}\\
w = \mathfrak{I}u + F'(u) &\qquad \text{in } \Omega \times (0,T), \label{eq:CH2intro} \\
u(x,0) = u_0(x) &\qquad \text{in } \real^d. \label{eq:iniCHintro}
\end{align}
Here above $\mathfrak{I}$ is a nonlocal integrodifferential operator with a non-integrable kernel (the examples to keep in mind are the fractional Laplacian and the $(s,q)$-Laplacian, and their regional counterparts), whereas $\mathfrak L$ is a linear invertible operator having a variational nature, like the classical or the (regional) fractional Laplacian.
A motivation for considering a nonlocal operator in place of the Laplacian, in equation \eqref{eq:CH2local}, can be traced back to the formulation of the physical model \cite{cahn1958free}. Roughly speaking, the reason is that nonlocal operators take into consideration also mid- and long-range interactions between particles, and this seems to better capture the behavior of the phenomena under consideration. The first rigorous derivation of a nonlocal Cahn-Hilliard model, based on stochastic arguments, can be found in \cite{giacomin1997phase}. We refer to Section~\ref{sec:nl_liter} for a brief overview of existing papers on the nonlocal Cahn-Hilliard problem, and their relationship with our results.

Typically, the function $F$ in \eqref{eq:intro_local_free_energy} represents a configuration potential with two wells, including possibly degenerate double-well potentials defined on bounded domains. In our setting we allow for more general potentials, not necessarily bounded from below (hence not necessarily of double-well-type) and decreasing to $-\infty$ in a controlled way. Moreover, $F$ can be singular, albeit with a certain form which still includes the usually considered potentials which we recall here below. To be precise, the classical choice for $F$ is the fourth-order polynomial
\begin{equation}\label{eq:cheetahlo}
F_{\rm pol}(z) \coloneqq \frac14 (1-z^2)^2, \quad z \in \real,
\end{equation}
with minima $z = \pm 1$ corresponding to the pure phases. It is well-known that, in view of the physical interpretation of the model, a more realistic description is given by the logarithmic double-well potential
\begin{equation}\label{eq:Flog}
F_{\rm log}(z) \coloneqq \frac{\theta}{2}\,\Big((1+z)\log(1+z) + (1-z)\log(1-z)\Big) + \frac{\theta_c}{2} - cz^2
\end{equation}
for $0 < \theta < \theta_c$ and $c > \theta/2$, which is defined on the bounded domain $(-1,1)$ and possesses two minima within this open interval. Another interesting example of $F$ which is covered by our analysis is the so-called double-obstacle potential, having the form
\begin{equation}\label{eq:Fob}
F_{\rm ob}(z) \coloneqq I_{[-1,1]}(z) + \frac12(1-z^2), \qquad
I_{[-1,1]}(z) \coloneqq
\begin{cases}
0 &\text{if } z \in [-1,1]\\
+\infty &\text{otherwise}.
\end{cases}
\end{equation}
In this latter case the derivative $F'_{\rm ob}$ has to be interpreted as the subdifferential $\partial F_{\rm ob}$ in the sense of convex analysis, and equation \eqref{eq:CH2intro} must be read as a differential inclusion.

In this setting, we then prove the existence of a weak solution to system \eqref{eq:CH1intro}--\eqref{eq:iniCHintro}. The main features and contributions of the paper can be summarized as follows:
\begin{itemize}
    \item we are able to consider a wide variety of potentials, without having to develop ad-hoc techniques for each one;
    \item the operator $\mathfrak L$ in equation \eqref{eq:CH1intro} can be either local, like the classical Laplacian, or nonlocal, like the fractional Laplacian, or a sum of the two. It can also be a higher order (integro)differential operator;
    \item the integrodifferential operators that we consider are quite general and can be either ``global'', when associated to homogeneous Dirichlet conditions, or regional;
    \item the operator $\mathfrak I$ in equation \eqref{eq:CH2intro} can be nonlinear, like the $(s,q)$-Laplacian with $q>2$;
    \item the variational argument that we employ in the proof of the existence results directly provides also energy estimates satisfied by the solutions;
    \item we prove the existence of a solution to the system in which $\mathfrak L=(-\Delta)^\sigma$ and $\mathfrak I=(-\Delta)^s$, combined with the homogeneous fractional Neumann conditions introduced in \cite{dipierro2017nonlocal}.
\end{itemize}


\subsection{Overview of the main results}

We now give a more detailed description of system \eqref{eq:CH1intro}--\eqref{eq:iniCHintro}. Yet, given its technical nature, we refrain here from providing the full functional theoretic framework needed to state our results, and we refer instead the reader to the appropriate sections for a rigorous treatment. In what follows, unless otherwise specified, $\Omega\subset\real^d$ is a bounded open set with Lipschitz boundary.

In place of the free energy \eqref{eq:intro_local_free_energy}, we consider a nonlocal energy having either the form
$$
\mathscr E(u)=\iint_{Q(\Omega)}\Phi(u(x)-u(y))K(x,y)\dd x\dd y+\int_\Omega F(u(x))\dd x,
$$
or
$$
\mathscr E_\Omega(u)=\iint_{\Omega^2}\Phi(u(x)-u(y))K_\Omega(x,y)\dd x\dd y+\int_\Omega F(u(x))\dd x,
$$
where $Q(\Omega)\coloneqq\real^{2d}\setminus\big((\real^d\setminus\Omega)\times(\real^d\setminus\Omega)\big)$, and $\Omega^2\coloneqq\Omega\times\Omega$. For the purposes of this introduction, the reader can assume that $\Phi(r)=|r|^q$, with $q\geq2$. The precise hypotheses are given in Section~\ref{sec:sqdef}, see in particular \eqref{eq:nameless}. Still, we stress that, in general, $\Phi$ need not be convex or even.

As anticipated above, the potential $F$ can be singular, provided that it has a suitable form. Roughly speaking, we require $F$ to be semiconvex (see Remark~\ref{rmk:chao}), meaning that it can be written as
$$
F=\Gamma+\Pi, \mbox{ where }\Gamma\mbox{ is convex and }\Pi(z) \coloneqq \int_0^z \pi(r)\dd r,
$$
with $\pi$ a globally Lipschitz continuous function. We also allow $F$ to decrease to $-\infty$, but only in a ``subcritical'' way. We refer to Hypothesis~\ref{stm:hyp}\,(i)--(iii) for the rigorous assumptions and to Remark~\ref{rmk:gen_hyp_potential} and Remark~\ref{rmk:barbieref} for some further comments and generalizations.

If we denote $\phi\coloneqq\Phi'$, the operator $\mathfrak I$ corresponding (see Lemma~\ref{lem:nonhouncomodino} and Lemma~\ref{lem:flareon}) to the double integrals in $\mathscr E$ and $\mathscr E_\Omega$ is then defined, respectively, as
\begin{equation}\label{eq:Iintro}
\langle \mathfrak I u,v\rangle\coloneqq \iint_{Q(\Omega)}\phi(u(x)-u(y))(v(x)-v(y))K(x,y)\dd x\dd y,
\end{equation}
and
\begin{equation}\label{eq:IOintro}
\langle \mathfrak I_\Omega u,v\rangle\coloneqq \iint_{\Omega^2}\phi(u(x)-u(y))(v(x)-v(y))K_\Omega(x,y)\dd x\dd y.
\end{equation}

The relevant properties of the associated functional spaces are described in Section~\ref{sec:HDbc} and in Section~\ref{sec:chem}.

\subsubsection{Homogeneous Dirichlet conditions}
Let us focus first on the operator $\mathfrak I$ given in \eqref{eq:Iintro}. On the kernel $K$ we impose two competing conditions (for the precise assumptions, see \eqref{eq:drogasintetica}--\eqref{eq:Kintegrability}), which are similar to those satisfied by (singular) L\'evy-type kernels (see, e.\,g., \cite{felsinger2015dirichlet} and Remark~\ref{rmk:KassmannLevy}). On the one hand, $K$ is assumed to be non-integrable in a neighborhood of the diagonal $\{x=y\}$. On the other hand, we assume that an appropriate renormalization of $K$ is integrable on $Q(\Omega)$. Roughly speaking, this second condition amounts to assuming that the double integral in the energy $\mathscr E(u)$ is finite when $u$ is a Lipschitz function (see Lemma~\ref{lem:AloeNatale}).
When $\phi$ is odd and $K$ is symmetric, the operator $\mathfrak I$ has also, at least formally, a pointwise representation as the integrodifferential operator
$$
\mathfrak I u(x)=\PV \int_{\real^d}\phi(u(x)-u(y))K(x,y)\dd y.
$$
The simplest and most widely studied example of such an operator is the $(s,q)$-Laplacian $(-\Delta)^s_q$, which corresponds to the choices $\phi(r)=|r|^{q-2}r$ and $K(x,y)=|x-y|^{-d-sq}$, with $s\in(0,1)$.
We refer the interested reader to the two very well-written surveys \cite{abatangelo2019getting,garofalo2017fractional} on the fractional Laplacian (obtained when $q=2$), and to the papers \cite{kuusi2015nonlocal,di2016local}, and the references cited therein, for the nonlinear case $q\not=2$.
We also mention \cite{vazquez2014recent} for an overview on evolution equations involving these operators.

We couple system \eqref{eq:CH1intro}--\eqref{eq:iniCHintro} with the homogeneous Dirichlet conditions
\begin{align}
u = 0 &\qquad \text{in } (\real^{d}\setminus\Omega)\times(0,T),\label{eq:DirIntro1}\\
w = 0 &\qquad \text{in } (\real^{d}\setminus\Omega)\times(0,T)\label{eq:DirIntro2}.
\end{align}
Given the nature of $\mathfrak I$, it is indeed necessary to prescribe $u$ on the whole complement $\real^d\setminus\Omega$, rather than only on the boundary $\partial\Omega$. From the functional theoretic point of view, we represent conditions \eqref{eq:DirIntro1} and \eqref{eq:DirIntro2} by considering only functions that belong to the space
\begin{equation}\label{eq:L20}
\cL^2_0 \coloneqq \left\lbrace v \in L^2(\real^d) \,:\, v = 0 \text{ a.\,e. in } \real^d\setminus\Omega \right\rbrace \quad \text{with } \|\cdot\|_{\cL^2_0} \coloneqq \|\cdot\|_{L^2(\real^d)}.
\end{equation}

The operator $\mathfrak L$ in equation \eqref{eq:CH1intro} is then a linear invertible operator $\mathfrak L:X_0\to(X_0)'$, defined on a Hilbert space $X_0$ that is densely embedded in $\cL^2_0$. The precise assumptions are given in Section~\ref{sec:dual}. This abstract formulation of equation \eqref{eq:CH1intro} is motivated by the possibility of considering both the cases $\mathfrak L=-\Delta$ and $\mathfrak L=(-\Delta)^\sigma$ in a unified way, within the same framework.

We can now state our first main result in the following non-formal guise. 

\begin{result}\label{stm:Metathm1}
    Given a function $u_0\in\cL^2_0$ with $\mathscr E(u_0)<+\infty$, there exists a triple of functions $(u,w,\zeta)$, with $u$ and $w$ satisfying \eqref{eq:DirIntro1} and \eqref{eq:DirIntro2} respectively, and $\zeta\in\partial\Gamma(u)$ almost everywhere in $\Omega\times(0,T)$, that solves, in an appropriate weak sense, the system
    \begin{align*}
\partial_t u + \mathfrak L w = 0 &\qquad\mbox{in }\Omega,\\
w = \mathfrak I u + \zeta + \pi(u) &\qquad\mbox{in }\Omega,
\end{align*}
almost everywhere in $(0,T)$. Moreover, $u(0) = u_0$  and
$$
\frac{1}{2}\int_0^t \|w(\tau)\|_{X_0}^2 \dd \tau +  \mathscr E(u(t)) \le \mathscr E(u_0)\qquad\mbox{for every }t\in[0,T).
$$
\end{result}

The rigorous definition of weak solution and statement of the existence result are given, respectively, in Definition~\ref{def:weaksolq} and Theorem~\ref{stm:existenceq}. We additionally observe that if the operator $\mathfrak I$ is ``strongly monotone'' (more precisely, if $\phi$ satisfies \eqref{eq:nameless1} instead of the lower bound in \eqref{eq:nameless}), then the solution is also unique.

We separately first consider the case corresponding to the choices $\mathfrak L=-\Delta$ and $\mathfrak I=(-\Delta)^s$, in Section~\ref{sec:LapFracLap}. The motivation is twofold: we believe that in this simpler, well-known, functional theoretic setting, the non-expert reader can better follow the essential steps of the proof, and it allows us to keep track of the dependence on the fractional index $s$ in the constants appearing in the various energy estimates. As a byproduct, we can then prove that, as $s\to 1^-$, solutions converge to a weak solution to the classical Cahn-Hilliard system \eqref{eq:CH1local}--\eqref{eq:CH2local}, see Theorem~\ref{stm:asymptotics} (and Remark~\ref{rmk:tartare} for some comments on related asymptotics results).

It is also worth mentioning that, by choosing $\mathfrak L$ to be the Riesz map associated to $\cL^2_0$, Metatheorem~\ref{stm:Metathm1} provides a solution to the Allen-Cahn equation, see Section~\ref{sec:AllenCahn}.

Moreover, $\mathfrak L$ can also be a higher order fractional Laplacian, as observed in Section~\ref{sec:high}.

\subsubsection{Regional operators}
Let us now consider the operator $\mathfrak I_\Omega$ given in \eqref{eq:IOintro}. The kernel $K_\Omega$ is assumed to satisfy conditions similar to the ones imposed on the kernel $K$ considered above, see \eqref{eq:KappaOm}--\eqref{eq:sticweak}. We still assume $K_\Omega$ to be singular in a neighborhood of the diagonal $\{x=y\}$, and we ask an appropriate renormalization to be, roughly speaking, locally summable in $\Omega^2$. We refer to Section~\ref{sec:kernels} for some interesting examples of regional kernels. In particular, assuming local summability, rather than summability in the whole $\Omega^2$ (formally, assuming \eqref{eq:sticweak} in place of \eqref{eq:stic}), allows us to consider also the kernel in (K2) defined in Section~\ref{sec:kernels}, which roughly corresponds to the choice of periodic boundary conditions for the $(s,q)$-Laplacian. Notice that, in this setting, the functions that we consider are only defined in $\Omega$ rather than on the whole space $\real^d$.

The basic example to keep in mind is the regional $(s,q)$-Laplacian, which has the pointwise representation
$$
\mathfrak I_\Omega u(x)=\PV\int_\Omega |u(x)-u(y)|^{q-2}(u(x)-u(y))\frac{\dd y}{|x-y|^{d+sq}}.
$$
We refer the interested reader to \cite{fall2022regional} and the references cited therein, for the case $q=2$. These regional operators arise naturally when considering ``global'' operators, like the fractional Laplacian, and, roughly speaking, not allowing jumps outside of $\Omega$. This can be made precise in the context of censored stable processes, see, e.\,g., \cite{bogdan2003censored}.

We then consider the system \eqref{eq:CH1intro}--\eqref{eq:iniCHintro}, with a regional operator $\mathfrak I=\mathfrak I_\Omega$, without adding other boundary conditions. 
This brings in some difficulties. The main issue consists of having to consider also non-trivial constant functions. Since these belong to the kernel of $\mathfrak I_\Omega$, variational arguments based on coercivity and compactness properties then become more delicate.
From the technical point of view, we deal with this complication by working in subspaces of functions that have fixed mass. Instead of $\cL^2_0$, we consider here the space $L^2(\Omega)$ and its subspace of functions having mass zero,
\begin{equation}\label{eq:L2punto}
\dot{L}^2(\Omega)\coloneqq\left\{v\in L^2(\Omega)\,:\,\mathfrak m(v)\coloneqq\frac{1}{|\Omega|}\int_\Omega v(x)\dd x=0\right\}.
\end{equation}

The definition of the operator $\mathfrak L$ in equation \eqref{eq:CH1intro} is quite technical, see Hypothesis~\ref{stm:lumpa}. Essentially, we consider $\mathfrak L$ to be a linear operator defined on a Hilbert space $X\subset L^2(\Omega)$ and we split $X=X_0\oplus\real$, with $X_0=X\cap \dot{L}^2(\Omega)$. Then, we assume that $\mathfrak L$ is invertible on $X_0$ while being trivial on the constant functions.
The reader should have in mind the case $\mathfrak L=-\Delta$ defined on $X=H^1(\Omega)$.

Under these hypotheses, our second main result can be stated as follows.

\begin{result}\label{stm:Metathm2}
Let $\Omega\subset\real^d$ be a bounded and connected open set with Lipschitz boundary. Given an appropriate function $u_0$, there exists a triple of functions $(u,w,\zeta)$, with $\zeta\in\partial\Gamma(u)$ almost everywhere in $\Omega\times(0,T)$, that solves, in an appropriate weak sense, the system
    \begin{align*}
\partial_t u + \mathfrak L w = 0 &\qquad\mbox{in }\Omega,\\
w = \mathfrak I_\Omega u + \zeta + \pi(u) &\qquad\mbox{in }\Omega,
\end{align*}
almost everywhere in $(0,T)$. Moreover, $u(0) = u_0$,
$$
\mathfrak m(u(t))=\mathfrak m(u_0)\qquad\mbox{for every }t\in[0,T),
$$
and
$$
\frac{1}{2}\int_0^t \|w(\tau)-\mathfrak m(w(\tau))\|_{X_0}^2 \dd \tau +  \mathscr E_\Omega(u(t)) \le \mathscr E_\Omega(u_0)\qquad\mbox{for every }t\in[0,T).
$$
\end{result}
Moreover, if the operator $\mathfrak I_\Omega$ is ``strongly monotone'' (more precisely, if $\phi$ satisfies \eqref{eq:nameless1} instead of the lower bound in \eqref{eq:nameless}), then $u$ is unique. On the other hand, the uniqueness of $w$ and $\zeta$ depends on the differentiability properties of $\Gamma$, see Remark~\ref{rmk:stj}\,(3). For the rigorous statement, see Definition~\ref{def:weaksolregL} and Theorem~\ref{stm:existenceregL}.

An interesting feature of regional operators is that they sometimes provide a way to reformulate problems involving ``global'' operators like the $(s,q)$-Laplacian. The simplest example is given by considering $\Omega=(0,1)^d$ and $\mathbb Z^d$-periodic functions. Appropriately splitting $\real^d$ into cubes, and taking advantage of the periodicity to perform a change of variables, one can then rewrite the $(s,q)$-Laplacian as a regional operator, with an appropriate kernel $K_\Omega$. We refer the reader to Section~\ref{sec:PeriodicBdSec}.

Another very important example is given by the fractional Laplacian, when associated to the fractional Neumann boundary condition introduced in \cite{dipierro2017nonlocal}. Motivated by variational considerations and by integration by parts formulas, the authors introduced the operator $\mathcal N_s u$ as a nonlocal counterpart to the normal derivative $\partial_\nu u$. Up to constants, this is defined as
$$
\mathcal N_s u(x)\coloneqq\int_\Omega\frac{u(x)-u(y)}{|x-y|^{d+2s}}\dd y\qquad\mbox{for }x\in\real^d\setminus\overline{\Omega}.
$$
In \cite{abatangelo2020remark} the author observed that if $\mathcal N_s u=0$ almost everywhere in $\real^d\setminus\Omega$, then the fractional Laplacian of $u$ in $\Omega$ can be rewritten as a regional operator,
\begin{equation}\label{eq:NeuIdeIntro}
(-\Delta)^s u(x)=\PV\int_\Omega (u(x)-u(y))K_\Omega^s(x,y)\dd y\quad\mbox{for every }x\in\Omega,
\end{equation}
where $K_\Omega^s$ denotes the kernel (K3) in Section~\ref{sec:kernels}. This observation was further explored in the context of weak solutions in \cite{audrito2022neumann}. Exploiting an estimate given in \cite[Proposition~2.1]{audrito2022neumann}, we prove in Section~\ref{sec:K3} that the kernel $K_\Omega^s$ fits into our functional theoretic framework.
As a consequence, we obtain the following existence result.

\begin{result}\label{stm:Metathm3}
Let $\Omega\subset\real^d$ be a bounded (not necessarily connected) open set with $C^2$ boundary. Then, given an appropriate function $u_0$, there exists a triple of functions $(u,w,\zeta)$, with $\zeta\in\partial\Gamma(u)$ almost everywhere in $\Omega\times(0,T)$, that solves, in an appropriate weak sense, the system
    \begin{align*}
\partial_t u + (-\Delta)^\sigma w = 0 &\qquad \text{in } \Omega\times(0,T),\\
w = (-\Delta)^s u + \zeta + \pi(u) &\qquad \text{in } \Omega\times(0,T),\\
\mathcal N_\sigma w = 0 &\qquad \text{in } (\real^{d}\setminus\Omega)\times(0,T),\\
\mathcal N_s u = 0 &\qquad \text{in } (\real^{d}\setminus\Omega)\times(0,T),\\
u(x,0) = u_0(x) &\qquad \text{in } \Omega.
\end{align*}
Moreover,
$$
\mathfrak m(u|_\Omega(t))=\mathfrak m(u_0|_\Omega)\qquad\mbox{for every }t\in[0,T),
$$
and appropriate energy estimates hold true.
\end{result}

The Cahn-Hilliard system in Metatheorem~\ref{stm:Metathm3} is interpreted in its natural weak formulation, given in Definition~\ref{def:weaksolNeu}. In Theorem~\ref{stm:Teo} we prove that this is equivalent to the weak formulation in the framework of regional operators, that is obtained via identity \eqref{eq:NeuIdeIntro}. Metatheorem~\ref{stm:Metathm3}, whose rigorous statement is given in Corollary~\ref{cor:NeuExi}, can then be obtained as a consequence of Metatheorem~\ref{stm:Metathm2}. We stress that $u$ is unique, while the uniqueness of $w$ and $\zeta$ depends on the differentiability properties of $\Gamma$ as for Metatheorem~\ref{stm:Metathm2}.

Other operators that can be written as regional operators that fit into our framework are certain spectral fractional Laplacians, as noted in Section~\ref{sec:geniodelletartarughe}.

\subsubsection{Scheme of the proof}
We now give an explanation of the main ideas involved in the proof of our existence results.

The proof of Metatheorem~\ref{stm:Metathm1} is inspired by the arguments in \cite{akagi2016fractional} and is based on the gradient flow structure of the Cahn-Hilliard system. However, rather than relying on the equations to derive the required energy estimates, as in \cite{akagi2016fractional}, here we take full advantage of the variational nature of the system. This allows us to consider much more general potentials and operators.

The rough idea of the proof can be described as follows.
First of all, we substitute the potential $F$ with the regularized potential $F_\lambda=\Gamma_\lambda+\Pi$, obtained by considering the Moreau regularization $\Gamma_\lambda$ in place of $\Gamma$. This procedure is well-known, and leaves us with a regular potential that satisfies nice approximating properties as $\lambda\to0$, see Section~\ref{sec:MorYos}.

As a second step, exploiting the linearity and invertibility of $\mathfrak L$, we formally rewrite \eqref{eq:CH1intro} as $w=-\mathfrak L^{-1}(\partial_t u)=-\partial_t(\mathfrak L^{-1}u)$. Substituting in \eqref{eq:CH2intro}, we are thus left with the equation
$$
\partial_t(\mathfrak L^{-1}u)+\mathfrak I u+F_\lambda'(u)=0.
$$
Then, the key observation is that the equation
$$
\mathfrak L^{-1}\Big(\frac{u-g}{\tau}\Big)+\mathfrak I u+F_\lambda'(u)=0
$$
is the Euler-Lagrange equation associated to the energy
$$
\mathscr E^\lambda_g(u)\coloneqq\frac{1}{2\tau}\|u-g\|^2_{(X_0)'}+\iint_{Q(\Omega)}\Phi(u(x)-u(y))K(x,y)\dd x\dd y+\int_\Omega F_\lambda(u(x))\dd x.
$$
This allows us to resort to the classical scheme given by De Giorgi's minimizing movements, to construct a solution. We discretize the time interval $(0,T)$ and we iteratively define the functions $u_n$ as minimizers of $\mathscr E^\lambda_{u_{n-1}}$, beginning with the initial data $u_0$.
In particular, the non-integrability of the kernel $K$ ensures both the validity of an appropriate fractional Poincar\'e inequality (see Proposition~\ref{prop:onthesofa}) and of the necessary compactness that, together with the growth assumption on the original potential $F$, allow us to prove the existence of such minimizers.

By appropriately gluing the functions $u_n$, and passing to the limit in the time-discretizing parameter, we find a solution $(u^\lambda,w^\lambda)$ of system \eqref{eq:CH1intro}--\eqref{eq:iniCHintro}, with $F_\lambda$ in place of $F$.

The final step is to show that these approximating solutions converge, as $\lambda\to0$, to a solution to \eqref{eq:CH1intro}--\eqref{eq:iniCHintro}.
Both convergence with respect to the time-discretizing parameter and convergence with respect to $\lambda$ are obtained as a consequence of uniform energy estimates that essentially follow from the definition of the functions $u_n$' as minimizers. We mention however that some care must be taken, in order to deal with the terms $\Gamma_\lambda'(u^\lambda)$, because of the possible singularity of $\Gamma$. This requires an additional $L^2$-estimate.

The rigorous argument is carried out in Section~\ref{sec:proofDirGenThm}.

The proof of Metatheorem~\ref{stm:Metathm2} follows the same scheme outlined above. The main difference is due to the appearance of the mass in the fractional Poincar\'e inequality, see Proposition~\ref{stm:lorenzanonce}.

To address this additional difficulty, we first show that the mass of $u$ must remain constant over time. Motivated by this observation, we temporarily restrict ourselves to consider functions $u$ whose mass is equal to $\mathfrak m(u_0)$, and functions $w$ having mass zero. Under these assumptions, we can carry out the time-discretization step of the above argument, ending up with an approximating pair $(u^\lambda,\omega^\lambda)$, with $\omega^\lambda$ having mass zero.

Then, the most delicate step consists in defining $w^\lambda$ by properly adjusting the mass of $\omega^\lambda$. Formally, this is done by defining $w^\lambda(t)\coloneqq \omega^\lambda(t)+\mathfrak m(F_\lambda(u^\lambda(t)))$. In order to pass to the limit $\lambda\to0$, we therefore need to prove uniform energy estimates for the term $\mathfrak m(F_\lambda(u^\lambda))$.
This is achieved by carefully adjusting the argument of \cite[Section~5]{kenmochi1995subdifferential}.

The rigorous proof is the content of Section~\ref{sec:ProofThmRegExi}.



\subsection{Related literature on existence for nonlocal Cahn-Hilliard problems}\label{sec:nl_liter}

We begin by noting that one of the first existence results for a rather general Cahn-Hilliard equation in Hilbert spaces was provided in \cite{kenmochi1995subdifferential}, where the authors reinterpreted the equation as a differential inclusion by viewing the sum of both $\mathfrak{I}u$ and $\gamma(u)$ as a subdifferential. It is thus implicit that, in this abstract formulation, one can take a nonlocal non-integrable operator like $\mathfrak I$ or $\mathfrak I_\Omega$ corresponding to the specific case in which $q=2$ and $\Phi$ is convex, hence obtaining the existence for the corresponding problem.

Explicitly nonlocal Cahn-Hilliard problems, taking into account mid- and long-range interactions, have been the object of a growing interest starting from the work \cite{giacomin1997phase}. We mention \cite{gajewski2003nonlocal,bates2005dirichlet,nec2008front,colli2012global} as the first occurrences of linear nonlocal operators defined as a convolution with symmetric and integrable kernels. In all these papers, the potential $F$ is assumed to be regular enough for their purposes.

Up to our knowledge, the first contribution where a linear operator defined via a non-integrable kernel was considered is \cite{abels2015cahn}. The form of such kernel $k$ is slightly more general than the one we consider here. However, $k$ is estimated both from above and below by $|x-y|^{-d-2s}$, with $s\in(1/2,1)$, and it is assumed to be $(d+2)$-times differentiable outside of the diagonal, with all derivatives being estimated from above by those of $|x-y|^{-d-2s}$. Additionally, they only consider a logarithmic-type potential and the operator $\mathfrak L$ in \eqref{eq:CH1intro} is the Laplacian.

While the setting in \cite{abels2015cahn} is only regional, in \cite{akagi2016fractional} the framework is ``global'', corresponding to homogeneous Dirichlet boundary conditions. Therein, the operators are $\mathfrak L = (-\Delta)^{\sigma}$ and $\mathfrak I = (-\Delta)^s$, with $\sigma,s \in (0,1)$, and the potential is of the type $F(u) = |u|^p/p-u^2/2$. Additionally, they consider various asymptotic limits and the behavior of certain stationary solutions. In \cite{akagi2019convergence} the same authors explore the long-time behavior of the system.

General regional linear operators other than the fractional Laplacian, defined via possibly non-integrable kernels, were considered in \cite{gal2017strong,gal2018doubly}. Our approach is rather different: while the abstract setting in the aforementioned papers is based on the theory of Dirichlet forms, we push the variational point of view to its limits without requiring the operators involved to be linear. Moreover, the explicit examples provided in \cite{gal2017strong} are ``very close'' to the usual regional fractional Laplacian. In \cite{gal2018doubly} the kernels are radially symmetric and satisfy a stronger integrability condition than ours. A priori, the non-integrability is instead more general since it is only imposed in an abstract sense by requiring the domain of the integrodifferential operator to be compactly embedded in $L^2(\Omega)$, but again, the examples provided are close to the standard ones. Additionally, the potentials considered in these two papers have to be smooth enough and satisfy suitable growth conditions.

A different notion of fractional operator, namely, a spectral one, was considered in \cite{colli2019well}, together with a potential having the same form as ours (but with at least quadratic growth). The problem considered is the Cahn-Hilliard equation, both of viscous or nonviscous type. The existence is obtained with a scheme similar to the one employed here, but with an additional regularization term to gain coercivity. The nonviscous case, with the specific choice of the operators defined as spectral fractional Laplacians (in $(0,1)^d$ with Dirichlet boundary conditions, or in a general $\Omega$ with Neumann boundary conditions) falls within our setting, as observed in Section~\ref{sec:geniodelletartarughe}. In this framework, we provide a more general existence result, since our potential satisfies much less restrictive growth conditions. Similar considerations hold for the paper \cite{gal2023separation}, where the operator in the second equation is the spectral fractional Laplacian corresponding to Neumann boundary conditions.

With the general variational approach employed in this paper we have thus extended, in a unified way, most of the previous literature, as far as existence results are concerned.
The case of homogeneous Dirichlet boundary conditions is mostly new, as the existence has been rigorously studied only in \cite{akagi2016fractional}.

The case of Neumann boundary conditions considered in Metatheorem~\ref{stm:Metathm3} is completely new. The existence of solutions with this kind of boundary condition is quite difficult to prove, and very few such results are available in the literature, even for simpler problems (like the heat equation).

Moreover, as just observed, in all the previous works the operators involved are linear, while the variational approach allows us to consider also nonlinear operators. Even for the Cahn-Hilliard problem in the local case, nonlinear operators in the second equation have almost never been considered. Indeed, to the best of our knowledge, the only papers concerning the Cahn-Hilliard equation with $q$-Laplacian are focused on stationary solutions \cite{takavc2009stationary,drabek2011manifolds}, or on the long time-behavior of certain solutions in the $1$-dimensional case \cite{folino2022generalized}.
In \cite{lisini2012cahn} the authors consider a local Cahn-Hilliard-type system with a nonlinearity in the first equation, arising from a possible coupling between $u$ and $w$, given by a mobility function $m(u)$. The existence argument is still based on minimizing movements, with a weighted transport metric in place of the $X_0'$-norm, and is suitable for mass preservation.

Finally, given the abstract formulation of equation \eqref{eq:CH1intro}, we can consider also operators like higher order fractional Laplacians, mixed local-nonlocal operators, and also spectral fractional operators.


\section{Preliminaries and assumptions}

In this section we collect the preliminary results, the definitions, and the assumptions that will be used throughout the paper.


\subsection{Reminiscences of duality}\label{sec:dual}

We recall here some well-known results concerning duality mappings.
Let $(X_0,(\,\cdot\,,\,\cdot\,)_0)$ be a Hilbert space, and consider a linear invertible operator $\mathfrak L:X_0\to(X_0)'$, which has a variational nature, meaning that there exists a bilinear map $B_{\mathfrak L}: X_0\times X_0\to \real$, symmetric,
$$
\mbox{bounded:}\quad |B_{\mathfrak L}(u,v)|\leq C \|u\|_0 \|v\|_0\quad\mbox{for every }u,v\in X_0,
$$
and
$$
\mbox{coercive:}\quad B_{\mathfrak L}(u,u)\geq c\|u\|_0^2\quad\mbox{for every }u\in X_0,
$$
such that $\mathfrak L u=B_{\mathfrak L}(u,\,\cdot\,)$ for every $u\in X_0$.
This defines a new inner product as
$$
(u,v)_{\mathfrak L} \coloneqq B_{\mathfrak L}(u,v)=\langle \mathfrak L u, v \rangle_0 \quad \mbox{for all } u,v \in X_0,
$$
hence the corresponding norm $\|\,\cdot\,\|_{\mathfrak L}$ by setting
$$
\|u\|_{\mathfrak L}^2 \coloneqq (u,u)_{\mathfrak L} \quad \mbox{for all } u \in X_0.
$$
The boundedness and coercivity of $B_{\mathfrak L}$ ensure that the norm $\|\,\cdot\,\|_{\mathfrak L}$ is equivalent to the original norm $\|\,\cdot\,\|_0$. 
This also defines an equivalent norm on the dual space $(X_0)'$ in the usual way, by setting
$$
\|M\|_{\mathfrak L^{-1}}\coloneqq \sup_{\substack{u\in X_0 \\ \|u\|_{\mathfrak L}\leq 1}}|Mu|.
$$
Note that $\mathfrak L$ is the Fréchet derivative of $\|\,\cdot\,\|_{\mathfrak L}^2/2$, and, by definition, it is also the Riesz map of the scalar product $(\,\cdot\, ,\,\cdot\,)_{\mathfrak L}$, so that
$$
\langle \mathfrak Lu, u \rangle_0=\|u\|_{\mathfrak L}^2=\|\mathfrak Lu\|_{\mathfrak L^{-1}}^2\quad\mbox{for every }u\in X_0,
$$
and
\begin{equation}\label{eq:lorenzace}
\|M\|_{\mathfrak L^{-1}}=\|\mathfrak L^{-1}M\|_{\mathfrak L}\quad\mbox{for every }M\in (X_0)'.
\end{equation}
Hence $\mathfrak L^{-1}$ is the Fréchet derivative of $\|\,\cdot\,\|_{\mathfrak L^{-1}}^2/2$, and $\mathfrak L$ is a linear isometry.

We make explicit the following computation, which we will use later. Given $M\in W^{1,2}(0,T;(X_0)')$ we have the identity
$$
\Big\langle M(t),\mathfrak L^{-1}\Big(\frac{\dd}{\dd t}M(t)\Big)\Big\rangle_0=\frac{1}{2}\frac{\dd}{\dd t}\|M(t)\|^2_{\mathfrak L^{-1}}
$$
for almost every $t\in(0,T)$. We begin by observing that there exists the Fr\'echet derivative $\frac{\dd}{\dd t}M(t)\in L(\real;(X_0)')$ for every $t\in(0,T)\setminus\Sigma$, with $|\Sigma|=0$.
Since $\mathfrak L$ is a linear isometry, for every $t\in(0,T)\setminus\Sigma$ there exists also the Fr\'echet derivative $\frac{\dd}{\dd t}(\mathfrak L^{-1}M)(t)\in L(\real;X_0)$, and $\frac{\dd}{\dd t}(\mathfrak L^{-1}M)(t)=\mathfrak L^{-1}\big(\frac{\dd}{\dd t}M(t)\big)$.
Moreover, for every $t\in(0,T)\setminus\Sigma$,
\begin{equation}\label{eq:sensata}\begin{aligned}
    \Big\langle M(t),\mathfrak L^{-1}\Big(\frac{\dd}{\dd t}M(t)\Big)\Big\rangle_0 &
    =\Big\langle M(t),\frac{\dd}{\dd t}(\mathfrak L^{-1}M)(t)\Big\rangle_0\\
    &
    =\Big\langle(\mathfrak L\circ\mathfrak L^{-1})M(t),\frac{\dd}{\dd t}(\mathfrak L^{-1}M)(t)\Big\rangle_0\\
    &
    =\Big((\mathfrak L^{-1}M)(t),\frac{\dd}{\dd t}(\mathfrak L^{-1}M)(t)\Big)_{\mathfrak L}\\
    &
    =\frac{1}{2}\frac{\dd}{\dd t}\|(\mathfrak L^{-1}M)(t)\|_{\mathfrak L}^2\\
    &
    =\frac{1}{2}\frac{\dd}{\dd t}\|M(t)\|_{\mathfrak L^{-1}}^2.
\end{aligned}\end{equation}

\begin{remark}\label{rmk:triple}
    In this paper we are interested in the case in which $X_0$ is a Hilbert space densely embedded either in $\cL^2_0$ or in $\dot{L}^2(\Omega)$, which are defined in \eqref{eq:L20} and \eqref{eq:L2punto}. In the first case, this yields the Hilbert triple $X_0\hookrightarrow \cL^2_0\simeq(\cL^2_0)'\hookrightarrow(X_0)'$, where the maps are as follows:
    \begin{itemize}
    \item $i:X_0\hookrightarrow \cL^2_0$ is the dense embedding we assumed to exist;
    \item the identification $\cL^2_0\simeq(\cL^2_0)'$ is given by the standard Riesz map $u\mapsto(u,\,\cdot\,)_{\cL^2_0}$;
    \item the embedding $(\cL^2_0)'\hookrightarrow(X_0)'$ is given by the adjoint of $i$, that is
    $$
    M\in (\cL^2_0)'\mapsto M\circ i\in (X_0)',
    $$
    which is injective since $i$ is a dense embedding.
    \end{itemize}
    In particular, we can interpret a function $u\in X_0$ as the element of the dual space $(\cL^2_0)'$ given by
    $$
    (i(u),v)_{\cL^2_0} \quad\mbox{for every }v\in \cL^2_0,
    $$
    and as an element of the dual space $(X_0)'$ by
    $$
    \langle u,v\rangle_0=(i(u),i(v))_{\cL^2_0} \quad\mbox{for every }v\in X_0.
    $$
    The same holds with $\cL^2_0$ replaced by $\dot{L}^2(\Omega)$ in the case of the dense embedding $i: X_0\hookrightarrow\dot{L}^2(\Omega)$.

As customary, in what follows we will forget about the embedding $i$ in the formulas.
\end{remark}


\subsection{The fractional Laplacian}

For $s \in (0,1)$ and $v$ in the Schwartz class $\mathcal{S}(\real^d)$ of the rapidly decaying functions at infinity, the \textit{fractional Laplacian} is defined as
$$
(-\Delta)^s v(x) \coloneqq C(s,d) \PV \int_{\real^d} \frac{v(x)-v(y)}{|x-y|^{d+2s}} \dd y.
$$
The constant $C(s,d)$ is such that
$$
\lim_{s \to 1^-} \frac{C(s,d)}{s(1-s)} = \frac{4d}{\omega_{d-1}},
$$
and is a normalization constant motivated by the equivalence with the pseudodifferential form of the fractional Laplacian, see, e.\,g., \cite[Proposition~3.3]{di2012hitchhiker}. For other equivalent definitions of the fractional Laplacian, see the interesting \cite{kwasnicki2017ten}. In this paper we will forget about the constant, except for its behavior as $(1-s)$ when $s\to1$, which we will need for Theorem~\ref{stm:asymptotics}. We will then consider
$$
\begin{aligned}
(-\Delta)^s v(x) &\coloneqq (1-s) \PV \int_{\real^d} \frac{v(x)-v(y)}{|x-y|^{d+2s}} \dd y\\
&
=\frac{1-s}{2} \PV \int_{\real^d} \frac{2v(x)-v(x+y)-v(x-y)}{|y|^{d+2s}} \dd y.
\end{aligned}
$$
The fractional Laplacian, as written in the second line above, is well-defined in the classical sense, and not only as a principal value, provided $v$ is regular enough, e.\,g., $v\in\mathcal S(\real^d)$.
We recall that $(-\Delta)^s$ tends to the classical Laplacian in a suitable way, as $s \to 1$.
\begin{lemma}[Proposition~4.4 in \cite{di2012hitchhiker}]\label{thm:deltasto1}
For any $\psi \in \mathcal{S}(\real^d)$, there holds
$$
\lim_{s \to 1^-} -(-\Delta)^s \psi = C(d)\Delta\psi,
$$
uniformly in $x\in\real^d$.
\end{lemma}

We shall prove a further result about the behavior of the weak fractional Laplacian $(-\Delta)^s$ as the index $s$ goes to $1$ in Lemma~\ref{lem:asympi}, which will play an important role in the proof of Theorem~\ref{stm:asymptotics}.

On the other hand, the natural variational setting for equations involving the operator $(-\Delta)^s$ is the space $W^{s,2}(\real^d) = H^s(\real^d)$, defined as the space of functions in $L^2(\real^d)$ having finite \textit{Gagliardo seminorm}
\begin{equation}\label{eq:patchouli}
[v]_{H^s(\real^d)}^2 \coloneqq \iint_{\real^{2d}} \frac{|v(x)-v(y)|^2}{|x-y|^{d+2s}} \dd x \dd y.
\end{equation}
More precisely, when considering the fractional Laplacian coupled with the boundary condition $u=0$ almost everywhere in $\real^d\setminus\Omega$, the natural functional space to consider for the variational framework is the following subspace of $H^s(\real^d)$:
$$
\cH^s_0 \coloneqq \left\lbrace v \in H^s(\real^d) \,:\, v = 0 \text{ a.\,e. in } \real^d\setminus\Omega \right\rbrace.
$$
Since in Theorem~\ref{stm:asymptotics} we are interested in the limit as $s\to1$, in place of the seminorm \eqref{eq:patchouli}, for functions in the subspace $\cH^s_0$ we consider its normalized version
$$
\|v\|^2_{\cH^s_0} \coloneqq (1-s) \iint_{Q(\Omega)} \frac{|v(x)-v(y)|^2}{|x-y|^{d+2s}} \dd x \dd y = (1-s)[v]_{H^s(\real^d)}^2,
$$
where $Q(\Omega) \coloneqq \real^{2d} \setminus \left((\real^d \setminus \Omega) \times (\real^d \setminus \Omega)\right)$. In light of the fractional Poincar\'{e} inequality, it is well-known that $\|\,\cdot\,\|_{\cH^s_0}$ defines a norm on the space $\cH^s_0$, see also \cite[Subsection~2.3]{akagi2016fractional}. Endowed with this norm, $\cH^s_0$ is a Hilbert space and the scalar product is
$$
(u,v)_{\cH^s_0} \coloneqq (1-s) \iint_{\real^{2d}} \frac{(u(x)-u(y))}{|x-y|^{d+2s}}\,(v(x)-v(y)) \dd x \dd y.
$$
In analogy to the classical case, one may consider the Hilbert triple
$$
\cH^s_0 \hookrightarrow \cL^2_0 \simeq (\cL^2_0)' \hookrightarrow (\cH^s_0)',
$$
where we recall that the space $\cL^2_0$, which was defined in \eqref{eq:L20}, is given by
$$
\cL^2_0 \coloneqq \left\lbrace v \in L^2(\real^d) \,:\, v = 0 \text{ a.\,e. in } \real^d\setminus\Omega \right\rbrace \quad \text{with } \|\cdot\|_{\cL^2_0} \coloneqq \|\cdot\|_{L^2(\real^d)}.
$$
The identification is made by means of the standard scalar product on $L^2(\real^d)$ (and therefore on its closed subspace $\cL^2_0$)
$$
(u,v)_{L^2} \coloneqq \int_{\real^d} u(x)v(x) \dd x \quad \text{for } u,v \in L^2(\real^d).
$$
On the other hand, it is easy to verify that $C^\infty_c(\Omega) \subset \cH^s_0$, hence the embedding $\cH^s_0 \hookrightarrow \cL^2_0$ is dense. Moreover, it is also compact by the fractional Poincar\'{e} inequality and \cite[Theorem~7.1]{di2012hitchhiker}.
Based on this functional framework, we can introduce the weak form of the fractional Laplacian $(-\Delta)^s$. More precisely, the operator $(-\Delta)^s : \cH^s_0 \to (\cH^s_0)'$ is defined by
\begin{equation}\label{eq:w-slapl}
\left\langle (-\Delta)^s v,\psi \right\rangle_{(\cH^s_0)'} \coloneqq (1-s)\iint_{\real^{2d}} \frac{(v(x)-v(y))}{|x-y|^{d+2s}}\,(\psi(x)-\psi(y)) \dd x \dd y = (v,\psi)_{\cH^s_0}
\end{equation}
for $v,\psi \in \cH^s_0$.


\subsection{Variational nonlocal operators}\label{sec:sqdef}

A more general type of integrodifferential operators, not necessarily (bi)linear or symmetric, commonly found in the literature, have the form given in \eqref{eq:Iintro}, which is the following:
$$
\langle\mathfrak I u,v\rangle \coloneqq
\iint_{Q(\Omega)}\phi(u(x)-u(y))(v(x)-v(y))K(x,y)\dd x\dd y.
$$
The simplest example is the $(s,q)$-Laplacian, which is obtained by considering $\phi(r) = |r|^{q-2}r$ and $K(x,y) = |x-y|^{-d-sq}$ (up to renormalization constants). This variational formulation may not have a pointwise counterpart. However, if $\phi$ is odd and $K$ is symmetric, i.\,e., $K(x,y)=K(y,x)$, then, at least formally, we can write
$$
\langle \mathfrak I u,v \rangle = 2 \int_{\real^d} v(x) \bigg(\int_{\real^d} \phi(u(x)-u(y))K(x,y)\dd y\bigg) \dd x,
$$
for $v \in C^\infty_c(\Omega)$. Then, the pointwise definition of the operator is
$$
\mathfrak I u(x) = \PV \int_{\real^d} \phi(u(x)-u(y))K(x,y)\dd y.
$$

More in general, we consider as $\phi$ a function with two-sided $(q-1)$-growth and we only require the kernel $K$ to satisfy two integrability hypotheses.
In particular, $\phi$ need not be odd and $K$ need not be symmetric.

The precise assumptions on $\phi$ are as follows. We require $\phi:\real\to\real$ to be a continuous function, with $\phi(0)=0$, and satisfying the bound
\begin{equation}\label{eq:nameless}
\frac{1}{\Lambda}|r|^q\leq\phi(r)r\leq\Lambda|r|^q \quad\mbox{for every }r\in\real,
\end{equation}
for some $q\geq2$ and $\Lambda\geq1$. Sometimes we will assume the following stronger variant of the lower bound
\begin{equation}\label{eq:nameless1}
\frac{1}{\Lambda}|r_1-r_2|^q\leq(\phi(r_1)-\phi(r_2))(r_1-r_2) \quad\mbox{for every }r_1,r_2\in\real.
\end{equation}
We define $\Phi:\real\to\real$ as $\Phi(z)=\int_0^z \phi(r)\dd r$ and we observe that 
\begin{equation}\label{eq:heaviness}
\frac{1}{\Lambda \, q}|r|^q\leq\Phi(r)\leq\frac{\Lambda}{q}|r|^q \quad\mbox{for every }r\in\real.
\end{equation}
\begin{remark}
Estimate \eqref{eq:nameless1} is satisfied in particular if $\phi\in C^1(\real)$ with
$$
\frac{2^{q-2}(q-1)}{\Lambda}|r|^{q-2}\leq\phi'(r) \quad\mbox{for every }r\in\real.
$$
If moreover,
$$
\phi'(r)\leq(q-1)\Lambda|r|^{q-2} \quad\mbox{for every }r\in\real,
$$
then also the upper bound in \eqref{eq:nameless} is satisfied.
\end{remark}

As for the \textit{kernel} $K$, we assume that
\begin{equation}\label{eq:drogasintetica}
K:\real^d\times\real^d\to[0,+\infty]\mbox{ is a Borel-measurable function},
\end{equation}
satisfying the following hypotheses:
\begin{align}
&\mbox{Singularity:}\quad \frac{1}{\Lambda}\chi_{B_\varrho}(x-y)\leq K(x,y)|x-y|^{d+sq}\quad\mbox{for every }x\neq y,\label{eq:Ksingh}\\[2pt]
&\mbox{Integrability:}\quad \iint_{Q(\Omega)}\min\{1,|x-y|^q\}K(x,y)\dd x\dd y<+\infty,\label{eq:Kintegrability}
\end{align}
for some $s\in(0,1)$ and some \textit{interaction radius} $\varrho>0$.

\begin{remark}\label{rmk:KassmannLevy}
    Hypotheses \eqref{eq:drogasintetica}--\eqref{eq:Kintegrability} are reminiscent of those satisfied by (singular) L\'evy-type kernels, which are commonly considered in the literature (for integrodifferential operators defined through L\'evy-type measures, see, e.\,g., \cite{dyda2020regularity}). To be precise, $K:\real^d\times\real^d\to[0,+\infty]$ is a L\'evy-type kernel if it is a measurable function satisfying
    $$
    K(x,y)\geq0,\qquad K(x,y)=K(y,x),\qquad\sup_{x\in\real^d}\int_{\real^d}\min\{1,|x-y|^2\}K(x,y)\dd y<+\infty.
    $$
    In particular, since $\Omega$ is bounded, the above $L^\infty$ condition implies (and is strictly stronger than) the $L^1$ condition \eqref{eq:Kintegrability}, with $q=2$. Regarding the Dirichlet problem for integrodifferential operators whose kernels satisfy similar integrability conditions, see, e.\,g., \cite{felsinger2015dirichlet}. As for the singularity, one usually requires $K$ to lie in the non-integrable side, meaning that
    $$
    K(x,y)\geq \mathcal K(x-y)\geq0,\qquad\mathcal K\not\in L^1(B_\varepsilon)\mbox{ for every }\varepsilon>0.
    $$
    As a compromise between full generality and readability, in this paper we limit ourselves to singularities of the form $\mathcal K(x-y)=|x-y|^{-d-\sigma}$, for some $\sigma>0$, thus imposing condition \eqref{eq:Ksingh}. This kind of singularity is enough to ensure the validity of a Poincar\'e inequality, see, e.\,g., Proposition~\ref{prop:onthesofa}, and the compactness of the embedding of an appropriate functional space in $\cL^2_0$, see Remark~\ref{rmk:Kqcompact}. Actually, we could consider more general singularities, e.\,g., in the spirit of \cite{correa2018nonlocal}.
\end{remark}

To motivate the generality of the above assumptions on the kernels, we provide some examples and observations in Section~\ref{sec:paprika}.

The resulting integrodifferential operators $\mathfrak I$ are of $(s,q)$-Laplacian type. Moreover, they have a variational nature as indeed, if we define the energy functional
\begin{equation}\label{eq:effe}
\mathfrak F(u) \coloneqq \iint_{Q(\Omega)}\Phi\big(u(x)-u(y)\big)K(x,y)\dd x\dd y,
\end{equation}
then, at least formally,
$$
\frac{\dd}{\dd \varepsilon}\Big|_{\varepsilon=0}\mathfrak F(u+\varepsilon v)=
\langle\mathfrak I u,v\rangle,
$$
thus $\mathfrak I$ is the first variation of $\mathfrak F$.
For the rigorous proof see Lemma~\ref{lem:nonhouncomodino}.

One can also define a regional version of these operators by considering, e.g., $K_\Omega$ to be the restriction to $\Omega^2$ of the kernel $K$, and setting
$$
\mathfrak F_\Omega(u) \coloneqq \iint_{\Omega^2}\Phi\big(u(x)-u(y)\big)K_\Omega(x,y)\dd x\dd y,
$$
which gives rise to
$$
\langle \mathfrak I_\Omega u,v \rangle = \iint_{\Omega^2}\phi(u(x)-u(y))(v(x)-v(y))K_\Omega(x,y)\dd x\dd y
$$
as its first variation (see Lemma~\ref{lem:flareon}). Important operators of this kind are the so-called censored stable processes, see, e.\,g., \cite{bogdan2003censored}. We refer the interested reader to \cite{fall2022regional}, both for its exhaustive Introduction, and for very interesting recent results.

More in general, the assumptions on the \textit{regional kernel} are as follows:
\begin{equation}\label{eq:KappaOm}
K_\Omega:\Omega\times\Omega\to[0,+\infty]\mbox{ is a Borel-measurable function}
\end{equation}
satisfying the lower bound
\begin{equation}\label{eq:uhu}
\mbox{Singularity:}\quad \frac{1}{\Lambda}\chi_{B_\varrho}(x-y)\leq K_\Omega(x,y)|x-y|^{d+sq}\quad\mbox{for every }x\neq y,
\end{equation}
for some $s \in (0,1)$ and some interaction radius $\varrho>0$, and
\begin{equation}\label{eq:sticweak}
\mbox{Integrability:}\quad \int_{\Omega'}\int_\Omega |x-y|^q \big(K_\Omega(x,y)+K_\Omega(y,x)\big)\dd x\dd y<+\infty\quad\mbox{for every }\Omega'\subset\subset\Omega,
\end{equation}
which is a weaker, local, version of
\begin{equation}\label{eq:stic}
\iint_{\Omega^2}|x-y|^q K_\Omega(x,y)\dd x\dd y<+\infty.
\end{equation}
The form of \eqref{eq:sticweak} is motivated by the possible non-symmetry of $K_\Omega$. Notice that \eqref{eq:sticweak} is equivalent to
$$
\iint_{\Omega^2\setminus(\Omega\setminus\Omega')^2}|x-y|^q K_\Omega(x,y)\dd x\dd y<+\infty\quad\mbox{for every }\Omega'\subset\subset\Omega.
$$
Compare assumptions \eqref{eq:sticweak}--\eqref{eq:stic} with \eqref{eq:Kintegrability}. The form of the latter ensures also the integrability at infinity of the kernel. Here instead, since $\Omega$ is bounded, we clearly only have
$$
\begin{aligned}
\iint_{\Omega^2\setminus(\Omega\setminus\Omega')^2} \min&\{1,|x-y|^q\} K_\Omega(x,y)\dd x\dd y\leq \iint_{\Omega^2\setminus(\Omega\setminus\Omega')^2} |x-y|^q K_\Omega(x,y)\dd x\dd y\\
    &
    \leq\max\{1,\textup{diam}(\Omega)^q\}\iint_{\Omega^2\setminus(\Omega\setminus\Omega')^2} \min\{1,|x-y|^q\} K_\Omega(x,y)\dd x\dd y.
\end{aligned}
$$
Some examples of such regional kernels can be found in Section~\ref{sec:kernels}.

The precise functional frameworks associated to these two kinds of operators is the object of the upcoming two sections.

\subsubsection{Homogeneous Dirichlet boundary conditions}\label{sec:HDbc}
Let $\phi$ be as above and the kernel $K$ as in \eqref{eq:drogasintetica}--\eqref{eq:Kintegrability}. We point out that, by \eqref{eq:heaviness}, $\mathfrak F(u)$, as defined in \eqref{eq:effe}, is finite if and only if
$$
\|u\|_{K,q} \coloneqq \left(\iint_{Q(\Omega)}|u(x)-u(y)|^q K(x,y)\dd x\dd y\right)^\frac{1}{q}<+\infty,
$$
and we define the space
$$
\cW^{K,q}_0(\Omega) \coloneqq \big\{u:\real^d\to\real\,:\,u=0\mbox{ a.\,e. in }\real^d\setminus\Omega, \ u\in L^q(\Omega)\mbox{ and }\|u\|_{K,q}<+\infty\big\},
$$
which, endowed with the norm $\|\,\cdot\,\|_{L^q(\Omega)}+\|\,\cdot\,\|_{K,q}$, turns out to be a Banach space.
We remark that, since $u=0$ almost everywhere in $\real^d\setminus\Omega$,
\begin{equation}\label{eq:domanis}\begin{aligned}
\|u\|_{K,q}&=\big\|(u(x)-u(y))K(x,y)^\frac{1}{q}\big\|_{L^q(\real^d\times\real^d)},\\
\mathfrak F(u)&=\iint_{\real^{2d}}\Phi(u(x)-u(y))K(x,y)\dd x\dd y.
\end{aligned}\end{equation}
By assumption \eqref{eq:Ksingh} on the kernel we have that
\begin{equation}\label{eq:oppercarita}
\|u\|_{K,q}^q\geq\frac{1}{\Lambda} \iint_{\real^{2d}\cap\{|x-y|<\varrho\}}\frac{|u(x)-u(y)|^q}{|x-y|^{d+sq}}\dd x\dd y,
\end{equation}
for every $u\in\cW_0^{K,q}(\Omega)$.
Furthermore, we have the following fractional Poincar\'e-type inequality.
This is surely well-known to the experts. Nevertheless, since the statement is not completely standard, we provide a proof in Section~\ref{sec:smoothie}.

\begin{proposition}[Fractional Poincar\'e inequality]\label{prop:onthesofa}
Let $s\in(0,1)$, $q\geq1$ and $\varrho>0$, and let $\Omega\subset\real^d$ be a bounded open set. There exists a constant $C(d,s,q,\Omega,\varrho)>0$ such that
    \begin{equation}\label{eq:gen_Poincare}
		\|u\|_{L^q(\Omega)}^q\leq C\iint_{\real^{2d}\cap\{|x-y|<\varrho\}}\frac{|u(x)-u(y)|^q}{|x-y|^{d+sq}}\dd x\dd y,
    \end{equation}
for every measurable function $u:\real^d\to\real$ such that $u=0$ almost everywhere in $\real^d\setminus\Omega$.
\end{proposition}

\begin{remark}
    Actually, the constant in \eqref{eq:gen_Poincare} can be chosen as $C=C(d,s,q,\textup{diam}\,\Omega,\varrho)$. Indeed, let $R \coloneqq \textup{diam}\,\Omega$ and $x_0$ be such that $\Omega \subset B_R(x_0)$. Define also the translated function $\tilde{u}(x) \coloneqq u(x+x_0)$. We apply Proposition~\ref{prop:onthesofa} to $\tilde{u}$ with $\Omega = B_R$, so that the constant $C$ is of the form $C(d,s,q,R,\varrho)$, and obtain
    $$
    \|\tilde{u}\|_{L^q(B_R)}^q\leq C\iint_{\real^{2d}\cap\{|x-y|<\varrho\}}\frac{|\tilde{u}(x)-\tilde{u}(y)|^q}{|x-y|^{d+sq}}\dd x\dd y = C\iint_{\real^{2d}\cap\{|x-y|<\varrho\}}\frac{|u(x)-u(y)|^q}{|x-y|^{d+sq}}\dd x\dd y,
    $$
    where the equality follows from a change of variable. Again by a change of variable and since $u=0$ almost everywhere in $B_R(x_0)\setminus\Omega$, we also have
    $$
    \|u\|_{L^q(\Omega)}^q
    = \|u\|_{L^q(B_R(x_0))}^q
    = \|\tilde{u}\|_{L^q(B_R)}^q.
    $$
\end{remark}

As a consequence of \eqref{eq:gen_Poincare}, we have the following bounds on norms
\begin{align*}
\|u\|_{L^q(\Omega)}^q+\iint_{\real^{2d}\cap\{|x-y|<\varrho\}}&\frac{|u(x)-u(y)|^q}{|x-y|^{d+sq}}\dd x\dd y\leq C\iint_{\real^{2d}\cap\{|x-y|<\varrho\}}\frac{|u(x)-u(y)|^q}{|x-y|^{d+sq}}\dd x\dd y\\
&
\leq C\Lambda\|u\|_{K,q}^q
\leq \|u\|_{L^q(\Omega)}^q+C\Lambda\|u\|_{K,q}^q\leq 2C\Lambda\|u\|_{K,q}^q
\end{align*}
for every measurable function $u:\real^d\to\real$ such that $u=0$ almost everywhere in $\real^d\setminus\Omega$. In particular, $\|\,\cdot\,\|_{K,q}$ is an equivalent norm on the space $\cW^{K,q}_0(\Omega)$.
In the specific case in which the kernel is the usual one, $|x-y|^{-d-sq}$, we denote
$$
\cW^{s,q}_0(\Omega)\coloneqq\cW^{K,q}_0(\Omega)=\big\{u:\real^d\to\real\,:\,u=0\mbox{ a.\,e. in }\real^d\setminus\Omega\mbox{ and }[u]_{W^{s,q}(\real^d)}<+\infty\big\}.
$$
Furthermore we have the equivalence of norms:
\begin{equation}\label{eq:pioggia}\begin{aligned}
[u]_{W^{s,q}(\real^d)}^q&\leq 
\bigg(\frac{2^{q-1}|\Omega|}{\varrho^{d+sq}}+\frac{2\mathscr H^{d-1}(\partial B_1)}{sq\varrho^{sq}}\bigg)\|u\|_{L^q(\Omega)}^q+\iint_{\real^{2d}\cap\{|x-y|<\varrho\}}\frac{|u(x)-u(y)|^q}{|x-y|^{d+sq}}\dd x\dd y\\
&
\leq C[u]_{W^{s,q}(\real^d)}^q,
\end{aligned}\end{equation}
for every $u\in\cW_0^{s,q}(\Omega)$.
As a consequence of \eqref{eq:Ksingh} and the first inequality in \eqref{eq:pioggia}, we have the embedding of $\cW^{K,q}_0(\Omega)$ in $\cW^{s,q}_0(\Omega)$.

\begin{remark}\label{rmk:Kqcompact}
The inequality
$$
\|u\|_{W^{s,q}(\Omega)}\leq C\|u\|_{K,q} \quad\mbox{for every }u\in\cW^{K,q}_0(\Omega)
$$
entails that $\cW^{K,q}_0(\Omega)$ is compactly embedded in $L^q(\Omega)$, hence in $\cL^2_0$, since $q\geq2$. Moreover,
$$
\eqref{eq:Kintegrability}\quad\implies\quad C^{0,1}_c(\Omega)\subset\cW^{K,q}_0(\Omega),
$$
as indeed,
$$
\begin{aligned}
    \|u\|_{K,q}^q &= \iint_{Q(\Omega)\cap\{|x-y|\leq1\}}|u(x)-u(y)|^q K(x,y)\dd x\dd y+\iint_{Q(\Omega)\cap\{|x-y|\geq1\}}|u(x)-u(y)|^q K(x,y)\dd x\dd y\\
    &
    \leq[u]_{C^{0,1}(\real^d)}^q \iint_{Q(\Omega)\cap\{|x-y|\leq1\}}|x-y|^q K(x,y)\dd x\dd y+2^q\|u\|_{L^\infty(\real^d)}^q\iint_{Q(\Omega)\cap\{|x-y|\geq1\}}K(x,y)\dd x\dd y\\
    &
    \leq C \iint_{Q(\Omega)}\min\{1,|x-y|^q\} K(x,y)\dd x\dd y.
\end{aligned}
$$
Hence $\cW^{K,q}_0(\Omega)$ is also densely embedded in $\cL^2_0$.

Actually, \eqref{eq:Kintegrability} is equivalent to requiring Lipschitz functions to have finite $\|\,\cdot\,\|_{K,q}$-norm, as proved in Lemma~\ref{lem:AloeNatale}.
\end{remark}

Furthermore, for every $u\in\cW_0^{K,q}(\Omega)$, we define the operator $\mathfrak I u\in(\cW_0^{K,q}(\Omega))'$ as
\begin{equation}\label{eq:nooooooo}
\langle\mathfrak I u,v\rangle_{K,q} \coloneqq
\iint_{Q(\Omega)}\phi(u(x)-u(y))(v(x)-v(y))K(x,y)\dd x\dd y.
\end{equation}
By \eqref{eq:nameless} and H\"older's inequality, this is indeed well defined for every $u,v\in \cW_0^{K,q}(\Omega)$. Similarly to \eqref{eq:domanis}, since $v=0$ almost everywhere in $\real^d\setminus\Omega$, we can write
$$
\langle\mathfrak I u,v\rangle_{K,q} \coloneqq
\iint_{\real^{2d}}\phi(u(x)-u(y))(v(x)-v(y))K(x,y)\dd x\dd y.
$$
This operator arises as the first variation of the energy $\mathfrak F$ in \eqref{eq:effe}, as shown in the next lemma.

\begin{lemma}\label{lem:nonhouncomodino}
    Let $u \in \cW_0^{K,q}(\Omega)$. Then,
    $$
    \frac{\dd}{\dd \varepsilon}\Big|_{\varepsilon=0}\mathfrak F(u+\varepsilon v)=
    \langle\mathfrak I u,v\rangle_{K,q}
    $$
    for every $v \in \cW_0^{K,q}(\Omega)$.
\end{lemma}

\begin{proof}
    Let $\varepsilon \in (-1,1)$. By the mean value theorem, there exists a function $\tilde{\lambda}_\varepsilon : \real \times \real \to [0,1]$ such that $\Phi(r_1+\varepsilon r_2)-\Phi(r_1) = \varepsilon\phi(r_1+\varepsilon\tilde{\lambda}_\varepsilon(r_1,r_2)r_2)\,r_2$ for every $r_1,r_2\in\real$. Setting
    $$
    \lambda_\varepsilon(x,y) \coloneqq \tilde\lambda_\varepsilon\big(u(x)-u(y),v(x)-v(y)\big) \ \mbox{ for every } x,y \in \real^d,
    $$
    we have
    $$
    \mathfrak F(u+\varepsilon v) - \mathfrak F(u) = \varepsilon \iint_{\real^{2d}} \phi\Big(u(x)-u(y) + \varepsilon\lambda_\varepsilon(x,y)\big(v(x)-v(y)\big)\Big)(v(x)-v(y)) K(x,y) \dd x\dd y.
    $$
    By hypothesis \eqref{eq:nameless} and Young's inequality, and the fact that $|\varepsilon|$ and $\lambda_\varepsilon$ are smaller than $1$, we have
    $$
    \begin{aligned}
    &\phi\Big(u(x)-u(y) + \varepsilon\lambda_\varepsilon(x,y)\big(v(x)-v(y)\big)\Big)(v(x)-v(y)) K(x,y) \\
    &\le \Lambda 2^{q-1}\Big(|u(x)-u(y)|^{q-1}+|v(x)-v(y)|^{q-1}\Big)|v(x)-v(y)|K(x,y) \\
    &= \Lambda 2^{q-1}\Big(|u(x)-u(y)|^{q-1}|v(x)-v(y)|+|v(x)-v(y)|^q\Big)K(x,y) \\
    &\le C\Big(|u(x)-u(y)|^q+|v(x)-v(y)|^q\Big)K(x,y) \eqqcolon \Xi(x,y).
    \end{aligned}
    $$
    Since $u,v \in \cW_0^{K,q}(\Omega)$, the function $\Xi$ belongs to $L^1(\real^{2d})$. As $\lambda_\varepsilon$ and $\varepsilon\to 0$, the claim follows by Lebesgue's dominated convergence theorem.
\end{proof}

\begin{remark}
Two sets of natural hypotheses, commonly found in the literature, that ensure the validity of \eqref{eq:Ksingh}--\eqref{eq:Kintegrability} are the following:
\begin{enumerate}
\item the usual type of kernel that is considered for problems involving the $(s,q)$-Laplacian (see, e.\,g., \cite{kuusi2015nonlocal,di2016local,korvenpaa2017fractional}, but also the seminal \cite{kassmann2009priori} for the case $q=2$) is assumed to satisfy
\begin{align}
&\mbox{Singularity:}\quad \frac{1}{\Lambda}\leq K(x,y)|x-y|^{d+sq}\leq\Lambda\quad\mbox{for }0<|x-y|<\varrho,\label{eq:Ksingular}\\[2pt]
&\mbox{Growth:}\quad 0\leq K(x,y)|x-y|^{d+\kappa}\leq\Lambda\quad\mbox{for }|x-y|\geq\varrho,\label{eq:Kgrowth}
\end{align}
for some $s\in(0,1)$ and $\varrho,\kappa>0$.

Besides the lower bound \eqref{eq:oppercarita}, the stronger assumptions \eqref{eq:Ksingular} and \eqref{eq:Kgrowth} imply also the upper bound
$$
\|u\|_{K,q}^q\leq\Lambda\left(\bigg(\frac{2^{q-1}|\Omega|}{\varrho^{d+\kappa}}+\frac{2\mathscr H^{d-1}(\partial B_1)}{\kappa\varrho^\kappa}\bigg)\|u\|_{L^q(\Omega)}^q+\iint_{\real^{2d}\cap\{|x-y|<\varrho\}}\frac{|u(x)-u(y)|^q}{|x-y|^{d+sq}}\dd x\dd y\right)
$$
for every $u\in\cW^{K,q}_0(\Omega)$. As a consequence, we have the equivalence of the norms
$$
\begin{aligned}
\frac{1}{\Lambda} \iint_{\real^{2d}\cap\{|x-y|<\varrho\}}&\frac{|u(x)-u(y)|^q}{|x-y|^{d+sq}}\dd x\dd y
\leq \|u\|_{K,q}^q
\leq \|u\|_{L^q(\Omega)}^q+\|u\|_{K,q}^q\\
&
\leq C\|u\|_{L^q(\Omega)}^q+\Lambda\iint_{\real^{2d}\cap\{|x-y|<\varrho\}}\frac{|u(x)-u(y)|^q}{|x-y|^{d+sq}}\dd x\dd y\\
&
\leq C \iint_{\real^{2d}\cap\{|x-y|<\varrho\}}\frac{|u(x)-u(y)|^q}{|x-y|^{d+sq}}\dd x\dd y,
\end{aligned}
$$
for every $u\in\cW^{K,q}_0(\Omega)$.

In particular, if the kernel $K$ satisfies the stronger hypotheses \eqref{eq:Ksingular}--\eqref{eq:Kgrowth}, we have a more explicit characterization of $\cW^{K,q}_0(\Omega)$, as indeed $\cW^{K,q}_0(\Omega)=\cW^{s,q}_0(\Omega)$. Hence, by \cite[Theorem~6]{fiscella2015density}, we also know that the space $C^\infty_c(\Omega)$ is dense in $\cW^{K,q}_0(\Omega)$.
\item Another type of kernel, which was considered in \cite{fiscella2015density} and whose properties are closely related to L\'evy measures, can be defined by picking $K : \real^d\to (0,+\infty]$ Borel-measurable such that
\begin{align*}
&\frac{1}{\Lambda} \leq K(x)|x|^{d+sq} \quad \mbox{for every }x \in \real^d \setminus\{0\}, \\[2pt]
&\int_{\real^d} \min\{1,|x|^q\} K(x)\dd x < +\infty,
\end{align*}
for some $s \in (0,1)$. By defining $\tilde K(x,y) \coloneqq K(x-y)$, hypothesis \eqref{eq:Ksingh} is clearly satisfied, whereas for \eqref{eq:Kintegrability} we have
\begin{equation}\label{eq:BananoPasquale}
    \sup_{y \in \real^d} \int_{\real^d} \min\{1,|x-y|^q\} \tilde K(x,y)\dd x < +\infty,
\end{equation}
since
$$
\int_{\real^d} \min\{1,|x-y|^q\} \tilde K(x,y)\dd x = \int_{\real^d} \min\{1,|x|^q\} K(x)\dd x \quad\mbox{for every } y\in\real^d.
$$
Since $\tilde K$ is symmetric in $x$ and $y$, and $\Omega$ is bounded, \eqref{eq:BananoPasquale} implies \eqref{eq:Kintegrability}. As we are precisely in the setting of \cite{fiscella2015density}, we have the density of $C^\infty_c(\Omega)$ in $\cW^{\tilde K,q}_0(\Omega)$.
\end{enumerate}
\end{remark}

\subsubsection{Regional nonlocal operators}\label{sec:chem}
As before, we consider a function $\phi$ satisfying the hypotheses of Section~\ref{sec:sqdef}, and in particular the bound \eqref{eq:nameless}, and a regional kernel $K_\Omega$ as in \eqref{eq:KappaOm}--\eqref{eq:sticweak}.

We define the seminorm
$$
\|u\|_{K_\Omega,q}\coloneqq\left(\iint_{\Omega^2}|u(x)-u(y)|^qK_\Omega(x,y)\dd x\dd y\right)^\frac{1}{q},
$$
and the functional space
$$
\cW^{K_\Omega,q}(\Omega)\coloneqq\left\{u\in L^q(\Omega)\,:\,\|u\|_{K_\Omega,q}<+\infty\right\},
$$
which, endowed with the norm $\|\,\cdot\,\|_{L^q(\Omega)}+\|\,\cdot\,\|_{K_\Omega,q}$, is a Banach space.

By \eqref{eq:uhu} we have the continuous embedding of $W^{K_\Omega,q}(\Omega)$ in $W^{s,q}(\Omega)$, as indeed
\begin{equation}\label{eq:oddish}
    \|u\|_{L^q(\Omega)}^q+[u]_{W^{s,q}(\Omega)}^q\leq \Lambda\|u\|_{K_\Omega,q}^q+\bigg(1+\frac{2^q|\Omega|}{\varrho^{d+sq}}\bigg)\|u\|_{L^q(\Omega)}^q.
\end{equation}
We also observe that
\begin{equation}\label{eq:lapras}\begin{split}
&
\eqref{eq:sticweak}\quad\implies\quad \left\{m+u\,:\, m\in\real, \ u\in C^{0,1}_c(\Omega)\right\}\subset\cW^{K_\Omega,q}(\Omega), \\
&\eqref{eq:stic}\quad\implies\quad C^{0,1}(\overline{\Omega})\subset\cW^{K_\Omega,q}(\Omega).
\end{split}\end{equation}
Indeed, if $u\in C^{0,1}_c(\Omega)$, then $\textup{supp}\,u\subset\Omega'$ for some open set $\Omega'\subset\subset\Omega$, hence
$$
\begin{aligned}
\|u+m\|_{K_\Omega,q}^q&=\iint_{\Omega^2\setminus(\Omega\setminus\Omega')^2} |u(x)-u(y)|^q K_\Omega(x,y)\dd x\dd y\\
    &
    \leq[u]_{C^{0,1}(\Omega)}^q\int_{\Omega'}\int_\Omega |x-y|^q \big(K_\Omega(x,y)+K_\Omega(y,x)\big)\dd x\dd y
\end{aligned}
$$
for every $m\in\real$.
In order to prove Metatheorem~\ref{stm:Metathm2}, we need to consider the closed subsets
$$
\widehat{\cW}^{K_\Omega,q}_m(\Omega)\coloneqq\left\{u\in\cW^{K_\Omega,q}(\Omega)\,:\,\frac{1}{|\Omega|}\int_\Omega u(x)\dd x=m\right\},
$$
with $m\in\real$. For $m=0$, the subset $\widehat{\cW}^{K_\Omega,q}_0(\Omega)$ is also a linear subspace, hence it is a Banach space.
From now on we will make use of the notation
$$
	\mathfrak m(u)\coloneqq \frac{1}{|\Omega|}\int_\Omega u(x)\dd x.
$$

The following fractional Poincar\'e inequality is probably well-known to the experts (see, e.\,g., \cite{ponce2004estimate} for similar results). Nevertheless, since it is not completely standard, we provide both an elementary proof and alternative versions in Appendix~\ref{sec:bubbletea}.

\begin{proposition}[Fractional Poincar\'e inequality]\label{stm:lorenzanonce}
Let $s\in(0,1)$, $q\geq 1$ and $\Omega\subset\real^d$ be a bounded and connected open set with Lipschitz boundary. Then, for every $\varrho >0$ there exists a constant $C=C(d,\Omega,s,q,\varrho)>0$ such that
\begin{equation}\label{eq:eandatavia}
\|u-\mathfrak m(u)\|_{L^q(\Omega)}^q\leq C\iint_{\Omega^2\cap\{|x-y|<\varrho\}}\frac{|u(x)-u(y)|^q}{|x-y|^{d+sq}}\dd x\dd y,
\end{equation}
for every measurable $u:\Omega\to\real$.
\end{proposition}

As a consequence, recalling \eqref{eq:uhu}, we have
\begin{equation}\label{eq:eevee}
\|u\|_{L^q(\Omega)}^q\leq 2^{q-1}\Lambda\, C\|u\|_{K_\Omega,q}^q+2^{q-1}|\Omega||m|^q,
\end{equation}
for every $u\in \widehat{\cW}^{K_\Omega,q}_m(\Omega)$. In particular, for $m=0$, the seminorm $\|\,\cdot\,\|_{K_\Omega,q}$ actually defines an equivalent norm on the space $\widehat{\cW}^{K_\Omega,q}_0(\Omega)$ (if $\Omega$ is connected).

\begin{remark}\label{rmk:porridge}
    As a consequence of the above considerations, we have the following dense and compact embeddings
    $$
    \cW^{K_\Omega,q}(\Omega)\subset L^2(\Omega)\quad\mbox{and}\quad\widehat{\cW}^{K_\Omega,q}_0(\Omega)\subset \dot{L}^2(\Omega),
    $$
    where $\dot{L}^2(\Omega)$ is defined in \eqref{eq:L2punto}. More precisely, the compactness follows from \eqref{eq:oddish} and the compact embedding of $W^{s,q}(\Omega)$ in $L^2(\Omega)$, while the density is ensured by \eqref{eq:lapras}, which is a consequence of \eqref{eq:sticweak}.
\end{remark}

If we define the energy functional $\mathfrak F_\Omega:\cW^{K_\Omega,q}(\Omega)\to[0,+\infty)$ as
\begin{equation}\label{eq:asti}
\mathfrak F_\Omega(u)\coloneqq\iint_{\Omega^2}\Phi(u(x)-u(y))K_\Omega(x,y)\dd x\dd y,
\end{equation}
then the operator $\mathfrak I_\Omega u\in(\cW^{K_\Omega,q}(\Omega))'$, which was defined in \eqref{eq:IOintro} in the variational form
\begin{equation}\label{eq:pruriginoso}
\langle \mathfrak I_\Omega u,v \rangle_{K_\Omega,q} = \iint_{\Omega^2}\phi(u(x)-u(y))(v(x)-v(y))K_\Omega(x,y)\dd x\dd y,
\end{equation}
corresponds to the first variation of $\mathfrak F_\Omega$.

\begin{lemma}\label{lem:flareon}
Let $u\in \cW^{K_\Omega,q}(\Omega)$. Then,
$$
    \frac{\dd}{\dd \varepsilon}\Big|_{\varepsilon=0}\mathfrak F_\Omega(u+\varepsilon v)=
    \langle\mathfrak I_\Omega u,v\rangle_{K_\Omega,q}
$$
for every $v\in \cW^{K_\Omega,q}(\Omega)$.
\end{lemma}

The proof is the same as that of Lemma~\ref{lem:nonhouncomodino}.

\begin{remark}\label{rmk:luca}
    Actually, in order to define $\mathfrak F_\Omega(u)$ it is enough to require $u:\Omega\to\real$ to be a measurable function. Nevertheless,
    $$
    \mathfrak F_\Omega(u)<+\infty\quad\implies\quad\|u\|_{K_\Omega,q}<+\infty\quad\implies\quad u\in L^q(\Omega),
    $$
    hence $u\in\cW^{K_\Omega,q}(\Omega)$. We only need to prove the second implication. For this, we can consider a finite covering of $\overline{\Omega}$ given by balls $B_{\varrho/2}(x_i)$, with $x_1,\dots,x_j\in\overline{\Omega}$. Then
    $$
    [u]_{W^{s,q}(\Omega\cap B_{\varrho/2}(x_i))}^q\leq \Lambda\|u\|_{K_\Omega,q}^q<+\infty,
    $$
    hence $u\in L^q(\Omega\cap B_{\varrho/2}(x_i))$ by \cite[Lemma~D.1.2]{lombardini2018minimization} (whose proof is a simple modification of the argument leading to inequality (8.3) in \cite{di2012hitchhiker}). Thus $u\in L^q(\Omega)$ as claimed.
\end{remark}

\subsubsection{A bridge between the global and regional frameworks}
In this section we show that the functional frameworks introduced in the previous two sections maintain deeper relations than superficially apparent.
Indeed, interestingly, when considering only functions that are equal to zero outside $\Omega$, it is possible to interpret a ``global'' operator (in the sense of \eqref{eq:nooooooo}) as a regional operator (of the type \eqref{eq:pruriginoso}) plus a weighted local term of order zero (see, e.\,g., \cite[formula~(2.3)]{bogdan2003censored} and \cite[formula~(2.2)]{correa2018nonlocal}).

More precisely, let $K:\real^d\times\real^d\to[0,+\infty]$ be Borel-measurable and satisfying \eqref{eq:Ksingh} and \eqref{eq:Kintegrability}, and consider its restriction $K_\Omega\coloneqq K|_{\Omega^2}$.
Clearly, $K_\Omega$ satisfies both \eqref{eq:uhu} and \eqref{eq:stic}.

Moreover, if $u:\real^d\to\real$ is such that $u=0$ almost everywhere in $\real^d\setminus\Omega$, then
$$
\begin{aligned}
\mathfrak F(u)&=\mathfrak F_\Omega(u)+\int_\Omega\Phi(u(x))\Bigg(\int_{\real^d\setminus\Omega}K(x,y)\dd y\Bigg)\dd x
+\int_{\real^d\setminus\Omega}K(x,y)\Bigg(\int_\Omega \Phi(-u(y))\dd y\Bigg)\dd x\\
&	
=\mathfrak F_\Omega(u)+\int_\Omega\Big(\Phi(u(x))\omega_1(x)+\Phi(-u(x))\omega_2(x)\Big)\dd x,
\end{aligned}
$$
where $\mathfrak F$ and $\mathfrak F_\Omega$ are as in \eqref{eq:effe} and \eqref{eq:asti} respectively, and
$$
\omega_1(x)\coloneqq \int_{\real^d\setminus\Omega}K(x,y)\dd y
\quad\mbox{and}\quad
\omega_2(x)\coloneqq \int_{\real^d\setminus\Omega}K(y,x)\dd y
$$
for every $x\in\Omega$.
As a consequence, if $\Phi$ is an even function, we obtain
$$
\mathfrak F(u)=\mathfrak F_\Omega(u)+\int_\Omega\Phi(u(x))\omega(x)\dd x,
$$
where $\omega\coloneqq\omega_1+\omega_2$.

In particular, this is true for the fractional seminorms, yielding
$$
\|u\|_{K,q}^q=\|u\|_{K_\Omega,q}^q+\int_\Omega|u(x)|^q\omega(x)\dd x
=\|u\|_{K_\Omega,q}^q+\|u\|_{L^q(\Omega,\dd\mu_\Omega)}^q,
$$
with $\mu_\Omega\coloneqq \omega\, \mathscr L^d$. This measure is sometimes referred to as the killing measure.

Therefore, given a measurable function $u:\real^d\to\real$ such that $u=0$ almost everywhere in $\real^d\setminus\Omega$, we have
$$
u\in\cW^{K,q}_0(\Omega)\quad\Longleftrightarrow\quad u|_\Omega\in \cW^{K_\Omega,q}(\Omega)\cap L^q(\Omega,\dd\mu_\Omega).
$$

Concerning the weight $\omega$, we point out the following properties.

Clearly, $\omega\geq0$, and \eqref{eq:Kintegrability} ensures that $\omega\in L^1_{\textup{loc}}(\Omega)$.
However, $\omega$ can be highly ``degenerate'', both in terms of integrability and of vanishing.
For example, \eqref{eq:Ksingh} implies that $\omega\not\in L^1(\Omega)$ if $sq>1$. More precisely, if $\partial\Omega$ is of class $C^2$, by arguing as in \cite[Lemma~2.1]{abatangelo2020remark} we find the estimate
$$
\omega(x)\geq c\,\textup{dist}(x,\partial\Omega)^{-sq}\quad\mbox{for every }x\in\Omega\mbox{ s.\,t. }\textup{dist}(x,\partial\Omega)<\varrho/2.
$$
On the other hand, if $K(x,y)=\chi_{B_\varrho}(x-y)|x-y|^{-d-sq}$, then
$$
\omega(x)=0\quad\mbox{for every }x\in\Omega\mbox{ s.\,t. }\textup{dist}(x,\partial\Omega)>\varrho.
$$


\subsection{Convex functions and their subdifferentials}

Let $\Gamma : \real \to (-\infty,+\infty]$ be a proper, convex, and lower semicontinuous function. Then $\gamma \coloneqq \partial\Gamma : \real \to 2^\real$, where
$$
\partial\Gamma(r_0) \coloneqq \Big\{\upsilon \in \real \, : \, \Gamma(r_0) \le \Gamma(r) + \upsilon(r_0-r) \ \ \mbox{for all } r \in \real\Big\}
$$
is the \textit{subdifferential} of $\Gamma$ at $r_0$, is a maximal monotone graph in $\real^2$.

We use the following notations:
$$
\begin{aligned}
	D(\Gamma) &\coloneqq \left\{r \in \real \, : \, \Gamma(r) < +\infty\right\} \ \mbox{for the \textit{effective domain} of } \Gamma,\\
	D(\gamma) = D(\partial\Gamma) &\coloneqq \left\{r \in \real \, : \, \gamma = \partial\Gamma \ne \emptyset\right\}.
\end{aligned}
$$

The basic properties of subdifferentials are summarized for later use in the following lemma. For a proof see e.\,g. \cite[Chapter~1]{barbu2010nonlinear}.

\begin{lemma}\label{lem:subdiff}
	The following properties hold:
	\begin{enumerate}[label=(\roman*),itemsep=5pt]
		\item if $\Gamma$ is differentiable at $r_0 \in \real$, then $\partial\Gamma(r_0) = \Gamma'(r_0)$;
		\item $D(\partial\Gamma)$ is a dense subset of $D(\Gamma)$;
		\item $\textup{Int}\,D(\Gamma) \subset D(\partial\Gamma)$.
	\end{enumerate}
\end{lemma}

Actually, in point (i) also the reverse implication holds true. More precisely:

\begin{lemma}\label{lem:viabrombeis}
    If $r_0\in\textup{Int}\, D(\Gamma)$ and there exists $\upsilon_0\in\real$ such that $\partial\Gamma(r_0)=\{\upsilon_0\}$, then the function $\Gamma$ is differentiable at $r_0$ and $\Gamma'(r_0)=\upsilon_0$.
\end{lemma}

A proof, of which we provide a sketch, can be found, e.\,g., in \cite[Section~9]{ambrosio2000geometric}. We recall that the \textit{first-order sub jet} of a lower semicontinuous function $f:\real\to\real$ at a point $r_0\in\real$ is defined as
$$
J^{1,-}f(r_0)\coloneqq\Big\{\upsilon\in\real\ : \ f(r_0)\leq f(r)+\upsilon(r_0-r)+\textup{o}(|r-r_0|)\Big\}.
$$
This can be interpreted as a local version of the subdifferential for functions which are not necessarily convex.
As mentioned in formula (60) of \cite[Section~9]{ambrosio2000geometric}, if $f$ is convex, then $J^{1,-}f(r_0)=\partial f(r_0)$. To conclude, if $f$ is (semi)convex, then $f$ is differentiable at $r_0$ if and only if $J^{1,-}f(r_0)$ is a singleton (see, e.\,g., \cite[Theorem~15\,(iii)]{ambrosio2000geometric}).

\subsubsection{Moreau-Yosida regularization}\label{sec:MorYos}
The function $\Gamma$ and its subdifferential $\gamma$ can be smoothly approximated resorting to the so-called \textit{Moreau-Yosida regularization}. This procedure is particularly useful when dealing with non-smooth potentials.

\begin{definition}[Moreau-Yosida regularization]
	Let $\lambda\in(0,1)$ be a small regularizing parameter.
	\begin{itemize}[itemsep=5pt]
		\item We approximate $\Gamma$ by its \textit{Moreau regularization}, that is, by the inf-convolution
		$$
		\Gamma_\lambda(r_0) \coloneqq \inf_{r\in\real}\left\{\Gamma(r)+\frac{1}{2\lambda}|r-r_0|^2\right\},
		$$
		for every $r_0\in\real$.
		\item We define the \textit{resolvent} $J_\lambda : \real \to \real$ of $\gamma$ as
		$$
		J_\lambda(r) \coloneqq (I + \lambda\gamma)^{-1}(r)
		$$
		for every $r \in \real$, and the \textit{Yosida approximation} $\gamma_\lambda : \real \to \real$ of $\gamma$ as
		$$
		\gamma_\lambda(r) \coloneqq \frac1\lambda (r - J_\lambda(r))
		$$
		for every $r \in \real$.
	\end{itemize}
\end{definition}
Actually, the inf-convolution does not require the convexity assumption to be defined.

In the next lemma we gather the main properties of the Moreau-Yosida regularization. The proofs can be found in \cite[Chapter~2]{barbu2010nonlinear} and \cite[Chapter~II]{brezis1973operateurs}.

\begin{lemma}\label{lemma:yosida}
	The functions $\Gamma_\lambda$ and $\gamma_\lambda$ are such that the following properties hold:
	\begin{enumerate}[label=(\roman*),itemsep=5pt]
		\item $\gamma_\lambda(r) \in \gamma(J_\lambda(r))$ for every $r\in\real$;
		\item $\gamma_\lambda$ is single-valued, monotone, and Lipschitz continuous with Lipschitz constant $1/\lambda$;
		\item $\Gamma_\lambda$ is convex, continuous, and differentiable with $\Gamma_\lambda'=(\partial\Gamma)_\lambda = \gamma_\lambda$.
	\end{enumerate}
\end{lemma}

\begin{remark}\label{rmk:hatg_growth}
By Lemma~\ref{lemma:yosida}\,(ii) and (iii), and the fact that $\Gamma_\lambda \le \Gamma$ for every $\lambda \in (0,1)$, we obtain that $\Gamma_\lambda$ has quadratic growth, that is,
$$
|\Gamma_\lambda(r)| \le \frac{1}{2\lambda}|r|^2 + |\Gamma_\lambda(0)| \le \frac{1}{2\lambda}|r|^2 + |\Gamma(0)|
$$
for every $r\in\real$.
\end{remark}

We report now two lemmas, whose proof is postponed to Appendix~\ref{sec:MYreg}. The first one gives a coercivity estimate for the regularized potential $\Gamma_\lambda + \Pi$, which will be crucial for proving existence of solutions to the Cahn-Hilliard system.

\begin{lemma}\label{lemma:scoponescientifico}
	Let $\Gamma : \real \to (-\infty,+\infty]$ be a proper lower semicontinuous function, let $\pi : \real \to \real$ be a $C_\pi$-Lipschitz function such that $\pi(0) = 0$ and let $\Pi(r) \coloneqq \int_0^r\pi(z)\dd z$.
 \begin{enumerate}
     \item Assume that
	$$
	\Gamma(r) + \Pi(r) \geq -a_1|r|^p - a_2 \qquad \mbox{for every } r \in \real,
	$$
	for some $a_1,a_2 \geq 0$ and $p \in (0,2)$. Then, there exist $\lambda_0, \alpha, \beta, a_3>0$ depending only on $C_\pi, a_1$ and $a_2$, such that
 \begin{equation}\label{eq:koala}
	\Gamma_\lambda(r) + \Pi(r) \geq -\alpha\lambda^\frac{1}{2} r^2 - a_3|r|^p - \beta,
	\end{equation}
	for every $\lambda \in (0,\lambda_0)$ and for every $r \in \real$.
 \item If, instead,
	$$
	\Gamma(r) + \Pi(r) \geq -a_1 r^2 - a_2 \qquad \mbox{for every } r \in \real,
	$$
	for some $a_1,a_2 \geq 0$, then, there exist $\lambda_0, \alpha, \beta>0$ depending only on $C_\pi, a_1$ and $a_2$, such that
\begin{equation}\label{eq:koalabis}
	\Gamma_\lambda(r) + \Pi(r) \geq  - 2a_1 r^2 -\alpha\lambda^\frac{1}{2} r^2 - \beta,
	\end{equation}
	for every $\lambda \in (0,\lambda_0)$ and for every $r \in \real$.
 \end{enumerate}
\end{lemma} 

Recalling Remark~\ref{rmk:hatg_growth}, if $\Gamma$ is also convex we further point out that
\begin{equation}\label{eq:pidgeon}
\Gamma_\lambda(r)+\Pi(r)\leq
|\Gamma(0)|+\Big(\frac{C_\pi}{2}+\frac{1}{2\lambda}\Big)|r|^2 \quad\mbox{for every }r\in\real,
\end{equation}
for every $\lambda\in(0,1)$.

The second lemma is a $\Gamma$-convergence-type result for $\Gamma_\lambda$, and will be employed to derive the energy estimates.

\begin{lemma}\label{lemma:parzialmenteserena}
	Let $\Gamma:\real\to[0,+\infty]$ be a proper lower semicontinuous function. Then,
	$$
	\Gamma = \Gamma\textrm{-}\liminf_{\lambda\to0} \Gamma_\lambda.
	$$
\end{lemma}


\subsection{Assumptions on the potential}

Throughout the paper, unless otherwise specified, we assume the following:

\begin{hypo}\label{stm:hyp}
The potential $F$ entering the system is represented by the sum $F=\Gamma+\Pi$, and we assume that
\begin{enumerate}[label=(\roman*)]
    \item $\Gamma : \real \to [0,+\infty]$ is a proper, convex, lower semicontinuous function, with $\Gamma(0) = 0$. Then, $\gamma\coloneqq \partial\Gamma : \real \to 2^\real$ is a maximal monotone graph, in the sense of convex analysis, such that $0 \in \gamma(0)$;
    \item $\pi : \real \to \real$ is $C_\pi$-Lipschitz continuous with $\pi(0) = 0$, and we set $\Pi(z) \coloneqq \int_0^z \pi(r) \dd r$. As a consequence,
    \begin{equation}\label{eq:ducks}
    \big|\Pi(r)\big|\leq\frac{C_\pi}{2}|r|^2 \quad\mbox{for every }r\in\real;
    \end{equation}
    \item $F$ satisfies the growth condition
    $$
    F(r)\geq -a_1 |r|^p-a_2 \quad\mbox{for every }r\in\real,
    $$
    for some $p \in (0,2)$ and for some $a_1,a_2\geq0$.
\end{enumerate}
\end{hypo}

\begin{remark}\label{rmk:gen_hyp_potential}
    Actually, Methateorem~\ref{stm:Metathm1} and Metatheorem~\ref{stm:Metathm2} (hence also Metatheorem~\ref{stm:Metathm3}) hold in the case of a more general potential:
    \begin{itemize}
    \item instead of the growth condition in Hypothesis~\ref{stm:hyp}\,(iii), we can assume that
    \begin{equation}\label{eq:tostasenzapesto}
    F(r)\geq -a_1r^2-a_2\qquad\mbox{for every }r\in\real,
    \end{equation}
    for some $a_2,\,a_2\geq0$. This requires some more care when performing the time-discretization step. Indeed, we need to check that the coercivity, namely an analogous of \eqref{eq:gargoyleq} or \eqref{eq:gargoyleregL}, still holds true.
    If $q>2$, then $a_1$ can be arbitrary (hence condition \eqref{eq:tostasenzapesto} is redundant, since $\Pi$ satisfies \eqref{eq:ducks} and $\Gamma\geq0$). A careful inspection of the proof of Lemma~\ref{lemma:neutralquokkaq} (respectively Lemma~\ref{lemma:neutralquokkaregL}) shows that this can be done by suitably using Young's and Poincar\'e's inequalities, estimate \eqref{eq:koalabis}, and the embedding $L^q(\Omega) \hookrightarrow L^2(\Omega)$.
    The same can be done for $q=2$, but in this case we need to consider $a_1$ small enough, e.\,g., $a_1<1/4C$, where $C = C_{s,\varrho,\Omega}$ is the optimal constant of the fractional Poincar\'e inequality. Actually, by being more careful in the proof of Lemma~\ref{lemma:scoponescientifico}\,(2), we can pick $a_1<1/C$. This can be seen by exploiting Young's inequality in the refined form $(a+b)^2 \le (1+\varepsilon)a^2 + (1+1/\varepsilon)b^2$ when proving \eqref{eq:linkpostali}. This is in line with the choice of the energy in \cite[formula~(89)]{akagi2016fractional}.
    Notice that, since $p<2$, Hypothesis~\ref{stm:hyp}\,(iii) ensures the validity of condition \eqref{eq:tostasenzapesto} with $a_1$ arbitrarily small, by Young's inequality. However, Hypothesis~\ref{stm:hyp}\,(iii) can be assumed independently of the integrodifferential operator, while \eqref{eq:tostasenzapesto} calls for a distinction between the cases $q = 2$ and $q>2$, as just observed. Hence, we chose to stick to Hypothesis~\ref{stm:hyp}\,(iii) throughout the paper.
    \item Furthermore, we could allow a dependence also on the variable $x\in\Omega$. More precisely, at essentially no cost, we can consider
    $$
    F(x,r)\coloneqq f(x)\Gamma(r)+\Pi(x,r),
    $$
    with $\Gamma$ as in Hypothesis~\ref{stm:hyp}\,(i), $f\in L^\infty(\Omega)$ such that $C_\divideontimes\geq f\geq C_\divideontimes^{-1}$ almost everywhere in $\Omega$, and $\Pi(x,r)\coloneqq\int_0^r \pi(x,z)\dd z$, where $\pi:\Omega\times\real\to\real$ is a Carath\'eodory function, Lipschitz in $r$, uniformly for almost every $x\in\Omega$, and such that $\pi(x,0)=0$ for almost every $x\in\Omega$. Moreover, we require $F$ to satisfy the growth condition in Hypothesis~\ref{stm:hyp}\,(iii), uniformly for almost every $x\in\Omega$. Considering a potential of this form would not bring any useful insight to the proof and the required changes in the computations are straightforward, hence, for the reader's convenience, we stick to Hypothesis~\ref{stm:hyp} throughout the paper.
    \end{itemize}
\end{remark}

\begin{remark}\label{rmk:barbieref}
    Notice that, under our assumptions, the potential $F$ is such that $F(0)=0$. However, this is not true for notable examples $F_{\rm pol}$, $F_{\rm log}$, and $F_{\rm ob}$ in \eqref{eq:cheetahlo}--\eqref{eq:Fob}. We can accommodate for this requirement by subtracting $F(0)$ and working with the potential $F-F(0)$. Indeed, on the one hand, after such a renormalization all three notable potentials can be written as in Hypothesis~\ref{stm:hyp}. On the other hand,  considering $F$ or $F+c$ for a fixed $c\in\real$ does not impact our analysis, as they both correspond to the same system; the only noticeable change lies in the energy, which changes by a factor $c|\Omega|$. The case of the logarithmic potential requires some more care but it can be considered by extending $F_{\rm log}$ to $\real$ as follows: $F_{\rm log}(\pm1) = \theta\log 2 + \theta_c/2 + c$, $F_{\rm log}(z) = +\infty$ if $|z|>1$.
\end{remark}

\begin{remark}\label{rmk:chao}
    A possible generalization of Hypothesis~\ref{stm:hyp}\,(i)--(ii) could be
    \begin{equation}\label{eq:semiconvexF}
        F:\real\to(-\infty,+\infty]\textrm{ is proper, lower semicontinuous, and \emph{semiconvex}},
    \end{equation}
    where by semiconvexity we mean that there exists $c\geq0$ such that
    \[
    r\mapsto F(r)+\frac{c}{2}|r|^2
    \]
    is a convex function on $\real$ (see, e.\,g., \cite[Definition~10]{ambrosio2000geometric}).
    To be precise, Hypothesis~\ref{stm:hyp}\,(i)--(ii) implies \eqref{eq:semiconvexF}, whereas assumption \eqref{eq:semiconvexF} implies that $F$ can be written as the sum $F=\Gamma+\Pi$, where $\Gamma$ and $\Pi$ satisfy Hypothesis~\ref{stm:hyp}\,(ii) and
    \begin{itemize}
        \item[(i)'] $\Gamma : \real \to [b,+\infty]$ is a proper, convex, lower semicontinuous function, with $\Gamma(a) = b$ for some $a,b \in \real$. Then, $\gamma\coloneqq \partial\Gamma : \real \to 2^\real$ is a maximal monotone graph, in the sense of convex analysis, such that $0 \in \gamma(a)$.
    \end{itemize}
    
    Indeed, suppose that $F=\Gamma+\Pi$, with $\Gamma$ and $\Pi$ satisfying respectively Hypothesis~\ref{stm:hyp}\,(i) and (ii). Since $\pi=\Pi'$ is $C_\pi$-Lipschitz in $\real$, the function $r\mapsto\Pi(r)+C_\pi|r|^2$ is convex, hence $F$ is semiconvex.
    
    On the other hand, let $F$ be as in \eqref{eq:semiconvexF}. By semiconvexity we know that
    \[
\Gamma(r)\coloneqq F(r)+\frac{c}{2}|r|^2
    \]
    is convex, provided $c>0$ is chosen big enough. Then, up to considering a larger $c$, we have that there exist $a,b\in\real$ such that
    \[
\Gamma(r)\geq\Gamma(a)=b\quad\mbox{for every }r\in\real.
    \]
    Moreover, the function $\Pi(r)\coloneqq-\frac{c}{2}|r|^2$ satisfies Hypothesis~\ref{stm:hyp}\,(ii) and, clearly, $F=\Gamma+\Pi$.
    
    To verify that assumption \eqref{eq:semiconvexF} is not equivalent to Hypothesis~\ref{stm:hyp}\,(i)--(ii), one can consider $F(r)=r$.
    
    Assuming \eqref{eq:semiconvexF} (equivalently, assumption (i)' and Hypothesis~\ref{stm:hyp}\,(ii)) would not affect our analysis in the regional case of Section~\ref{sec:pizzallecipolle}, except for making some computations a little more cumbersome.
    However, in the Dirichlet case of Section~\ref{sec:spritz} our argument does not directly carry through under assumption (i)' because of the choice of the test function in \eqref{eq:basketcase}, which requires $0 \in \gamma(0)$. We thus stick to Hypothesis~\ref{stm:hyp}\,(i)--(ii).

    We also mention that semiconvexity (usually referred to as $\lambda$-\textit{convexity}, with $\lambda \in \real$) is standard in the literature regarding gradient flows and the minimizing movements scheme, see, e.\,g., \cite[Section~2.4]{ambrosio2005gradient} and \cite{santambrogio2017euclidean}.
\end{remark}


\section{Homogeneous Dirichlet boundary conditions}\label{sec:spritz}

In this section we consider the problem
\begin{align}
\partial_t u - \Delta w = 0 &\qquad \text{in } \Omega \times (0,T), \label{eq:CH1}\\
w = (-\Delta)^s u + F'(u) &\qquad \text{in } \Omega \times (0,T), \label{eq:CH2} \\
u(x,0) = u_0(x) &\qquad \text{in } \real^d, \label{eq:iniCH} \\
u = 0 &\qquad \text{in } (\real^{d}\setminus\Omega)\times(0,T), \label{eq:bouCH1}\\
w = 0 &\qquad \text{on } \partial\Omega\times(0,T). \label{eq:bouCH2}
\end{align}


\subsection{Main results}\label{sec:LapFracLap}

We consider the Hilbert space
$$
\cH^1_0 \coloneqq \left\lbrace v \in H^1(\real^d) \,:\, v = 0 \text{ a.\,e. in } \real^d\setminus\Omega \right\rbrace,
$$
endowed with the standard $H^1$-norm. This can be identified with $H^1_0(\Omega)$, by extending the functions by $0$ outside of $\Omega$ (see, e.\,g., \cite[Theorems~4.4 and 4.5]{heinonen2006nonlinear}). In particular \eqref{eq:bouCH2} is satisfied in the trace sense by every function $w\in\cH^1_0$. We consider the Laplacian operator $-\Delta:\cH^1_0\to(\cH^1_0)'$ defined in the usual distributional way, via
\begin{equation}\label{eq:nevicava}
\langle-\Delta u,v\rangle_{(\cH^1_0)'}=\int_{\real^d}\nabla u\cdot\nabla v\dd x=(u,v)_{\cH^1_0}.
\end{equation}
We are thus in the setting of Section~\ref{sec:dual}, where the coercivity is provided by Poincar\'e's inequality. Hence, from now on we will consider $\cH^1_0$ to be equipped with the equivalent scalar product \eqref{eq:nevicava} and its induced norm, as customary.

Next, we introduce a weak (energy) formulation of system \eqref{eq:CH1}--\eqref{eq:bouCH2}.
\begin{definition}[Solution to the fractional Cahn-Hilliard system]\label{def:weaksol}
Let $T>0$ be fixed. We say that $(u_s,w_s,\zeta_s)$ is a weak solution to the fractional Cahn-Hilliard system \eqref{eq:CH1}--\eqref{eq:bouCH2} associated with the initial datum $u_0 \in \cH^s_0$ if
$$
\begin{aligned}
u_s &\in L^\infty(0,T;\cH^s_0) \cap W^{1,2}(0,T;(\cH^1_0)'), \\
w_s &\in L^2(0,T;\cH^1_0),\\
\zeta_s &\in L^2(0,T;\cL^2_0), \qquad \zeta_s\in\gamma(u_s) \quad \text{a.\,e. in } \Omega \times (0,T),
\end{aligned}
$$
and $(u_s,w_s,\zeta_s)$ satisfies the following weak formulation of \eqref{eq:CH1}--\eqref{eq:CH2}:
\begin{align}
\partial_t u_s - \Delta w_s = 0 &\qquad \text{in } (\cH^1_0)' \label{eq:CHw1}\\
w_s = (-\Delta)^s u_s + \zeta_s + \pi(u_s) &\qquad \text{in } (\cH^s_0)', \label{eq:CHw2}
\end{align}
almost everywhere in $(0,T)$, with $u_s(0) = u_0$ almost everywhere in $\real^d$.
\end{definition}

\begin{remark}\label{rmk:AubinLions}
We point out that, by the Aubin-Lions-Simon compactness lemma (see \cite[Theorem~5]{simon1986compact}), we have
$$
L^\infty(0,T;\cH^s_0) \cap W^{1,2}(0,T;(\cH^1_0)') \subset C([0,T];\cL^2_0).
$$
\end{remark}

We introduce the energy functional $E_s:\cH^s_0\to(-\infty,+\infty]$ defined as
\begin{equation}\label{eq:deardeer}
E_s(v) \coloneqq \frac{1}{2}\|v\|_{\cH^s_0}^2 + \io F(v(x)) \dd x.
\end{equation}

The first step is to prove the following existence and uniqueness result.

\begin{theorem}[Existence and uniqueness]\label{stm:existence}
Let Hypothesis~\ref{stm:hyp} be satisfied, and let $s \in (0,1)$. Assume
$$
u_0 \in \cH^s_0, \qquad \Gamma(u_0) \in L^1(\real^d).
$$
Then the fractional Cahn-Hilliard system \eqref{eq:CH1}--\eqref{eq:bouCH2} admits a unique solution $(u_s,w_s,\zeta_s)$, in the sense of Definition~\ref{def:weaksol}, and the function $\zeta_s$ is uniquely determined by $(u_s,w_s)$. Moreover, the following energy estimates hold:
\begin{align}
\frac{1}{2}\int_0^t \|w_s(\tau)\|_{\cH^1_0}^2 \dd \tau + E_s(u_s(t)) \le E_s(u_0) &\qquad \mbox{for } \ 0 \le  t < T,\label{eq:en-est}\\[3pt]
E_s(u_s(t)) \le E_s(u_s(\tau)) &\qquad \mbox{for } \ 0 \le \tau \le t < T.\label{eq:est-en}
\end{align}
\end{theorem}

We observe that, since $T$ is an arbitrary fixed positive number, the solution to \eqref{eq:CH1}--\eqref{eq:bouCH2} is unique and satisfies \eqref{eq:est-en}, we can extend the solution up to $T=+\infty$. More precisely, we have the following:

\begin{coro}[Global in time solution]
Let the hypotheses of Theorem~\ref{stm:existence} hold. Then, there exists a unique triple $(u_s,w_s,\zeta_s)$, which satisfies Definition~\ref{def:weaksol} for every $T>0$.
\end{coro}


\begin{remark}
    In the case of the polynomial potential $F_\textrm{pol}$ in \eqref{eq:cheetahlo}, a possible splitting (up to the additive constant $1/4$, see Remark~\ref{rmk:barbieref}) is given by $\Gamma(z)=z^4/4$ and $\Pi(z)=-z^2/2$. Since $\Gamma$ is differentiable, we have $\zeta_s=\gamma(u_s)=u_s^3$ and equation \eqref{eq:CHw2} reads
    \begin{equation}\label{eq:padrozio}
    w_s=(-\Delta)^s u_s+u_s^3-u_s\qquad\mbox{in }(\cH^s_0)'.
    \end{equation}
    A priori, this equation makes sense for test functions $v$ that belong to $L^4(\Omega)\cap\cH^s_0$. This is indeed what we would obtain by directly proving the existence of a solution, without employing the Yosida approximation (but still following the other steps of the argument, within the functional space $L^4(\Omega)\cap\cH^s_0$). On the other hand, obtaining the solution via Theorem~\ref{stm:existence} tells us that $u_s^3\in L^2(0,T;\cL^2_0)$ ensuring that \eqref{eq:padrozio} does indeed make sense. Moreover, the energy estimate \eqref{eq:en-est} implies also that $u_s\in L^\infty(0,T;L^4(\Omega))$. This kind of mismatch in the functional spaces of the test functions for equation \eqref{eq:CHw2} appears also for more general regular potentials, such as those considered in \cite{akagi2016fractional}.
\end{remark}

Then, we come to the behavior of solutions as $s$ tends to $1$.
For this we also need an analogous definition of solution to the classical (local) problem.

\begin{definition}[Solution to the classical Cahn-Hilliard system]\label{def:weaksol_local}
Let $T>0$ be fixed. We say that $(u,w,\zeta)$ is a weak solution to the classical Cahn-Hilliard system \eqref{eq:CH1local}--\eqref{eq:CH2local}, \eqref{eq:iniCH}--\eqref{eq:bouCH2} associated with the initial datum $u_0 \in \cH^1_0$ if
$$
\begin{aligned}
u &\in L^2(0,T;\cH^1_0)\cap W^{1,2}(0,T;(\cH^1_0)'), \\
w &\in L^2(0,T;\cH^1_0), \\
\zeta & \in L^2(0,T;\cL^2_0), \qquad \zeta\in\gamma(u) \quad \text{a.\,e. in } \Omega \times (0,T),
\end{aligned}
$$
and $(u,w,\zeta)$ satisfies the following weak formulation of \eqref{eq:CH1local}--\eqref{eq:CH2local}:
\begin{align}
\partial_t u - \Delta w = 0 &\qquad \text{in } (\cH^1_0)', \label{eq:CHw1_local}\\
w = -\Delta u + \zeta + \pi(u) &\qquad \text{in } (\cH^1_0)', \label{eq:CHw2_local}
\end{align}
almost everywhere in $(0,T)$, with $u(0) = u_0$ almost everywhere in $\Omega$.
\end{definition}

We can now state the asymptotics result.

\begin{theorem}[Asymptotics]\label{stm:asymptotics}
Let Hypothesis~\ref{stm:hyp} be satisfied, and let us consider $u_0 \in \cH^1_0$, with $\Gamma(u_0)\in L^1(\real^d)$.
Let $(u_s,w_s,\zeta_s)$ denote the corresponding sequence of unique solutions to \eqref{eq:CH1}--\eqref{eq:bouCH2}, according to Definition~\ref{def:weaksol}. Then there exists a triple of limit functions $(u,w,\zeta)$, which is a weak solution to the local equation in the sense of Definition~\ref{def:weaksol_local} (with a dimensional constant multiplying the Laplacian in \eqref{eq:CHw2_local}), and
\begin{equation}\label{eq:aiuto}
    \begin{array}{rcll}
	w_s &\rightharpoonup& w &\quad \text{weakly in } L^2(0,T;\cH^1_0), \\[2pt]
	u_s &\to& u &\quad \text{in } C([0,T];(\cH^1_0)'), \\[2pt]
 u_s(t) &\to& u(t) &\quad \text{strongly in } \cL^2_0\text{ for every }t\in(0,T), \\[2pt]
\partial_t u_s &\rightharpoonup& \partial_t u &\quad \text{weakly in } L^2(0,T;(\cH^1_0)'), \\[2pt]
    \zeta_s &\rightharpoonup& \zeta &\quad \text{weakly in } L^2(0,T;\cL^2_0).
\end{array}
\end{equation}
\end{theorem}

\begin{remark}
More in general, we can consider a sequence $\{s_k\}_{k \in \nat} \subset (0,1)$ such that $s_k \to 1$ as $k\to\infty$, a sequence $\{u_{0,k}\} \subset \cH^{s_k}_0$, and a function $u_0\in\cH^1_0,$ satisfying
\begin{equation*}
\sup_{k\in\nat} \left(\|u_{0,k}\|_{\cH^{s_k}_0} + \|\Gamma(u_{0,k})\|_{L^1(\real^d)}\right) < +\infty, \qquad u_{0,k} \to u_0 \text{ strongly in } \cL^2_0.
\end{equation*}
\end{remark}

\begin{remark}\label{rmk:tartare}
The study of the asymptotic behavior of nonlocal functionals has been widely considered in the literature, in various forms. We limit ourselves to mentioning the seminal \cite{bourgain2001another}, for the limit as $s\to1$ of the $W^{s,q}(\Omega)$-seminorm and related energies, which initiated an extensive research. In particular, such energies, having the form
\begin{equation}\label{eq:diegonaska}
\iint_{\Omega^2}\omega\left(\frac{|u(x)-u(y)|}{|x-y|}\right)\rho_\varepsilon(x-y)\dd x\dd y
\end{equation}
with $\omega:[0,+\infty)\to[0,+\infty)$ continuous and $\{\rho_\varepsilon\}\subset L^1(\real^d)$ converging to $\delta_0$ as $\varepsilon\searrow0$, were further investigated in \cite{ponce2004new}. 
In the context of Cahn-Hilliard equations, defined by considering \eqref{eq:diegonaska} in place of the usual (local) Dirichlet integral in the free energy, the limiting behavior of the solutions has been studied in several papers, see, e.\,g., the latest \cite{davoli2023local}
and the references cited therein.

On the other hand, the asymptotic limits of the fractional Laplacian are quite natural (and to be expected, recall its definition as a pseudodifferential operator). As $s\to1$, it converges to the classical Laplacian, see, e.\,g., \cite[Proposition~4.4\,(ii)]{di2012hitchhiker} for an elementary proof. Instead, as $s\to0$, it converges to the identity, see, e.\,g., \cite[Proposition~4.4\,(i)]{di2012hitchhiker}. This second limiting behavior has been considered within the context of Cahn-Hilliard equations in \cite{akagi2016fractional}. Theorem~\ref{stm:asymptotics} can thus be considered, in some sense, complementary to those results.

As for the asymptotics of (the eigenvalues of) $(-\Delta)^s_q$, we refer the interested reader to \cite{brasco2015stability}.
\end{remark}

\subsubsection{Proof of Theorem~\ref{stm:existence}}\label{sec:exince}
This section contains the proof of existence of a unique solution $(u_s,w_s,\zeta_s)$ to the fractional Cahn-Hilliard equation with homogeneous Dirichlet boundary conditions. Since $s$ is fixed in this Section, from now on we drop it as a subscript. We subdivide it in different steps.

\textbf{Step 1: uniqueness.}
Suppose that $(u_1,w_1,\zeta_1)$ and $(u_2,w_2,\zeta_2)$ are weak solutions of \eqref{eq:CH1}--\eqref{eq:bouCH2} in the sense of Definition~\ref{def:weaksol} with the same initial datum $u_0$. Here we are assuming that the potential $F=\Gamma +\Pi$ satisfies Hypothesis~\ref{stm:hyp}\,(i) and (ii), but not necessarily (iii). If we denote $\tilde{u}=u_1-u_2$ and $\tilde{w}=w_1-w_2$, then
\begin{equation}\label{eq:plagiarism}
\partial_t \tilde{u}-\Delta \tilde{w}=0\textrm{ in }(\cH^1_0)',
\qquad
\tilde{w}=(-\Delta)^s\tilde{u}+\zeta_1-\zeta_2+\pi(u_1)-\pi(u_2)\textrm{ in }(\cH^s_0)',
\end{equation}
for every $t\in(0,T)\setminus\Sigma$, with $\Sigma$ of measure zero. Fixing from now on $t\in (0,T)\setminus \Sigma$, applying $(-\Delta)^{-1}:(\cH^1_0)'\to \cH^1_0$ to both sides of the first equation in \eqref{eq:plagiarism}, we obtain
$$
(-\Delta)^{-1}\frac{\dd}{\dd t} \tilde{u}+\tilde{w}=0\textrm{ in }\cH^1_0.
$$
Recalling Remark~\ref{rmk:triple}, evaluating $(\tilde{u},\,\cdot\,)_{\cL^2_0}\in(\cH^1_0)'$ at $0$, exploiting this equation and identity \eqref{eq:sensata}, we get
$$\begin{aligned}
0=\Big(\tilde{u},(-\Delta)^{-1}\frac{\dd}{\dd t} \tilde{u}+\tilde{w}\Big)_{\cL^2_0}&
=
\Big\langle \tilde{u},(-\Delta)^{-1}\frac{\dd}{\dd t}\tilde{u}\Big\rangle_{(\cH^1_0)'}+\langle \tilde{u},\tilde{w}\rangle_{(\cH^1_0)'}\\
&
=
\frac{1}{2} \frac{\dd}{\dd t}\|\tilde{u}\|^2_{(\cH^1_0)'}+
\langle \tilde{u},\tilde{w}\rangle_{(\cH^1_0)'}.
\end{aligned}
$$
Test now the second equation of \eqref{eq:plagiarism} by $\tilde{u}$. Then, by the monotonicity of $\gamma$ and Hypothesis~\ref{stm:hyp}\,(ii), it follows that
\begin{align*}
\langle \tilde{w},\tilde{u}\rangle_{(\cH^s_0)'} &= \|\tilde{u}\|^2_{\cH^s_0}
+ \io (\zeta_1-\zeta_2)\tilde{u} \dd x
+ \io (\pi(u_1)-\pi(u_2))\tilde{u} \dd x\\
&
\geq \|\tilde{u}\|^2_{\cH^s_0}-C_\pi\|\tilde{u}\|_{\cL^2_0}^2.
\end{align*}
Once again, by Remark~\ref{rmk:triple}
\begin{equation}\label{eq:key-id}
\left\langle \tilde{w},\tilde{u} \right\rangle_{(\cH^s_0)'}= \int_{\real^d} \tilde{u}(x)\tilde{w}(x)\dd x = \left\langle \tilde{u},\tilde{w} \right\rangle_{(\cH^1_0)'},
\end{equation}
hence we obtain
$$
\frac{1}{2} \frac{\dd}{\dd t}\|\tilde{u}\|^2_{(\cH^1_0)'}
+\|\tilde{u}\|^2_{\cH^s_0}
\leq C_\pi\|\tilde{u}\|_{\cL^2_0}^2\leq \frac12\|\tilde{u}\|^2_{\cH^s_0}+C\|\tilde{u}\|^2_{(\cH^1_0)'},
$$
for some constant $C\geq 0$. Here we used also Ehrling's lemma (see \cite[Lemma~8]{simon1986compact}), i.\,e., for each $\varepsilon > 0$ there exists a constant $C_\varepsilon \ge 0$ such that
\begin{equation}\label{eq:ehrling}
\|v\|_{\cL^2_0} \le \varepsilon\|v\|_{\cH^s_0} + C_\varepsilon\|v\|_{(\cH^1_0)'} \quad \text{ for all } v \in \cH^s_0.
\end{equation}
Thus we get
$$
\frac{1}{2} \frac{\dd}{\dd t}\|\tilde{u}\|^2_{(\cH^1_0)'}
+\frac{1}{2}\|\tilde{u}\|^2_{\cH^s_0}
\leq C\|\tilde{u}\|^2_{(\cH^1_0)'},
$$
for every $t\in(0,T)\setminus\Sigma$.
Integrating on $(0,t^*)$ and exploiting the fundamental theorem of calculus, we obtain
$$
\|\tilde{u}(t^*)\|^2_{(\cH^1_0)'}\leq \|\tilde{u}(0)\|^2_{(\cH^1_0)'}+2C\int_0^{t^*} \|\tilde{u}(t)\|^2_{(\cH^1_0)'}\dd t,
$$
for every $t^*\in(0,T)$.
Since $\tilde{u}(0)=0$ almost everywhere in $\real^d$, by Gr\"onwall's inequality in integral form, this implies that $\tilde{u}(t^*)=0$ almost everywhere in $\real^d$, for every $t^*\in[0,T).$ In turn, by the first equation in \eqref{eq:plagiarism}, this entails that $\tilde{w}(t^*)=0$ almost everywhere in $\real^d$, for every $t^*\in[0,T)$. Finally, by the second equation in \eqref{eq:plagiarism}, we can conclude that $\zeta_1(t^*)=\zeta_2(t^*)$ almost everywhere in $\real^d$, for every $t^*\in[0,T)$.

\textbf{Notation.} Since $\zeta$ is uniquely determined by $(u,w)$, from now on by solution we will mean the pair $(u,w)$ instead of the triple $(u,w,\zeta)$.

\textbf{Step 2: approximation.}
For every $\lambda \in (0,1)$, let $\Gamma_\lambda : \real \to [0,+\infty)$ be the Moreau regularization of $\Gamma$ as defined in Section~\ref{sec:MorYos}, and let $\gamma_\lambda \coloneqq \Gamma_\lambda'$, which has Lipschitz constant $1/\lambda$. We consider the approximated problem
\begin{align}
	\partial_t u^\lambda - \Delta w^\lambda = 0 &\qquad \text{in } (\cH^1_0)', \label{eq:CH1_lambda}\\
	w^\lambda = (-\Delta)^s u^\lambda + \gamma_\lambda(u^\lambda) + \pi(u^\lambda) &\qquad \text{in } (\cH^s_0)', \label{eq:CH2_lambda}
\end{align}
almost everywhere in $(0,T)$, and
\begin{equation}
	u^\lambda(x,0) = u_0(x) \qquad \text{in } \Omega. \label{eq:iniCH_lambda}
\end{equation}

We need to show that the approximating system \eqref{eq:CH1_lambda}--\eqref{eq:iniCH_lambda} admits a solution for every $\lambda > 0$ fixed. Inspired by \cite[Section~4]{akagi2016fractional}, we proceed by time-discretization, with appropriate adjustments to deal with the Yosida approximation.

\textbf{Step 3: time-discretization.}
Let $N \in \nat$ and let $\tau = \tau/N \coloneqq T/N$ be the time step. For each $\lambda > 0$ fixed, we carry out the following discretization of \eqref{eq:CH1_lambda}--\eqref{eq:iniCH_lambda}:
\begin{align}
	\frac{u^\lambda_n - u^\lambda_{n-1}}{\tau} - \Delta w^\lambda_n = 0 &\qquad \text{in } (\cH^1_0)', \label{eq:CH1_n}\\
	w^\lambda_n = (-\Delta)^s u^\lambda_n + \gamma_\lambda(u^\lambda_n) + \pi(u^\lambda_n) &\qquad \text{in } (\cH^s_0)', \label{eq:CH2_n}
\end{align}
for $n = 1, ..., N$ with initial condition $u^\lambda_0 = u_0$. In order to show existence of $(u^\lambda_n,w^\lambda_n)$ satisfying \eqref{eq:CH1_n}--\eqref{eq:CH2_n}, we introduce the functional $E^\lambda_n : \cH^s_0 \to \real$ given by
$$
E^\lambda_n(u) \coloneqq \frac{\tau}{2} \left\|\frac{u - u^\lambda_{n-1}}{\tau}\right\|_{(\cH^1_0)'}^2 + \frac{1}{2}\|u\|_{\cH^s_0}^2 + \io \left(\Gamma_\lambda(u) + \Pi(u)\right) \dd x,
$$
and we consider a sequence of minimizers constructed by iteration.
 In order to do this, we first define the energy
$$
E^\lambda(u)\coloneqq \frac{1}{2}\|u\|_{\cH^s_0}^2 + \io \left(\Gamma_\lambda(u) + \Pi(u)\right) \dd x,
$$
for every $u\in\cH^s_0$, and we point out the following two preliminary results.

\begin{lemma}\label{lemma:neutralquokka}
The functional $E^\lambda$ is lower semicontinuous with respect to the $L^2$-convergence. Moreover
\begin{equation}\label{eq:gargoyle}
E^\lambda(u)\geq \frac{1}{4}\|u\|_{\cH^s_0}^2-\beta_0|\Omega| \qquad\mbox{for every }u\in\cH^s_0,
\end{equation}
for $\lambda\in(0,\bar\lambda)$, 
where $\beta_0=\beta+8C_{\textrm{P}} c(p)a_3^2$ and $\bar\lambda=\bar\lambda(C_{\textrm{P}},a_1,a_2,C_\pi)>0$ is small enough.
\end{lemma}

\begin{proof}
In order to prove the lower semicontinuity, we consider $v_k\in\cL^2_0$ such that $v_k\to v\in \cL^2_0$ strongly in $L^2(\Omega)$ and we observe that by the $L^2$-convergence and \eqref{eq:pidgeon}, 
$$
\io \big(\Gamma_\lambda(v_k)+\Pi(v_k)\big)\dd x\ \xrightarrow{k\to\infty}\ 
\io \big(\Gamma_\lambda(v)+\Pi(v)\big)\dd x.
$$
On the other hand, by the pointwise convergence and Fatou's lemma, we obtain
\begin{equation}\label{eq:pandasaviour}
\|v\|_{\cH^s_0}^2\leq\liminf_{k\to\infty}\|v_k\|_{\cH^s_0}^2.
\end{equation}
Indeed, suppose by contradiction that \eqref{eq:pandasaviour} does not hold true. If $v\in\cH^s_0$, then, we can find $\varepsilon>0$ and a subsequence $v_{k_h}$ such that
\begin{equation}\label{eq:platipus}
\|v_{k_h}\|_{\cH^s_0}^2\leq\|v\|_{\cH^s_0}^2-\varepsilon,
\end{equation}
for every $k_h$. The strong $L^2$-convergence implies that there exist a subsequence $\{v_{k_{h_j}}\}$ such that $v_{k_{h_j}}\to v$ almost everywhere in $\real^d$, as $j\to\infty$. However, by Fatou's lemma,
$$
\|v\|_{\cH^s_0}^2\leq\liminf_{j\to\infty}\|v_{k_{h_j}}\|_{\cH^s_0}^2,
$$
which gives a contradiction with \eqref{eq:platipus}. If instead $v\not\in\cH^s_0$, then we can repeat the same argument, by substituting the right-hand side of \eqref{eq:platipus} with $1/\varepsilon$, thus proving \eqref{eq:pandasaviour}.

As for \eqref{eq:gargoyle}, by \eqref{eq:koala} we have
$$
E^\lambda(u)\geq \frac{1}{2}\|u\|_{\cH^s_0}^2-\alpha\lambda^\frac{1}{2}\|u\|_{\cL^2_0}^2-\beta|\Omega|-a_3\|u\|_{L^p(\Omega)}^p,
$$
and we recall the fractional Poincar\'e inequality (see, e.\,g., \cite[Proposition~2.5]{brasco2015stability} together with the density of $C^\infty_c(\Omega)$ in $\cH^s_0$ proved in \cite[Theorem~6]{fiscella2015density})
$$
\|u\|_{\cL^2_0}^2\leq C_{\textrm{P}}\|u\|_{\cH^s_0}^2,
$$
with $C_{\textrm{P}}=C_{\textrm{P}}(d,2,\Omega)>0$, and 
for every $u\in\cH^s_0$. Moreover, exploiting also Young's inequality,
\begin{align*}
a_3\|u\|_{L^p(\Omega)}^p&\leq \varepsilon\|u\|_{\cL^2_0}^2+\frac{c(p)a_3^2}{\varepsilon}|\Omega|\\
&
\leq \varepsilon C_{\textrm{P}}\|u\|_{\cH^s_0}^2 +\frac{c(p)a_3^2}{\varepsilon}|\Omega|\\
&
\leq\frac{1}{8}\|u\|_{\cH^s_0}^2+8C_{\textrm{P}} c(p)a_3^2|\Omega|.
\end{align*}
Therefore, for every $\lambda>0$ sufficiently small, \eqref{eq:gargoyle} holds true.
\end{proof}

\begin{lemma}\label{lemma:emma}
Given $\lambda\in(0,\bar\lambda)$ and any $g\in\cL^2_0$ we define the energy $E_g^\lambda:\cH^s_0\to\real$ by setting
$$
E_g^\lambda(u) \coloneqq \frac{1}{2\tau}\|u-g\|^2_{(\cH^1_0)'} + E^\lambda(u).
$$
Then, there exists at least one function $u_\star\in\cH^s_0$ such that
$$
E_g^\lambda(u_\star)=\inf_{u\in\cH^s_0}E_g^\lambda(u).
$$
\end{lemma}

\begin{proof}
First of all we observe that $E_g^\lambda(u)\in\real$ for every $u\in\cH^s_0$ by \eqref{eq:pidgeon}. Moreover, by \eqref{eq:gargoyle} the infimum is finite. Then we can consider a minimizing sequence, i.\,e., $\{v_k\}\subset\cH^s_0$ such that
$$
\lim_{k\to\infty}E^\lambda_g(v_k) = \inf_{u\in\cH^s_0}E_g^\lambda(u),
$$
and, by \eqref{eq:gargoyle}, we have
$$
\|v_k\|_{\cH^s_0}^2\leq4\beta_0|\Omega|+4\inf_{u\in\cH^s_0}E_g^\lambda(u)+4
$$
for every $k$ big enough. Together with the fractional Poincar\'e inequality, this implies that $\{v_k\}$ is bounded in $H^s(\Omega)$, which is compactly embedded in $L^2(\Omega)$ (see, e.\,g., \cite[Theorem~7.1]{di2012hitchhiker}), hence there exists $u_\star\in\cH^s_0$ such that
$$
v_k\to u_\star\qquad\mbox{strongly in }L^2(\Omega),
$$
up to a subsequence, that we do not relabel. By the $L^2$-convergence and the lower semicontinuity of $E^\lambda$ proved in Lemma~\ref{lemma:neutralquokka},
$$
E_g^\lambda(u_\star)\leq\liminf_{k\to\infty}E^\lambda_g(v_k)=
\inf_{u\in\cH^s_0}E_g^\lambda(u),
$$
concluding the proof.
\end{proof}

We now proceed to the construction of the solution $(u_n^\lambda,w_n^\lambda)$ to \eqref{eq:CH1_n}--\eqref{eq:CH2_n}.

As a first step, we apply Lemma~\ref{lemma:emma} to pick a minimizer $u_1^\lambda\in\cH^s_0$ of $E_1^\lambda=E_{u_0}^\lambda$.
Then, $u_n^\lambda\in\cH^s_0$ is chosen iteratively as a minimizer of $E_n^\lambda=E_{u_{n-1}}^\lambda$.
Let us recall that, according to Section~\ref{sec:dual}, $(-\Delta)^{-1} : (\cH^1_0)'\to \cH^1_0$ is the inverse duality map, whence $(-\Delta)^{-1}$ coincides with the Fréchet derivative of the functional
$$
v \in (\cH^1_0)'\mapsto \frac12 \|v\|_{(\cH^1_0)'}^2.
$$
Thus, since $E^\lambda_n$ is of class (at least) $C^1$ in $\cH^s_0$, and recalling \eqref{eq:key-id}, we obtain
$$
(-\Delta)^{-1}\left(\frac{u^\lambda_n - u^\lambda_{n-1}}{\tau}\right) + (-\Delta)^s u^\lambda_n + \gamma_\lambda(u^\lambda_n) + \pi(u^\lambda_n) = 0 \quad \text{ in } (\cH^s_0)'.
$$
Setting
\begin{equation}\label{eq:wnl}
w^\lambda_n \coloneqq \Delta^{-1}\left(\frac{u^\lambda_n - u^\lambda_{n-1}}{\tau}\right) \in \cH^1_0,
\end{equation}
we have constructed a solution $(u^\lambda_n,w^\lambda_n)$ to \eqref{eq:CH1_n}--\eqref{eq:CH2_n}.


\textbf{Step 4: estimates independent of $\tau$.}
We observe that
\begin{equation}\label{eq:mail}
E^\lambda(u)\leq E^\lambda_n(u)\qquad\mbox{and}\qquad
E^\lambda_n(u^\lambda_{n-1})=E^\lambda(u^\lambda_{n-1}),
\end{equation}
for every $n\geq1$. Since $\Gamma_\lambda\leq\Gamma$, we also have
\begin{equation}\label{eq:ail}
E^\lambda(u)\leq E_s(u),
\end{equation}
with $E_s$ as defined in \eqref{eq:deardeer}.
Moreover, by the definition of $u^\lambda_n$ as minimizer,
\begin{equation}\label{eq:mai}
E^\lambda_n(u^\lambda_n)\leq E^\lambda_n(u^\lambda_{n-1}).
\end{equation}
By using \eqref{eq:mail} and \eqref{eq:mai}, we find
$$
E^\lambda(u^\lambda_n)\leq E^\lambda_n(u^\lambda_n)\leq E^\lambda_n(u^\lambda_{n-1})=E^\lambda(u^\lambda_{n-1}),
$$
and hence
$$
\frac{\tau}{2} \left\|\frac{u_n^\lambda - u^\lambda_{n-1}}{\tau}\right\|_{(\cH^1_0)'}^2+E^\lambda(u^\lambda_n)-E^\lambda(u^\lambda_{n-1})\leq0,
$$
for every $n\geq1$.
By summing up, and exploiting \eqref{eq:ail}, we obtain
\begin{equation}\label{eq:underclasshero}
\frac{\tau}{2}\sum_{n=1}^{\bar{n}} \left\|\frac{u_n^\lambda - u^\lambda_{n-1}}{\tau}\right\|_{(\cH^1_0)'}^2 + E^\lambda(u_{\bar{n}}^\lambda)\le E_s(u_0).
\end{equation}
Recalling the definition of $w_n^\lambda$ given in \eqref{eq:wnl}, we observe that, by \eqref{eq:lorenzace}
\begin{equation}\label{eq:itendini}
\|w_n^\lambda\|_{\cH^1_0}^2= \left\|\frac{u_n^\lambda - u^\lambda_{n-1}}{\tau}\right\|_{(\cH^1_0)'}^2.
\end{equation}
Therefore
\begin{equation}\label{eq:perdopo}
\frac{\tau}{2}\sum_{n=1}^{\bar{n}} \|w_n^\lambda\|_{\cH^1_0}^2 + E^\lambda(u_{\bar{n}}^\lambda)\le E_s(u_0),
\end{equation}
for every $\bar{n}\geq1$. Recalling estimate \eqref{eq:gargoyle}, we deduce
\begin{equation}\label{eq:esti2_n}
2\tau\sum_{n=1}^{N} \|w^\lambda_n\|_{\cH^1_0}^2 + \max_n \|u^\lambda_n\|_{\cH^s_0}^2\le 8\big( E_s(u_0)+\beta_0|\Omega|\big),
\end{equation}
and also, by \eqref{eq:itendini},
\begin{equation}\label{eq:esti3_n}
\tau\sum_{n=1}^{N} \left\|\frac{u^\lambda_n - u^\lambda_{n-1}}{\tau}\right\|_{(\cH^1_0)'}^2 \le 4\big( E_s(u_0)+\beta_0|\Omega|\big).
\end{equation}
We now construct piecewise constant interpolants of $\{u_n\}$ and $\{w_n\}$, and the piecewise linear interpolant of $\{u_n\}$. More precisely, we define $t_n \coloneqq n\tau=nT/N$ for $n=0,\dots,N$, and
\begin{equation}\label{eq:inter}
\bar{u}^\lambda_\tau(t) \equiv u^\lambda_n, \quad \hat{u}^\lambda_\tau(t) = \frac{t-t_{n-1}}{\tau} u^\lambda_n + \frac{t_n-t}{\tau} u^\lambda_{n-1} \qquad \text{for } t \in [t_{n-1},t_n),
\end{equation}
and $\bar{w}^\lambda_\tau$ analogously. 
Note that
\begin{equation}\label{eq:bar_hat}
\begin{aligned}
\|\bar{u}^\lambda_\tau(t) - \hat{u}^\lambda_\tau(t)\|_{(\cH^1_0)'} &= \left\|u^\lambda_n - \left(\frac{t-t_{n-1}}{\tau}u^\lambda_n + \frac{t_n-t}{\tau}u_{n-1}\right)\right\|_{(\cH^1_0)'} \\
&= \frac{t_n-t}{\tau}\|u^\lambda_n - u^\lambda_{n-1}\|_{(\cH^1_0)'} \stackrel{\eqref{eq:esti3_n}}{\le} 
C\sqrt{\tau} \quad \text{for all } t \in [t_{n-1},t_n).
\end{aligned}\end{equation}
These functions satisfy
\begin{align}
	\partial_t\hat{u}^\lambda_\tau - \Delta \bar{w}^\lambda_\tau = 0 &\qquad \text{in } (\cH^1_0)', \label{eq:CHbar1}\\
	\bar{w}^\lambda_\tau = (-\Delta)^s \bar{u}^\lambda_\tau + \gamma_\lambda(\bar{u}^\lambda_\tau) + \pi(\bar{u}^\lambda_\tau) &\qquad \text{in } (\cH^s_0)',\label{eq:CHbar2}
\end{align}
almost everywhere in $(0,T)$,
and by \eqref{eq:esti2_n}--\eqref{eq:esti3_n}
\begin{equation}\label{eq:esti1_tau}
\int_0^T \|\bar{w}^\lambda_\tau(t)\|_{\cH^1_0}^2 \dd t + \sup_{t \in [0,T]} \|\bar{u}^\lambda_\tau(t)\|_{\cH^s_0}^2
\le C, \qquad \int_0^T \|\partial_t\hat{u}^\lambda_\tau(t)\|_{(\cH^1_0)'}^2 \dd t \le C,
\end{equation}
which entails
\begin{equation}\label{eq:esti2_tau}
\sup_{t \in [0,T]} \|\hat{u}^\lambda_\tau(t)\|_{\cH^s_0}^2 \le C,
\end{equation}
where $C$ is a constant depending on the initial datum $u_0$ but not on $\tau$ nor $\lambda$.

It is also convenient to point out that from \eqref{eq:perdopo} we obtain
\begin{equation}\begin{aligned}\label{eq:perdono}
\frac{1}{2}\int_0^t \|\bar{w}^\lambda_\tau(r)\|_{\cH^1_0}^2\dd r
+E^\lambda(\bar{u}^\lambda_\tau(t))
\leq E_s(u_0),
\end{aligned}\end{equation}
for every $t\in(0,T)$.


\textbf{Step 5: convergence as $\tau \to 0$.}
From \eqref{eq:bar_hat} and the estimates \eqref{eq:esti1_tau} and \eqref{eq:esti2_tau} established above, there exists a (non-relabeled) subsequence of $\tau \to 0$ (equivalently, $N\to\infty$) such that
\begin{equation}\label{eq:conv}
\begin{array}{rcll}
	\bar{w}^\lambda_\tau &\rightharpoonup& w^\lambda &\quad \text{weakly in } L^2(0,T;\cH^1_0), \\[2pt]
	\bar{u}^\lambda_\tau &\stackrel{*}{\rightharpoonup}& u^\lambda &\quad \text{weakly-star in } L^\infty(0,T;\cH^s_0), \\[2pt]
 \hat{u}^\lambda_\tau &\stackrel{*}{\rightharpoonup}& u^\lambda &\quad \text{weakly-star in } L^\infty(0,T;\cH^s_0), \\[2pt]
	\partial_t\hat{u}^\lambda_\tau &\rightharpoonup& \partial_t u^\lambda &\quad \text{weakly in } L^2(0,T;(\cH^1_0)'),
\end{array}
\end{equation}
for some
$$
u^\lambda \in L^\infty(0,T;\cH^s_0)\cap W^{1,2}(0,T;(\cH^1_0)'), \quad w^\lambda \in L^2(0,T;\cH^1_0).
$$
As a consequence of the Aubin-Lions-Simon compactness lemma and the compact embeddings $\cH^s_0 \hookrightarrow \cL^2_0 \hookrightarrow (\cH^1_0)'$, we obtain the stronger convergence
\begin{equation}\label{eq:conv_ALS}
\hat{u}^\lambda_\tau \to u^\lambda \quad \text{strongly in } C([0,T];\cL^2_0).
\end{equation}
Estimate \eqref{eq:bar_hat}, along with \eqref{eq:conv_ALS}, yield
$$
\sup_{t \in [0,T]} \|\bar{u}^\lambda_\tau(t) - u^\lambda(t)\|_{(\cH^1_0)'} \le \sup_{t \in [0,T]} \|\bar{u}^\lambda_\tau(t) - \hat{u}^\lambda_\tau(t)\|_{(\cH^1_0)'} + \sup_{t \in [0,T]} \|\hat{u}^\lambda_\tau(t) - u^\lambda(t)\|_{(\cH^1_0)'} \to 0.
$$
Moreover, from Ehrling's lemma \eqref{eq:ehrling} and \eqref{eq:esti1_tau}, for any $\varepsilon > 0$ we can take $C_\varepsilon \ge 0$ such that
\begin{align*}
\|\bar{u}^\lambda_\tau(t) - u^\lambda(t)\|_{\cL^2_0} &\le \varepsilon \|\bar{u}^\lambda_\tau(t) - u^\lambda(t)\|_{\cH^s_0} + C\varepsilon \|\bar{u}^\lambda_\tau(t) - u^\lambda(t)\|_{(\cH^1_0)'} \\
&\le \varepsilon C + C\varepsilon \sup_{r \in [0,T]}\|\bar{u}^\lambda_\tau(r) - u^\lambda(r)\|_{(\cH^1_0)'} \quad \text{for any } t \in [0,T].
\end{align*}
There follows that
$$
\bar{u}^\lambda_\tau \to u^\lambda \quad \text{strongly in } L^\infty(0,T;\cL^2_0),
$$
hence the Lipschitz continuity of $\gamma_\lambda$ (see Lemma~\ref{lemma:yosida}\,(ii)) and $\pi$ yields
$$
\gamma_\lambda(\bar{u}^\lambda_\tau) \to \gamma_\lambda(u^\lambda), \ \pi(\bar{u}^\lambda_\tau) \to \pi(u^\lambda) \quad \text{strongly in } L^\infty(0,T;\cL^2_0).
$$
By the second convergence in \eqref{eq:conv}, we also obtain
$$
\int_0^T \vartheta(t)\langle(-\Delta)^s\bar{u}^\lambda_\tau(t),v\rangle_{(\cH^s_0)'}\dd t\xrightarrow{\tau\to0} 
\int_0^T \vartheta(t)\langle(-\Delta)^s u^\lambda(t),v\rangle_{(\cH^s_0)'}\dd t,
$$
for every $\vartheta\in L^1(0,T)$ and every $v\in\cH^s_0$.
Therefore, since by \eqref{eq:CHbar1}--\eqref{eq:CHbar2} we have
\begin{equation*}\begin{aligned}
&\int_0^T \left(\io \partial_t \hat{u}^\lambda_\tau(x,t)v_1(x)+\nabla \bar{w}_\tau^\lambda(x,t)\cdot\nabla v_1(x) \dd x\right)\vartheta(t) \dd t = 0,\\
&
\int_0^T \left(\langle(-\Delta)^s \bar{u}_\tau^\lambda(t),v_2\rangle_{(\cH_0^s)'} + \io \big(\gamma_\lambda(\bar{u}_\tau^\lambda(x,t))+\pi(\bar{u}_\tau^\lambda(x,t))-\bar{w}_\tau^\lambda(x,t)\big)v_2(x) \dd x\right)\vartheta(t) \dd t=0,
\end{aligned}\end{equation*}
for every $\vartheta\in L^2(0,T)$, $v_1\in\cH^1_0$ and $v_2\in\cH^s_0$, passing to the limit as $\tau\to0$, we obtain
\begin{equation}\label{eq:CHlambdaint}\begin{aligned}
&\int_0^T \left(\io \partial_t u^\lambda(x,t)v_1(x)+\nabla w^\lambda(x,t)\cdot\nabla v_1(x) \dd x\right)\vartheta(t) \dd t = 0,\\
&
\int_0^T \left(\langle(-\Delta)^s u^\lambda(t),v_2\rangle_{(\cH_0^s)'} + \io \big(\gamma_\lambda(u^\lambda(x,t))+\pi(u^\lambda(x,t))-w^\lambda(x,t)\big)v_2(x) \dd x\right)\vartheta(t) \dd t=0.
\end{aligned}\end{equation}
Thus, by the fundamental lemma of the calculus of variations for Bochner spaces, $(u^\lambda,w^\lambda)$ is a weak solution to \eqref{eq:CH1_lambda}--\eqref{eq:iniCH_lambda}. We observe that, by applying the results that we proved in Step~1 of Section~\ref{sec:exince}, the solution is also unique.

More precisely, we have proved that there exists a subsequence $\tau_{N_h}\to0$ for which the convergences in \eqref{eq:conv} hold true for a pair of limit functions $(u^\lambda,w^\lambda)$ that, a priori, could depend on the subsequence $\{\tau_{N_h}\}_{h\in\nat}$. Actually, by Step~1 of Section~\ref{sec:exince}, applied here with $\Gamma$ replaced by $\Gamma_\lambda$, which is still convex by Lemma~\ref{lemma:yosida}\,(iii), the pair $(u^\lambda,w^\lambda)$ is the unique weak solution to \eqref{eq:CH1_lambda}--\eqref{eq:iniCH_lambda}. Thus any such subsequence converges to the same limit, so that the convergence does indeed hold true as $\tau_N=T/N\to0$.

\textbf{Step 6: estimates independent of $\lambda$.}
First of all, note that estimates \eqref{eq:esti1_tau}--\eqref{eq:esti2_tau} are preserved in the limit as $\tau \to 0$. We therefore have
\begin{align}\label{eq:lest}
\int_0^T \|w^\lambda(t)\|_{\cH^1_0}^2 \dd t \le C, \quad
\int_0^T \|\partial_t u^\lambda(t)\|_{(\cH^1_0)'}^2 \dd t \le C, \quad
\sup_{t \in [0,T]} \|u^\lambda(t)\|_{\cH^s_0}^2 \le C,
\end{align}
where $C$ is a constant independent of $\lambda$.

Furthermore, since $\gamma_\lambda(0)=0$ and $\gamma_\lambda$ is Lipschitz continuous, we have $\gamma_\lambda(u_\lambda)\in\cH^s_0$, hence we can use it as a test function in equation \eqref{eq:CH2_lambda}. Integrating on $(0,T)$ we obtain
\begin{equation}\label{eq:basketcase}\begin{aligned}
\int_0^T\io &|\gamma_\lambda(u^\lambda(x,t))|^2 \dd x\dd t\\
&
+ \int_0^T(1-s)\iint_{\real^{2d}} \frac{u^\lambda(x,t)-u^\lambda(y,t)}{|x-y|^{d+2s}}\big(\gamma_\lambda(u^\lambda(x,t))-\gamma_\lambda(u^\lambda(y,t))\big) \dd x\dd y\dd t\\
&\qquad\quad
= \int_0^T\io \big(w^\lambda(x,t)-\pi(u^\lambda(x,t))\big)\gamma_\lambda(u^\lambda(x,t)) \dd x\dd t.
\end{aligned}\end{equation}
By the monotonicity of $\gamma_\lambda$, the second term on the left-hand side is nonnegative. Thus, applying Young's inequality, recalling that $\pi$ is Lipschitz, and exploiting \eqref{eq:lest}, we obtain the estimate
\begin{equation}\label{eq:lest2}
\|\gamma_\lambda(u^\lambda)\|_{L^2(0,T;\cL^2_0)} \leq C,
\end{equation}
uniformly in $\lambda$.

Finally, by passing to the limit $\tau\to0$ in \eqref{eq:perdono}, recalling the lower semicontinuity of $E_\lambda$ proved in Lemma~\ref{lemma:neutralquokka}, we obtain
\begin{equation}\begin{aligned}\label{eq:scusa}
\frac{1}{2}\int_0^t \|w^\lambda(r)\|_{\cH^1_0}^2\dd r
+E^\lambda(u^\lambda(t))
\leq E_s(u_0),
\end{aligned}\end{equation}
for every $t\in(0,T)$.

\textbf{Step 7: convergence as $\lambda \to 0$.}
Estimates \eqref{eq:lest} and \eqref{eq:lest2} imply that (up to considering a subsequence that we do not relabel)
\begin{equation}\label{eq:kilt}
\begin{array}{rcll}
	w^\lambda &\rightharpoonup& w &\quad \text{weakly in } L^2(0,T;\cH^1_0), \\[2pt]
	u^\lambda &\stackrel{*}{\rightharpoonup}& u &\quad \text{weakly-star in } L^\infty(0,T;\cH^s_0), \\[2pt]
	\partial_t u^\lambda &\rightharpoonup& \partial_t u &\quad \text{weakly in } L^2(0,T;(\cH^1_0)'), \\[2pt]
    \gamma_\lambda(u^\lambda) &\rightharpoonup& \zeta &\quad \text{weakly in } L^2(0,T;\cL^2_0),
\end{array}
\end{equation}
for some
$$
u \in L^\infty(0,T;\cH^s_0)\cap W^{1,2}(0,T;(\cH^1_0)'), \quad w \in L^2(0,T;\cH^1_0),\quad \zeta \in L^2(0,T;\cL^2_0).
$$
As a consequence of the Aubin-Lions-Simon compactness lemma and the compact embeddings $\cH^s_0 \hookrightarrow \cL^2_0 \hookrightarrow (\cH^1_0)'$, we obtain the stronger convergence
\begin{equation}\label{eq:strong}
u^\lambda \to u \quad \text{strongly in } C([0,T];\cL^2_0),
\end{equation}
up to considering another subsequence, that we do not relabel.

The second convergence in \eqref{eq:kilt} is again enough to obtain the convergence of the fractional Laplacians, as indeed
$$
\int_0^T \vartheta(t)\langle(-\Delta)^s u^\lambda, v\rangle_{(\cH^s_0)'}\dd t
\xrightarrow{\lambda\to0}
\int_0^T \vartheta(t)\langle(-\Delta)^s u, v\rangle_{(\cH^s_0)'}\dd t,
$$
for every $v\in \cH^s_0$ and every $\vartheta\in L^1(0,T)$.
Hence, passing to the limit as $\lambda\to0$ in \eqref{eq:CHlambdaint}, we see that the triple $(u,w,\zeta)$ satisfies the following equations
\begin{equation}\label{eq:malsai}\begin{aligned}
&\int_0^T \left(\io \partial_t u(x,t)v_1(x)+\nabla w(x,t)\cdot\nabla v_1(x) \dd x\right)\vartheta(t) \dd t = 0,\\
&
\int_0^T \left(\langle(-\Delta)^s u(t),v_2\rangle_{(\cH_0^s)'} + \io \big(\zeta(x,t)+\pi(u(x,t))-w(x,t)\big)v_2(x) \dd x\right)\vartheta(t) \dd t=0,
\end{aligned}\end{equation}
for every $v_1\in \cH^1_0$, $v_2\in\cH^s_0$ and every $\vartheta\in L^2(0,T)$.
Therefore the triple $(u,w,\zeta)$ satisfies equations \eqref{eq:CHw1}--\eqref{eq:CHw2} almost everywhere in $(0,T)$.

We now prove that the last convergence in \eqref{eq:kilt} and \eqref{eq:strong} ensure that
\begin{equation}\label{eq:barba}
\zeta(x,t) \in \gamma(u(x,t))\quad\mbox{for almost every }(x,t) \in \Omega \times (0,T).
\end{equation}

Indeed, the operator $\Upsilon \subset L^2(0,T;\cL^2_0) \times L^2(0,T;\cL^2_0)$ defined by
$$
\Upsilon \coloneqq \left\lbrace [u,\xi] \in L^2(0,T;\cL^2_0) \times L^2(0,T;\cL^2_0) \, : \, \xi(x,t) \in \gamma(u(x,t)) \, \text{ for a.\,e. } (x,t) \in \Omega \times (0,T) \right\rbrace
$$
is maximal monotone, and for $\lambda > 0$ and $u \in L^2(0,T;\cL^2_0)$ we can define the Yosida approximation of $\Upsilon$ as
$$
\Upsilon_\lambda(u) \coloneqq \frac1\lambda\big(I - (I+\lambda\Upsilon)^{-1}\big)(u).
$$
Note that
$$
((I+\lambda\Upsilon)^{-1}u)(x,t) = J_\lambda(u(x,t)), \quad \Upsilon_\lambda(u)(x,t) = \gamma_\lambda(u(x,t)), \quad \text{for a.\,e. } (x,t) \in \Omega \times (0,T).
$$
Since $u^\lambda \to u$ strongly in $C([0,T];\cL^2_0)$, and therefore in $L^2(0,T;\cL^2_0)$, we also have that $(I+\lambda\Upsilon)^{-1}u^\lambda \to u$ strongly in $L^2(0,T;\cL^2_0)$. Moreover, $\Upsilon_\lambda(u^\lambda) \rightharpoonup \zeta$ weakly in $L^2(0,T;\cL^2_0)$. Lemma~\ref{lemma:yosida}\,(i) yields $\Upsilon_\lambda(u^\lambda) \in \Upsilon((I+\lambda\Upsilon)^{-1}u^\lambda)$. Therefore, by the strong-weak closure of the maximal monotone operator $\Upsilon$ in $L^2(0,T;\cL^2_0) \times L^2(0,T;\cL^2_0)$ (see \cite[Lemma~2.3]{barbu2010nonlinear}), we conclude that $\zeta \in \Upsilon(u)$, that is, $\zeta(x,t) \in \gamma(u(x,t))$ for almost every $(x,t) \in \Omega \times (0,T)$.

We provide an alternative, more explicit and self-contained proof of \eqref{eq:barba} in Appendix~\ref{sec:MYreg}.

By arguing as in the end of Step~5 of Section~\ref{sec:exince}, and exploiting once again the uniqueness of the solution, we see that there is indeed no need to consider subsequences.

\textbf{Step 8: energy estimates.}
We relabel the triple obtained in the previous step as $(u_s,w_s,\zeta_s)$.
In order to prove \eqref{eq:en-est} and \eqref{eq:est-en}, we only need to pass to the limit in \eqref{eq:scusa}. For this we exploit the convergences proved in Step~7 of Section~\ref{sec:exince} and the following lemma.

\begin{lemma}
    We have
    $$
    E_s\leq\Gamma\textrm{-}\liminf_{\lambda\to0}E^\lambda,
    $$
    where the functionals are taken as defined on $\cH^s_0$, and the $\Gamma$-convergence is considered in the $\cL^2_0$-sense.
\end{lemma}

\begin{proof}
Let $u,u_k\in\cH^s_0$ such that $u_k\to u$ strongly in $\cL^2_0$ and $\lambda_k\to0$.
By \eqref{eq:ducks}, we have
$$
\io\Pi(u)\dd x= \lim_{k\to\infty}\io\Pi(u_k)\dd x.
$$
Moreover, by Lemma~\ref{lemma:parzialmenteserena} and arguing as in the proof of \eqref{eq:pandasaviour}, we obtain
$$
\frac{1}{2}\|u\|_{\cH^s_0}^2+\io\Gamma(u)\dd x\leq\liminf_{k\to\infty} \left(\frac{1}{2}\|u_k\|_{\cH^s_0}^2+\io\Gamma_{\lambda_k}(u_k)\dd x\right),
$$
thus concluding the proof.
\end{proof}

Thus, by passing to the limit in \eqref{eq:scusa}, we obtain
$$
\frac{1}{2}\int_0^t \|w_s(r)\|_{\cH^1_0}^2\dd r
+E_s(u_s(t))
\leq E_s(u_0),
$$
for every $t\in[0,T)$,
which is \eqref{eq:en-est}, and also implies \eqref{eq:est-en} with $\tau=0$.
The general form of \eqref{eq:est-en} is then a consequence of the uniqueness of the solution to \eqref{eq:CHw1}--\eqref{eq:CHw2}, which implies that the solution in the time interval $(0,T-\tau)$, corresponding to the initial datum $u_s(\tau)$, is actually $(u_s(\,\cdot\,+\tau),w_s(\,\cdot\,+\tau))$. More precisely, this argument yields the stronger estimate
$$
\frac{1}{2}\int_\tau^t \|w_s(r)\|_{\cH^1_0}^2\dd r
+E_s(u_s(t))
\leq E_s(u_s(\tau))
\qquad\mbox{for }0\leq\tau\leq t <T.
$$

For later use, it is also convenient to point out that, by the weak convergence of $\gamma_\lambda(u^\lambda)$ to $\zeta_s$, and the weak lower semicontinuity of the $L^2$-norm, we can pass to the limit in \eqref{eq:lest2}, obtaining
$$
\|\zeta_s\|_{L^2(0,T;\cL^2_0)}\leq 2\Big(\|w_s\|_{L^2(0,T;\cL^2_0)}+C_\pi \, \sqrt{T}\sup_{t\in[0,T]}\|u_s(t)\|_{\cL^2_0}\Big).
$$
Then, by \eqref{eq:en-est}, recalling Hypothesis \ref{stm:hyp}\,(iii) and using Young's inequality and the fractional Poincar\'e inequality, we obtain
\begin{align*}
    \sup_{t\in[0,T]}\|u_s(t)\|_{\cL^2_0}\leq C_{\textrm{P}}\sup_{t\in[0,T]}\|u_s(t)\|_{\cH^s_0}\leq C \big(E_s(u_0)+|\Omega|\big),
\end{align*}
which implies
\begin{equation}\label{eq:serveunnome}
\|\zeta_s\|_{L^2(0,T;\cL^2_0)}\leq C \big(1+\sqrt{T}\big)E_s(u_0)+\beta_1|\Omega|.
\end{equation}
We stress that the constants $C$ and $\beta_1$ do not depend on $u_0,s,T$.

\subsubsection{Proof of Theorem~\ref{stm:asymptotics}}
In order to prove the theorem, we need the following auxiliary result.

\begin{lemma}\label{lem:asympi}
Let $\{s_k\}_{k\in\nat} \subset (0,1)$ be a sequence such that $s_k \nearrow 1$ as $k\to\infty$. Let $\{v_k\} \subset \cL^2_0$ be a sequence such that $v_k \in \cH^{s_k}_0$ for all $k \in \nat$. Moreover, let us assume that $v_k$ converges to some function $v \in \cH^1_0$ strongly in $\cL^2_0$. Then we have
$$
\left\langle (-\Delta)^{s_k} v_k,\psi \right\rangle_{(\cH^s_0)'} \to C(d) \langle -\Delta v,\psi\rangle_{(\cH^1_0)'} \quad\textrm{for }\psi\in C^\infty_c(\Omega).
$$
\end{lemma}

\begin{proof}
We choose a test function $\psi \in C_c^\infty(\Omega)$. We have
$$
\begin{aligned}
\left\langle (-\Delta)^{s_k} v_k,\psi \right\rangle_{(\cH^s_0)'} &= (1-s_k) \iint_{\real^{2d}} \frac{(v_k(x)-v_k(y))}{|x-y|^{d+2s_k}}\,(\psi(x)-\psi(y)) \dd x \dd y \\
&= \io v_k(x) \left(-(1-s_k) \int_{\real^d} \frac{\psi(x+z)+\psi(x-z)-2\psi(x)}{|z|^{d+2s_k}} \dd z \right) \dd x \\
&= \io \left(v_k(x)-v(x)\right) (-\Delta)^{s_k}\psi(x) \dd x + \io v(x) \, (-\Delta)^{s_k}\psi(x) \dd x \\
&=: T_k' + T_k''.
\end{aligned}
$$
By the uniform convergence in Lemma~\ref{thm:deltasto1},
the first term is such that
$$
|T_k'| \le \io \left|v_k(x)-v(x)\right| \left|(-\Delta)^{s_k}\psi(x) \right| \dd x \le C(d,\psi)\,\|v_k-v\|_{\cL^2_0},
$$
while, 
$$
\Big|T_k'' -C(d) \io v(x) \, (-\Delta\psi(x)) \dd x \Big|\leq \io v(x)\big\|C(d)\Delta\psi-(-\Delta)^{s_k}\psi\big\|_{C(\real^d)}\dd x.
$$
By using also Lebesgue's dominated convergence theorem, we obtain
\begin{equation}\label{eq:bonsai}\begin{aligned}
&\Big|\left\langle (-\Delta)^{s_k} v_k,\psi \right\rangle_{(\cH^s_0)'} -C(d) \langle -\Delta v,\psi\rangle_{(\cH^1_0)'}\Big|\\
&\qquad\quad\leq C(d,\psi)\,\|v_k-v\|_{\cL^2_0}+|\Omega|^\frac{1}{2}\|v\|_{\cL^2_0}\big\|C(d)\Delta\psi-(-\Delta)^{s_k}\psi\big\|_{C(\real^d)}\to 0,
\end{aligned}\end{equation}
concluding the proof.
\end{proof}

We begin by observing that, by the Bourgain, Brezis, and Mironescu's asymptotics results, see \cite{bourgain2001another}, since $u_0\in\cH^1_0$, we have
\begin{equation}\label{eq:ADPM}
\limsup_{s\to1} E_s(u_0)\le C(d) \left(\frac{1}{2}\|u_0\|_{\cH^1_0}^2+\io F(u_0(x))\dd x\right)\eqqcolon E_1(u_0).
\end{equation}

Since the constant $C$ in \eqref{eq:en-est} does not depend on $s$, and $u_0\in\cH^1_0$, by arguing as in the proof of \eqref{eq:gargoyle} with $E_s$ in place of $E^\lambda$, we obtain
\begin{equation}\label{eq:en-esti}\begin{aligned}
\limsup_{s\to1}\left(\int_0^T\|w_s(t)\|^2_{\cH^1_0} \dd t + \sup_{t\in[0,T]}E_s(u_s(t))\right) 
\le C (E_1(u_0)+1).
\end{aligned}\end{equation}
Here above we have also exploited the fact that the Poincar\'e constant $C_{\textrm{P}}$, which comes from \cite[Proposition~2.5]{brasco2015stability}, does not depend on $s$. Recalling Hypothesis~\ref{stm:hyp}\,(iii), and arguing again as in the proof of \eqref{eq:gargoyle}, yields
\begin{equation}\label{eq:orbweaver}
\|u_s(t)\|_{\cL^2_0}^2\leq C_{\textrm{P}} \|u_s(t)\|^2_{\cH^s_0} \le C(E_1(u_0)+1) \quad \text{for every } t \in (0,T),
\end{equation}
uniformly for $s$ close to $1$.
Since $\partial_t u_s$ and $w_s$ are related by equation~\eqref{eq:CHw1}, we have also the estimate
\begin{equation}\label{eq:redback}
\|\partial_t u_s\|_{L^2(0,T;(\cH^1_0)')}^2\leq C.
\end{equation}
These two estimates do not depend on $s\in(0,1)$, hence there exists $u\in C([0,T];(\cH^1_0)')$ such that (up to a subsequence that we do not relabel)
\begin{equation}\label{eq:pianoscassato}
u_s\to u\quad\mbox{in }C([0,T];(\cH^1_0)').
\end{equation}
This is indeed a consequence of the Aubin-Lions-Simon compactness lemma, applied here with the compact embedding of $\cL^2_0$ in $(\cH^1_0)'$.

By \cite[Lemma~3.10]{brasco2015stability}, for every $t \in (0,T)$ there exist a function $v^t \in \cH^1_0$ and a sequence $\{s_h(t)\}_{h\in\nat}$ such that $s_h(t)\nearrow1$ and
\begin{equation}\label{eq:aestrongly}
u_{s_h(t)}(t) \to v^t \quad \text{strongly in } \cL^2_0.
\end{equation}
By \eqref{eq:pianoscassato} and the uniqueness of the limit, we conclude that $v^t=u(t)$, and the convergence in \eqref{eq:aestrongly} does not hold only for the sequence $s_h(t)$ (which could depend on $t$), but actually for the whole sequence for which we have the limit in \eqref{eq:pianoscassato}.
By \eqref{eq:serveunnome} and \eqref{eq:ADPM}, we obtain
$$
\|\zeta_s\|_{L^2(0,T;\cL^2_0)}\leq C,
$$
uniformly in $s\in(0,1)$. This, together with \eqref{eq:redback} and \eqref{eq:en-esti}, yields the following convergences
\begin{equation*}
\begin{array}{rcll}
\partial_tu_s &\rightharpoonup& \partial_t u &\quad \text{weakly in } L^2(0,T;(\cH^1_0)'), \\[2pt]
w_s&\rightharpoonup& w&\quad\mbox{weakly in }L^2(0,T;\cH^1_0), \\[2pt]
\zeta_s&\rightharpoonup&\zeta&\quad\mbox{weakly in }L^2(0,T;\cL^2_0),
\end{array}
\end{equation*}
once again up to a subsequence that we do not relabel.
We have thus obtained the convergences in \eqref{eq:aiuto} (for a certain sequence $s_h\nearrow1$).
Moreover, by \eqref{eq:orbweaver}, \eqref{eq:aestrongly} and Lebesgue's dominated convergence theorem, we can conclude that
\begin{equation}\label{eq:magpie}
u_s\to u\quad\mbox{strongly in }L^2(0,T;\cL^2_0).
\end{equation}
By the Lipschitz continuity of $\pi$, \eqref{eq:bonsai}, and the convergence in \eqref{eq:magpie} we conclude that
$$\begin{aligned}
&\pi(u_s)\to\pi(u)\quad\mbox{strongly in }L^2(0,T;\cL^2_0),\\
&\int_0^T \vartheta(t)\langle(-\Delta)^s u_s(t),v\rangle_{(\cH^s_0)'}\dd t
\to C(d)\int_0^T \vartheta(t)\langle-\Delta u(t),v\rangle_{(\cH^1_0)'}\dd t,
\end{aligned}$$
for every $\vartheta\in L^2(0,T)$ and $v\in C^\infty_c(\Omega)$. We can thus let $s\to1$ in \eqref{eq:malsai}, obtaining
\begin{equation*}\begin{aligned}
&\int_0^T \left(\io \partial_t u(x,t)v_1(x)+\nabla w(x,t)\cdot\nabla v_1(x) \dd x\right)\vartheta(t) \dd t = 0,\\
&
\int_0^T \left(C(d)\langle -\Delta u(t),v_2\rangle_{(\cH_0^1)'} + \io \big(\zeta(x,t)+\pi(u(x,t))-w(x,t)\big)v_2(x) \dd x\right)\vartheta(t) \dd t=0,
\end{aligned}\end{equation*}
for every $v_1\in \cH^1_0$, $v_2\in C^\infty_c(\Omega)$ and every $\vartheta\in L^2(0,T)$. By density, we can actually consider $v_2\in\cH^1_0$. Moreover, we point out that arguing as in the end of Step~7 of Section~\ref{sec:exince} and exploiting \cite[Lemma~2.3]{barbu2010nonlinear} ensures that $\zeta(x,t)\in\gamma(u(x,t))$ for almost every $(x,t)\in\Omega\times(0,T)$. Thus, by applying the fundamental lemma of the calculus of variations for Bochner spaces, we can conclude that $(u,w,\zeta)$ is a solution to the classical Cahn-Hilliard system in the sense of Definition~\ref{def:weaksol_local}.
It is well-known that such a solution is unique (for a simple proof, one can adapt, e.\,g., the argument in Step~1 of Section~\ref{sec:exince}). As a consequence, the limits considered above hold for the whole sequence $s\to1$. This concludes the proof of the asymptotics result.


\subsection{A more general framework}

In place of the Laplacian, we can consider a more general (possibly nonlocal) operator $\mathfrak L$, in \eqref{eq:CH1}.
Moreover, in \eqref{eq:CH2} we can substitute the fractional Laplacian with a nonlocal integrodifferential operator $\mathfrak I$ and solve the system
\begin{align}
\partial_t u + \mathfrak L w = 0 &\qquad \text{in } \Omega \times (0,T), \label{eq:CH1q}\\
w = \mathfrak I u + F'(u) &\qquad \text{in } \Omega \times (0,T), \label{eq:CH2q} \\
u(x,0) = u_0(x) &\qquad \text{in } \real^d, \label{eq:iniCHq} \\
u = 0 &\qquad \text{in } (\real^{d}\setminus\Omega)\times(0,T), \label{eq:bouCH1q}\\
w = 0 &\qquad \text{in } (\real^{d}\setminus\Omega)\times(0,T). \label{eq:bouCH2q}
\end{align}
We are thus in the setting of Metatheorem~\ref{stm:Metathm1}.

Concerning equation \eqref{eq:CH1q}, the key assumption is that $\mathfrak L:X_0\to(X_0)'$ is a linear invertible operator defined on a Hilbert space $X_0$ densely embedded in $\cL^2_0$, which we recall is defined in \eqref{eq:L20}.
Moreover, albeit not strictly necessary for the solvability of the above system, in order for condition \eqref{eq:bouCH2q} to make sense, the space $X_0$ should also have a meaningful notion of ``boundary datum'' (either in the local or nonlocal sense).
More precisely, $\mathfrak L$ must have a variational nature, being associated to an equivalent norm on $X_0$. For the details we refer to Section~\ref{sec:dual}.

The operator $\mathfrak I$ that we consider in equation \eqref{eq:CH2q} is of $(s,q)$-fractional Laplacian type, as defined in \eqref{eq:nooooooo}, with the corresponding functional framework of Section~\ref{sec:HDbc}. For the precise details we refer to Section~\ref{sec:sqdef}.
The potential $F$ satisfies Hypothesis~\ref{stm:hyp}.

\subsubsection{Existence of the solution}\label{sec:waywardson}
With the aforementioned hypotheses, weak solutions to the corresponding Cahn-Hilliard system are defined as follows.

\begin{definition}[Solution to the generalized fractional Cahn-Hilliard system]\label{def:weaksolq}
Let $T>0$ be fixed. We say that $(u,w,\zeta)$ is a weak solution to the Cahn-Hilliard system \eqref{eq:CH1q}--\eqref{eq:bouCH2q} associated with the initial datum $u_0 \in \cW^{K,q}_0(\Omega)$ if
$$
\begin{aligned}
u &\in L^\infty(0,T;\cW^{K,q}_0(\Omega)) \cap W^{1,2}(0,T;(X_0)'), \\
w &\in L^2(0,T;X_0),\\
\zeta &\in L^2(0,T;\cL^2_0), \qquad \zeta\in\gamma(u) \quad \text{a.\,e. in } \Omega \times (0,T),
\end{aligned}
$$
and $(u,w,\zeta)$ satisfies the following weak formulation of \eqref{eq:CH1q}--\eqref{eq:CH2q}:
\begin{align}
\partial_t u + \mathfrak L w = 0 &\qquad \text{in } (X_0)' \label{eq:CHw1q}\\
w = \mathfrak I u + \zeta + \pi(u) &\qquad \text{in } (\cW^{K,q}_0(\Omega))', \label{eq:CHw2q}
\end{align}
almost everywhere in $(0,T)$, with $u(0) = u_0$ almost everywhere in $\real^d$.
\end{definition}

By Remark~\ref{rmk:Kqcompact}, $\cW^{K,q}_0(\Omega)$ is densely and compactly embedded in $\cL^2_0$.
Moreover, by assumption, $X_0$ is densely embedded in $\cL^2_0$. We then have the embeddings $\cW^{K,q}_0(\Omega)\hookrightarrow\cL^2_0\simeq(\cL^2_0)'\hookrightarrow(X_0)'$, which, as in Remark~\ref{rmk:AubinLions}, entails
$$
L^\infty(0,T;\cW^{K,q}_0(\Omega)) \cap W^{1,2}(0,T;(X_0)') \subset C([0,T];\cL^2_0).
$$

From a rigorous point of view, throughout this section we should have kept track of the dense embedding $i : X_0 \hookrightarrow \cL^2_0$. However, with a slight abuse of notation, we tacitly identify $w$ with $i(w)$ for the reader's convenience.

\begin{theorem}[Existence and uniqueness]\label{stm:existenceq}
Let Hypothesis~\ref{stm:hyp} be satisfied. Assume
$$
u_0 \in \cW^{K,q}_0(\Omega), \qquad \Gamma(u_0) \in L^1(\real^d).
$$
Then the Cahn-Hilliard system \eqref{eq:CH1q}--\eqref{eq:bouCH2q} admits a solution $(u,w,\zeta)$, in the sense of Definition~\ref{def:weaksolq} for every $T>0$. Moreover, defining the energy
$$
\mathscr E(u) \coloneqq \iint_{Q(\Omega)}\Phi(u(x)-u(y))K(x,y)\dd x\dd y + \io F(u(x))\dd x,
$$
the following energy estimate holds:
\begin{equation}\label{eq:en-estq}
\frac{1}{2}\int_0^t \|w(\tau)\|_{X_0}^2 \dd \tau +  \mathscr E(u(t)) \le \mathscr E(u_0)\qquad\mbox{for every }t\geq 0.
\end{equation}
Furthermore, if $\phi$ satisfies condition \eqref{eq:nameless1}, then the solution is unique, and
\begin{equation}\label{eq:est-enq}
\mathscr E(u(t)) \le \mathscr E(u(\tau)) \qquad \mbox{for } \ 0 \le \tau \le t.
\end{equation}
\end{theorem}

\subsubsection{Proof of Theorem~\ref{stm:existenceq}}\label{sec:proofDirGenThm}
In order to prove Theorem~\ref{stm:existenceq}, we need the analogues of Lemma~\ref{lemma:neutralquokka} and Lemma~\ref{lemma:emma}. For this we introduce the regularized energy
$$
\mathscr E^\lambda(u) \coloneqq \mathfrak F(u) + \io\big( \Gamma_\lambda(u)+\Pi(u)\big)\dd x
$$
for $\lambda\in(0,1)$, with $\mathfrak F$ as defined in \eqref{eq:effe}.

\begin{lemma}\label{lemma:neutralquokkaq}
The functional $\mathscr E^\lambda$ is lower semicontinuous with respect to the $L^2$-convergence. Moreover
	\begin{equation}\label{eq:gargoyleq}
		\mathscr E^\lambda(u) \geq \frac{1}{2\Lambda q}\|u\|_{K,q}^q-\beta_0|\Omega|\qquad\mbox{for every }u\in\cW^{K,q}_0(\Omega),
	\end{equation}
	for $\lambda\in(0,\bar\lambda)$, 
	where $\beta_0>0$ and $\bar\lambda\in(0,1)$ depend only on $d,\Lambda,s,\varrho,q,\Omega,a_1,a_2$ and $C_\pi$.
\end{lemma}

\begin{proof}
For the semicontinuity, we only need to check that if $\{v_k\}\subset\cL^2_0$ is such that $v_k\to v\in\cL^2_0$ strongly in $L^2(\Omega)$, then
$$
	\|v\|_{K,q}^q \leq \liminf_{k\to\infty}\|v_k\|_{K,q}^q.
$$
This follows by arguing as in the proof of \eqref{eq:pandasaviour}, exploiting Fatou's Lemma and recalling that $\Phi$ is continuous and nonnegative. The estimate \eqref{eq:gargoyleq} follows again by using the fractional Poincar\'e inequality \eqref{eq:gen_Poincare} and the fact that $q\geq2$. Indeed, recalling \eqref{eq:koala} and using Young's inequality,
$$
\begin{aligned}
    \mathscr E^\lambda(u)&\geq \frac{1}{\Lambda q}\|u\|_{K,q}^q-\alpha\lambda^\frac{1}{2}\|u\|_{L^2(\Omega)}^2-a_3\|u\|_{L^p(\Omega)}^p-\beta|\Omega|\\
    &
    \geq \frac{1}{\Lambda q}\|u\|_{K,q}^q-C(\lambda^\frac{1}{2}+\varepsilon)\|u\|_{L^2(\Omega)}^2-c(\varepsilon)|\Omega|\\
    &
    \geq \frac{1}{\Lambda q}\|u\|_{K,q}^q-C(\lambda^\frac{1}{2}+\varepsilon)\iint_{\real^{2d}\cap\{|x-y|<\varrho\}}\frac{|u(x)-u(y)|^q}{|x-y|^{d+sq}}\dd x\dd y-c(\varepsilon)|\Omega|\\
    &
    \geq \frac{1}{\Lambda q}\|u\|_{K,q}^q-\frac{C(\lambda^\frac{1}{2}+\varepsilon)}{\Lambda}\|u\|_{K,q}^q-c(\varepsilon)|\Omega|\\
    &
    \geq \frac{1}{2\Lambda q}\|u\|_{K,q}^q-\beta_0|\Omega|,
\end{aligned}
$$
as claimed.
\end{proof}

The coerciveness proved in Lemma~\ref{lemma:neutralquokkaq}, together with the compact embedding of $\cW^{K,q}_0(\Omega)$ in $\cL^2_0$ (see Remark~\ref{rmk:Kqcompact}), are enough to prove the following existence result, by arguing as in the proof of Lemma~\ref{lemma:emma}.

\begin{lemma}\label{lemma:emmaq}
	Given $\lambda\in(0,\bar\lambda)$, $\tau>0$ and any $g\in\cL^2_0$, we define the energy $\mathscr E_g^\lambda:\cW^{K,q}_0(\Omega)\to\real$ by setting
	$$
	\mathscr E_g^\lambda(u) \coloneqq \frac{1}{2\tau}\|u-g\|^2_{(X_0)'} + \mathscr E^\lambda(u).
	$$
	Then, there exists at least one function $u_\star\in\cW^{K,q}_0(\Omega)$ such that
	$$
	\mathscr E_g^\lambda(u_\star) = \inf_{u\in\cW^{K,q}_0(\Omega)}\mathscr E_g^\lambda(u).
	$$
 In particular, $u_\star$ satisfies
 \begin{equation}\label{eq:moreover}
\mathfrak L^{-1}\left(\frac{u_\star-g}{\tau}\right)+\mathfrak I u_\star+\gamma_\lambda(u_\star)+\pi(u_\star)=0\qquad\mbox{in }(\cW^{K,q}_0(\Omega))'.
 \end{equation}
\end{lemma}

The Euler-Lagrange equation \eqref{eq:moreover} is obtained by recalling that $\mathfrak L^{-1}$ is the Fréchet derivative of $\|\,\cdot\,\|_{\mathfrak L^{-1}}^2/2$ and applying Lemma~\ref{lem:nonhouncomodino}.

We subdivide the remaining part of the proof into the following steps.
We stress that the lower bound in \eqref{eq:nameless} is not sufficient to ensure the uniqueness of the solution, since, a priori, the operator $\mathfrak I$ is not linear. In any case, this does not affect the proof of the existence.

\textbf{Step 1: time-discretization.}
We begin by fixing $T>0$ and we proceed as in the proof of Theorem~\ref{stm:existence}.
We consider $N\in\nat$ and $\tau=\tau_N \coloneqq T/N$. Arguing as in Step~3 of Section~\ref{sec:exince}, we first apply Lemma~\ref{lemma:emmaq} to pick a minimizer $u^\lambda_1\in\cW^{K,q}_0(\Omega)$ of the functional
$$
\mathscr E_1^\lambda(u) \coloneqq \frac{\tau}{2}\left\|\frac{u-u_0}{\tau}\right\|^2_{(X_0)'} + \mathscr E^\lambda(u).
$$
Then we proceed iteratively, choosing $u_n^\lambda\in\cW^{K,q}_0(\Omega)$ to be a minimizer of
\begin{equation}\label{eq:carryon}
\mathscr E_n^\lambda(u) \coloneqq \frac{\tau}{2}\left\|\frac{u-u_{n-1}}{\tau}\right\|^2_{(X_0)'} + \mathscr E^\lambda(u),
\end{equation}
and defining $w_n^\lambda\in X_0$ as in \eqref{eq:wnl}, i.\,e.,
$$
w_n^\lambda=(-\mathfrak L)^{-1}\left(\frac{u^\lambda_n - u^\lambda_{n-1}}{\tau}\right).
$$
Since $(-\mathfrak L)^{-1}=-\mathfrak L^{-1}$, and recalling \eqref{eq:moreover}, the pair $(u_n^\lambda,w_n^\lambda)$ satisfies
\begin{align*}
	\frac{u^\lambda_n - u^\lambda_{n-1}}{\tau} + \mathfrak L w^\lambda_n = 0 &\qquad \text{in } (X_0)',\\
	w^\lambda_n = \mathfrak I u^\lambda_n + \gamma_\lambda(u^\lambda_n) + \pi(u^\lambda_n) &\qquad \text{in } (\cW^{K,q}_0(\Omega))',
\end{align*}
for every $n=1,\ldots,N$. Moreover we have the counterpart of the energy estimate \eqref{eq:underclasshero}, that is
\begin{equation}\label{eq:somesay}
\frac{\tau}{2}\sum_{n=1}^{\bar{n}} \left\|\frac{u_n^\lambda - u^\lambda_{n-1}}{\tau}\right\|_{(X_0)'}^2 + \mathscr E^\lambda(u_{\bar{n}}^\lambda)\le \mathscr E(u_0).
\end{equation}
We also recall that, by \eqref{eq:lorenzace} we have
$$
\|w_n^\lambda\|_{X_0}=\left\|\frac{u_n^\lambda - u^\lambda_{n-1}}{\tau}\right\|_{(X_0)'}.
$$

\textbf{Step 2: convergence as $\tau\to0$.}
We then define the interpolating functions as in \eqref{eq:inter}. From \eqref{eq:somesay} we deduce the analogues of \eqref{eq:esti1_tau} and \eqref{eq:esti2_tau}, that lead to the convergence of the interpolants as $\tau\to0$, as in \eqref{eq:conv} and \eqref{eq:conv_ALS}, via compactness, to a pair $(u^\lambda,w^\lambda)$ such that
$$
u^\lambda \in L^\infty(0,T;\cW^{K,q}_0(\Omega))\cap W^{1,2}(0,T;(X_0)'), \quad w^\lambda \in L^2(0,T;X_0).
$$
More precisely,
\begin{equation*}
\begin{array}{rcll}
	\bar{w}^\lambda_\tau &\rightharpoonup& w^\lambda &\quad \text{weakly in } L^2(0,T;X_0), \\[2pt]
	\hat{u}^\lambda_\tau &\to& u^\lambda &\quad \text{strongly in } C([0,T];\cL^2_0), \\[2pt]
 \bar{u}^\lambda_\tau &\to& u^\lambda &\quad \text{strongly in } L^\infty(0,T;\cL^2_0), \\[2pt]
	\partial_t\hat{u}^\lambda_\tau &\rightharpoonup& \partial_t u^\lambda &\quad \text{weakly in } L^2(0,T;(X_0)'),
\end{array}
\end{equation*}
However, here the convergence only holds true for a specific subsequence $\{\tau_{N_h}\}$, since we cannot rely on the uniqueness of solutions.

For the convergence of $\mathfrak I \bar{u}^\lambda_\tau$ we adapt the classical argument for $(-\Delta)^s_q$ (see, e.\,g. \cite[Proof of Theorem~3.4]{baroncini2018continuity}).
We begin by defining
$$
\xi_h^\lambda(x,y,t)\coloneqq\phi\Big(\bar{u}_{\tau_{N_h}}^\lambda(x,t)-\bar{u}_{\tau_{N_h}}^\lambda(y,t)\Big)K(x,y)^\frac{1}{q'},
$$
where $q'$ is the conjugate exponent of $q$, i.\,e., $q'=q/(q-1)$.
Then, the $L^{q'}(\real^{2d}\times(0,T))$-norm of $\xi_h^\lambda$ is uniformly bounded as a consequence of the upper bound in \eqref{eq:nameless}, the coercivity in \eqref{eq:gargoyleq}, and the energy estimate \eqref{eq:somesay}. Indeed,
\begin{equation}\label{eq:uno}\begin{aligned}
\|\xi_h^\lambda\|_{L^{q'}(\real^{2d}\times(0,T))}^{q'}&
=\int_0^T\iint_{\real^{2d}}\Big|\phi\Big(\bar{u}_{\tau_{N_h}}^\lambda(x,t)-\bar{u}_{\tau_{N_h}}^\lambda(y,t)\Big)\Big|^{q'}K(x,y)\dd x\dd y\dd t\\
&\leq\Lambda^{q'}\int_0^T\|\bar{u}_{\tau_{N_h}}^\lambda(t)\|_{K,q}^q\dd t\leq CT\big(\mathscr E(u_0)+1\big),
\end{aligned}\end{equation}
uniformly in $h$ and $\lambda$. Since $\phi$ is continuous and $\bar{u}_{\tau_{N_h}}^\lambda(x,t)\to u^\lambda(x,t)$ pointwise almost everywhere in $\real^d\times(0,T)$, we have
\begin{equation}\label{eq:due}
\xi_h^\lambda(x,y,t)\to \phi\big(u^\lambda(x,t)-u^\lambda(y,t)\big)K(x,y)^\frac{1}{q'}\eqqcolon\xi^\lambda(x,y,t),
\end{equation}
for almost every $(x,y,t)\in\real^{2d}\times(0,T)$.
Thus, by \cite[Theorem~8.59]{ziemer2017modern}, the boundedness in \eqref{eq:uno} and the convergence in \eqref{eq:due} imply the weak convergence
\begin{equation*}
\xi_h^\lambda\rightharpoonup\xi^\lambda\quad\mbox{weakly in }L^{q'}(\real^{2d}\times(0,T)).
\end{equation*}
Therefore we have
\begin{equation}\label{eq:slowpoke}\begin{aligned}
\int_0^T\vartheta&(t)\langle\mathfrak I \bar{u}_{\tau_{N_h}}^\lambda(t),v\rangle_{K,q}\dd t\\
&=
\int_0^T\vartheta(t)\iint_{\real^{2d}}\xi_h^\lambda(x,y,t)\big(v(x)-v(y)\big)K(x,y)^\frac{1}{q}\dd x\dd y\dd t
\xrightarrow{h\to\infty}\int_0^T\vartheta(t)\langle\mathfrak I u^\lambda(t),v\rangle_{K,q}\dd t,
\end{aligned}\end{equation}
for every $\vartheta\in L^q(0,T)$ and $v\in\cW^{K,q}_0(\Omega)$.
Indeed, if we set
$$
\Theta(x,y,t)\coloneqq \vartheta(t)\big(v(x)-v(y)\big)K(x,y)^\frac{1}{q},
$$
then
$$
\|\Theta\|_{L^{q}(\real^{2d}\times(0,T))}^{q}=\int_0^T |\vartheta(t)|^q\iint_{\real^{2d}}|v(x)-v(y)|^q K(x,y)\dd x\dd y\dd t=\|\vartheta\|_{L^q(0,T)}^q \|v\|_{K,q}^q.
$$
Hence, by passing to the limit $h\to\infty$, the limit functions $(u^\lambda,w^\lambda)$ satisfy the system
\begin{equation}\label{eq:belnome}\begin{aligned}
&\int_0^T \left(\io \partial_t u^\lambda(x,t)v_1(x) \dd x + \langle \mathfrak L w^\lambda(t),v_1 \rangle_0 \right)\vartheta(t) \dd t = 0,\\
&
\int_0^T \left(\langle\mathfrak I u^\lambda(t),v_2\rangle_{K,q} + \io \big(\gamma_\lambda(u^\lambda(x,t))+\pi(u^\lambda(x,t))-w^\lambda(x,t)\big)v_2(x) \dd x\right)\vartheta(t) \dd t=0,
\end{aligned}\end{equation}
for every $\vartheta\in L^q(0,T)$, $v_1\in X_0$ and $v_2\in\cW^{K,q}_0(\Omega)$.

\textbf{Step 3: convergence as $\lambda\to0$.}
Passing to the limit $\tau\to0$ also yields the analogues of \eqref{eq:lest}.
We then need to show the validity of \eqref{eq:lest2}. Since $\gamma_\lambda(0)=0$ and $\gamma_\lambda$ is Lipschitz continuous, we have $\gamma_\lambda(u^\lambda)\in\cW^{K,q}_0(\Omega)$, hence we can choose $v_2=\gamma_\lambda(u^\lambda)$ and $\vartheta\equiv 1$ in the second equation of \eqref{eq:belnome}, obtaining the analogue of \eqref{eq:basketcase}. Since $\phi(r)r\geq0$, and $\gamma_\lambda$ is increasing, we see that
$$
\int_0^T\iint_{\real^{2d}} \phi(u^\lambda(x,t)-u^\lambda(y,t))\big(\gamma_\lambda(u^\lambda(x,t))-\gamma_\lambda(u^\lambda(y,t))\big)K(x,y) \dd x\dd y\dd t\geq0.
$$
This yields \eqref{eq:lest2}, as wanted. 
These estimates are enough to ensure the passage to the limit of $u^\lambda,\gamma_\lambda(u^\lambda)$ and $w^\lambda$, for a certain subsequence $\lambda_k\to0$, obtaining
\begin{equation*}
\begin{array}{rcll}
	w^{\lambda_k} &\rightharpoonup& w &\quad \text{weakly in } L^2(0,T;X_0), \\[2pt]
	u^{\lambda_k} &\to& u &\quad \text{strongly in } C([0,T];\cL^2_0), \\[2pt]
	\partial_t u^{\lambda_k} &\rightharpoonup& \partial_t u &\quad \text{weakly in } L^2(0,T;(X_0)'),\\[2pt]
    \gamma_{\lambda_k}(u^{\lambda_k}) &\rightharpoonup& \zeta &\quad \text{weakly in } L^2(0,T;\cL^2_0).
\end{array}
\end{equation*}

We are left to show that the triple $(u,w,\zeta)$ satisfies \eqref{eq:CHw1q}--\eqref{eq:CHw2q}. 
We observe that by \eqref{eq:uno} and \eqref{eq:due} we obtain
$$
\|\xi^{\lambda_k}\|_{L^{q'}(\real^{2d}\times(0,T))}^{q'}\leq C\,T(\mathscr E(u_0)+1),
$$
uniformly in $k$. Moreover,
$$
\xi^{\lambda_k}(x,y,t)\to \phi(u(x,t)-u(y,t))K(x,y)^\frac{1}{q'},
$$
for almost every $(x,y,t)\in\real^{2d}\times(0,T)$.
Hence,
\begin{equation}\label{eq:repetto}
\xi^{\lambda_k}\rightharpoonup \phi(u(x,t)-u(y,t))K(x,y)^\frac{1}{q'}\quad\mbox{weakly in }L^{q'}(\real^{2d}\times(0,T)).
\end{equation}
We can thus pass to the limit $\lambda_k\to0$ in \eqref{eq:belnome} and use the fundamental lemma of the calculus of variations for Bochner spaces, to obtain \eqref{eq:CHw1q}--\eqref{eq:CHw2q} for almost every $t\in(0,T)$, as wanted.
Moreover we point out that, arguing as in the end of Step~7 of Section~\ref{sec:exince}, ensures that $\zeta(x,t)\in\gamma(u(x,t))$ for almost every $(x,t)\in\Omega\times(0,T)$.
Finally, the energy estimate \eqref{eq:en-estq} can be derived from \eqref{eq:somesay} as in Step~8 of Section~\ref{sec:exince}.

\textbf{Step 4: global in time solution.}
So far we have proved that for every $T>0$ there exists (at least) one solution $(u,w)$ of \eqref{eq:CHw1q}--\eqref{eq:CHw2q}. In order to construct a global in time solution we can proceed as follows. Fix $T=1$ and let $(u_1,w_1,\zeta_1)$ be a weak solution corresponding to the initial datum $u_0$. Then, consider a weak solution $(u_2,w_2,\zeta_2)$ corresponding to the initial datum $u_1\big(\frac{3}{4}\big)$.
We proceed iteratively by 
choosing $(u_{k+1},w_{k+1},\zeta_{k+1})$ to be a weak solution corresponding to the initial datum $u_{k}\big(\frac{3k}{4}\big)$. Finally, given $t\geq0$, we define
\begin{equation}\label{eq:superattak}
(u(t),w(t),\zeta(t)) \coloneqq \big(u_k(t-3k/4),w_k(t-3k/4),\zeta(t-3k/4)\big),
\end{equation}
for $k$ such that $t\in\big[\frac{3k}{4},\frac{3(k+1)}{4}\big)$.
This is indeed a global in time solution.

\textbf{Step 5: uniqueness.}
The stronger condition \eqref{eq:nameless1} ensures that
\begin{equation}\label{eq:Okinawa}
\frac{1}{\Lambda}\|u_1-u_2\|_{K,q}^q\leq\langle\mathfrak I u_1-\mathfrak I u_2,u_1-u_2\rangle_{K,q},
\end{equation}
for every $u_1,u_2\in\cW^{K,q}_0(\Omega)$.
Thus, if $(u_1,w_1,\zeta_1)$ and $(u_2,w_2,\zeta_2)$ are two weak solutions in the sense of Definition~\ref{def:weaksolq}, associated to the same initial datum $u_0$, we obtain
$$
\langle w_1-w_2,u_1-u_2\rangle_{K,q}
\geq \frac{1}{\Lambda}\|u_1-u_2\|_{K,q}^q- C_\pi\|u_1-u_2\|_{\cL^2_0}^2,
$$
hence, using Ehrling's lemma,
\begin{equation}\label{eq:fantasmaformaggino}
\frac{\dd}{\dd t}\|u_1-u_2\|^2_{(X_0)'}+\|u_1-u_2\|^q_{K,q}\leq C \|u_1-u_2\|^2_{(X_0)'},
\end{equation}
for almost every $t\in(0,T)$.
This is enough to conclude the desired uniqueness, via Gr\"onwall's inequality. As a consequence, we also infer that $\zeta_1=\zeta_2$.

The energy estimate \eqref{eq:est-enq} follows by arguing as in Step~8 of Section~\ref{sec:exince}.


\subsection{Examples of admissible kernels}\label{sec:paprika}

In order to show the generality of our assumptions, we provide some meaningful and non-standard examples of admissible kernels for the operator $\mathfrak I$ appearing in equation \eqref{eq:CHw2q}.
For a fixed $q$, we say that $K$ is an admissible kernel if there exist $\varrho>0$ and $s\in (0,1)$ for which $K$ satisfies \eqref{eq:drogasintetica}--\eqref{eq:Kintegrability}.

We first observe that, if $K$ is admissible and $a$ is
a bounded Borel-measurable function $a:\real^d\times\real^d\to[\underline{a},\overline{a}]$, with $0<\underline{a}\leq\overline{a}$, then $a K$ is still an admissible kernel.

Moreover, if $K_1$ and $K_2$ are kernels satisfying \eqref{eq:drogasintetica}--\eqref{eq:Kintegrability}, with $s_1,\varrho_1$ and $s_2,\varrho_2$ respectively, then their sum $K_1+K_2$ is admissible. More in general, given a finite number of (not necessarily disjoint) measurable sets $A_1,\dots,A_n\subset\real^{2d}$ such that
$$
\{(x,y)\in\real^{2d}\,:\,|x-y|<\varrho\}\subset\bigcup_{i=1}^n A_i,
$$
for some $\varrho>0$, and $K_i$ is a kernels satisfying \eqref{eq:drogasintetica}--\eqref{eq:Kintegrability}, with $s_i,\varrho_i$, for every $i=1,\dots,n$, then $\chi_{A_1}K_1+\dots+\chi_{A_n}K_n$ is also admissible.

Some concrete examples are the following:
\begin{itemize}
    \item the kernel $K_{s_1,s_2}(x,y)=|x-y|^{-d-s_1q}+|x-y|^{-d-s_2q}$, which allows us to consider the sum of the operators $(-\Delta)^{s_1}_q+(-\Delta)^{s_2}_q$ within our functional framework, by interpreting it as a single operator $\mathfrak I$ given by $\phi(r)=|r|^{q-2}r$ and $K_{s_1,s_2}$;
    \item we can distinguish different intensities for the interactions confined within a certain region $A\subset\real^d$ and the long-ranged ones, by considering $K(x,y)=|x-y|^{-d-s_1q}\chi_{A^2}(x,y)+|x-y|^{-d-s_2q}\chi_{\real^{2d}\setminus A^2}(x,y)$;
    \item more in general, the interaction intensity can depend on the specific pair of points $(x,y)$, namely $K(x,y)=|x-y|^{-d-s(x,y)q}$, with $s:\real^d\times\real^d\to[s_0,s_1]\subset(0,1)$ Borel-measurable.
\end{itemize}


\subsection{Examples of admissible operators $\mathfrak L$}

Since the operator $\mathfrak L$ has an abstract and very general formulation, besides usual examples like the (fractional) Laplacian, we can consider some unusual ones. These examples are the content of the upcoming subsections.

\subsubsection{Fractional diffusion}
A particular example of operator that we can consider in equation \eqref{eq:CH1q} is given by  $\mathfrak L = \mathfrak I$ as in \eqref{eq:nooooooo}, with the choice $\Phi(r) = |r|^2/2$.
More precisely, we consider a kernel $K_\star$, that is a Borel-measurable function $K_\star:\real^d\times\real^d\to[0,+\infty]$, satisfying the following hypotheses:
\begin{align*}
&\mbox{Singularity:}\quad \frac{1}{\Lambda_\star}\chi_{B_{\varrho_\star}}(x-y)\leq K_\star(x,y)|x-y|^{d+2\sigma}\quad\mbox{for every }x \neq y, \\[2pt]
&\mbox{Integrability:}\quad \iint_{Q(\Omega)}\min\{1,|x-y|^2\}K_\star(x,y)\dd x\dd y<+\infty,
\end{align*}
for some $\sigma\in(0,1)$ and $\Lambda_\star\ge1$. Then, we define $X_0\coloneqq\cW^{K_\star,2}_0(\Omega)$ and $\mathfrak L:X_0\to(X_0)'$ by setting
$$
\langle\mathfrak L u,v\rangle_0\coloneqq \iint_{\real^{2d}}\big(u(x)-u(y)\big)\big(v(x)-v(y)\big)K_\star(x,y)\dd x\dd y\quad\mbox{for }u,v\in X_0,
$$
and we consider system \eqref{eq:CH1q}--\eqref{eq:bouCH2q}.

Therefore, Theorem~\ref{stm:existenceq} can be applied in the particular case in which $\mathfrak L=(-\Delta)^\sigma$ is the fractional Laplacian (thus extending the existence result proved in \cite{akagi2016fractional}).

\subsubsection{Generalized fractional Allen-Cahn}\label{sec:AllenCahn}
As a byproduct of the strategy that we employed to prove existence of solutions to the Cahn-Hilliard system, we are able to derive analogous results for the Allen-Cahn equation
\begin{equation}\label{eq:tristero}
\begin{aligned}
\partial_t u+\mathfrak I u+F'(u)=0 &\qquad \text{in } \Omega \times (0,T), \\
u(x,0) = u_0(x) &\qquad \text{in } \real^d, \\
u = 0 &\qquad \text{in } (\real^{d}\setminus\Omega)\times(0,T).
\end{aligned}
\end{equation}
Indeed, this can be interpreted as a Cahn-Hilliard system of the form \eqref{eq:CH1q}--\eqref{eq:bouCH2q}, by considering $X_0=\cL^2_0$, and $\mathfrak L$ as the standard Riesz map of $\cL^2_0$, that is $\mathfrak L u=(u,\,\cdot\,)_{\cL^2_0}$.
Thus, Theorem~\ref{stm:existenceq} yields the existence of a solution to \eqref{eq:tristero}. 

Clearly, a more direct way to prove existence consists in discretizing in time, considering the energy
$$
\mathscr E_n^\lambda(u) \coloneqq \frac{\tau}{2}\left\|\frac{u-u_{n-1}}{\tau}\right\|^2_{\cL^2_0} + \mathscr E^\lambda(u)
$$
in place of \eqref{eq:carryon}, defining iteratively $u_n^\lambda\in\cW^{K,q}_0(\Omega)$ as its minimizer, and carrying out the approximation argument of Section~\ref{sec:waywardson} (forgetting about $w$).

\subsubsection{Higher order operators}\label{sec:high}
In the first equation we can consider higher order (integro)differential operators, in particular the fractional Laplacian of order $s \in (0,+\infty) \setminus \nat$, which is obtained by generalizing the integrodifferential definition of the usual fractional Laplacian, see \cite[formula~(1.2)]{abatangelo2018loss}. As an example, we give a brief overview in the case $s \in (1,2)$. For $u \in C^4(\real^d) \cap L^{\infty}(\real^d)$ we can define
$$
(-\Delta)^s u(x) \coloneqq C(s,d) \PV \int_{\real^2} \frac{u(x+2y)-4u(x+y)+6u(x)-4u(x-y)+u(x-2y)}{|y|^{d+2s}} \dd y,
$$
where $C(s,d)$ is the normalization constant appearing in the equivalence with the pseudodifferential form of the fractional Laplacian, see \cite[formula~(2.1)]{saldana2018fractional}. This operator has a corresponding variational framework, see \cite{abatangelo2018loss}. To be precise, for $\Omega \subset \real^d$ bounded open set with Lipschitz boundary, we consider the space
$$
\cH^s_0 \coloneqq \left\lbrace u \in H^s(\real^d) \, : \, u \equiv 0 \mbox{ a.\,e.\ in } \real^d\setminus\Omega \right\rbrace,
$$
where, as usual,
$$
H^s(\real^d) \coloneqq \left\lbrace u \in L^2(\real^d) \, : \, (1+|\xi|^2)^{\frac{s}{2}}\mathscr F(u) \in L^2(\real^d) \right\rbrace.
$$
We then have the bilinear operator $B : \cH^s_0 \times \cH^s_0 \to \real$ given by \cite[formulas~(3.2)--(3.4)]{saldana2018fractional}, namely,
$$
\begin{aligned}
B(u,v) &= \int_{\real^d} |\xi|^{2s}\mathscr Fu(\xi)\mathscr Fv(\xi) \dd\xi \\
&= C_1 \int_{\real^d}\int_{\real^d} \frac{(\nabla u(x)-\nabla u(y))\cdot(\nabla v(x)-\nabla v(y))}{|x-y|^{d+2(s-1)}} \dd x\dd y \\
&= C_2 \int_{\real^d}\int_{\real^d} \frac{\big(2u(x)-u(x+y)-u(x-y)\big)\big(2v(x)-v(x+y)-v(x-y)\big)}{|y|^{d+2s}} \dd x\dd y.
\end{aligned}
$$
As observed in \cite[Section~3.1]{abatangelo2018loss}, the bilinear form gives rise to a scalar product which turns $\cH^s_0$ into a Hilbert space, hence we are in the setting of Section~\ref{sec:dual}. Additionally, since $\cH^s_0$ is densely embedded in $\cL^2_0$, we can take $\mathfrak L w = (-\Delta)^s w = B(w,\,\cdot\,)$. The same holds for $s>2$.

\medskip

On the same note we observe that, if the potential $F$ is such that $F \equiv 0$, then the Cahn-Hilliard system formally turns into a system-version of the heat equation for the operator $\mathfrak L(\mathfrak I)$. In the particular case when $\mathfrak L = (-\Delta)^\sigma$ and $\mathfrak I = (-\Delta)^s$, with $\sigma + s > 1$, Theorem~\ref{stm:existenceq} solves the heat equation for the higher order fractional Laplacian $(-\Delta)^{\sigma+s}$ in this system form.

\subsubsection{Sum of operators}
Consider two operators $\mathfrak L_1:X_0^1\to(X_0^1)'$ and $\mathfrak L_1:X_0^2\to(X_0^2)'$ that both fit into the framework of Section~\ref{sec:dual}, where $X_0^1$ and $X_0^2$ are Hilbert spaces densely embedded in $\cL^2_0$. If there exists a continuous embedding $i:X^1_0\hookrightarrow X_0^2$, we can clearly interpret $\mathfrak L_2$ as an operator defined on $X_0^1$ by setting
$$
\tilde{B}_{\mathfrak L_2}:X_0^1\times X_0^1\to\real,\qquad  \tilde{B}_{\mathfrak L_2}(u,v)\coloneqq B_{\mathfrak L_2}(i(u),i(v)).
$$
We can then consider the operator $\mathfrak L\coloneqq\mathfrak L_1+\mathfrak L_2:X_0^1\to(X_0^1)'$ via the associated bilinear form $B_{\mathfrak L}\coloneqq B_{\mathfrak L_1}+\tilde{B}_{\mathfrak L_2}$. Indeed, this map is bounded, since
$$
\begin{aligned}
|B_{\mathfrak L}(u,v)|&\leq |B_{\mathfrak L_1}(u,v)|+|\tilde{B}_{\mathfrak L_2}(u,v)|\leq C_1\|u\|_{0,1}\|v\|_{0,1}+C_2\|i(u)\|_{0,2}\|i(v)\|_{0,2}\\
&\leq C_1\|u\|_{0,1}\|v\|_{0,1}+C_2\|i\|^2_{L(X_0^1,X_0^2)}\|u\|_{0,1}\|v\|_{0,1},
\end{aligned}
$$
and coercive, since
$$
B_{\mathfrak L}(u,u)\geq B_{\mathfrak L_1}(u,u)\geq c_1\|u\|^2_{X_0^1}.
$$

Interesting explicit examples are:
$$
\mathfrak L=-\Delta+(-\Delta)^\sigma:\cH^{\max\{1,\sigma\}}_0\to\Big(\cH^{\max\{1,\sigma\}}_0\Big)'
$$
and
$$
\mathfrak L=(-\Delta)^{\sigma_1}+(-\Delta)^{\sigma_2}:\cH^{\max\{\sigma_1,\sigma_2\}}_0\to\Big(\cH^{\max\{\sigma_1,\sigma_2\}}_0\Big)',
$$
with $\sigma,\sigma_1,\sigma_2\in(0,+\infty)\setminus\mathbb N$.


\section{Regional nonlocal operators}\label{sec:pizzallecipolle}

In this section we address Metatheorem~\ref{stm:Metathm2}, thus considering generalized versions of the following Cahn-Hilliard system
\begin{align}
\partial_t u -\Delta w = 0 &\qquad \text{in } \Omega \times (0,T), \label{eq:CH1reg}\\
w = \mathfrak I_\Omega u + F'(u) &\qquad \text{in } \Omega \times (0,T), \label{eq:CH2reg} \\
u(x,0) = u_0(x) &\qquad \text{in } \Omega, \label{eq:iniCHreg}\\
\frac{\partial w}{\partial\nu}=0
&\qquad \text{on } \partial\Omega \times (0,T). \label{eq:bouCHreg}
\end{align}
Here above, $\mathfrak I_\Omega$ is a variational nonlocal operator of regional type, that can be interpreted in the weak sense of \eqref{eq:pruriginoso}.

\begin{remark}\label{rmk:gut}
Differently from the Dirichlet case, the Neumann boundary condition \eqref{eq:bouCHreg} ensures that the mean of $u$ remains constant in time. Indeed, formally integrating \eqref{eq:CH1reg} yields
$$
\frac{\dd}{\dd t}\int_\Omega u(x,t)\dd x=\int_\Omega\Delta w(x,t)\,dx=\int_{\partial\Omega}\frac{\partial w}{\partial\nu}(x,t)\dd \mathcal H^{d-1}_x=0.
$$
\end{remark}


\subsection{Assumptions and main result}

For the sake of simplicity, in this section we further assume $\Omega\subset\real^d$ to be connected. This hypothesis is needed exclusively to ensure the validity of the fractional Poincar\'e inequality when requiring a bound from below on the kernel $K_\Omega$ only in a neighborhood of the diagonal (see Proposition~\ref{stm:lorenzanonce}). Still, more general situations are allowed, if we require some compatibility between the geometry of $\Omega$ and the choice of interaction radius $\varrho>0$. For the details we refer to Section~\ref{sec:melange}.

Formally, in light of Remark~\ref{rmk:gut}, the time derivative $\partial_t u$ as an element of $(H^1(\Omega))'$ is such that
$$
\{\mbox{constant functions}\}\subset\textup{Ker}(\partial_t u).
$$
This ensures the consistency of the weak formulation of equation \eqref{eq:CH1reg} with respect to the Neumann boundary condition \eqref{eq:bouCHreg} as
$$
\langle \partial_t u, v\rangle_{(H^1(\Omega))'}+\int_\Omega\nabla u\cdot\nabla v\dd x=0,
$$
for every $v\in H^1(\Omega)$.
	
This observation motivates the formulation of a more general Cahn-Hilliard-type system, in which we consider a linear operator $\mathfrak L:X\to X'$ defined on a Hilbert space $X\subset L^2(\Omega)$. We recall that the mean map $\mathfrak m: L^2(\Omega)\to\real$, which is linear, induces the splitting
\begin{equation}\label{eq:umpa}
    L^2(\Omega)=\real\oplus \dot{L}^2(\Omega)=\{\mbox{constant functions}\}\oplus \dot{L}^2(\Omega).
\end{equation}
\begin{hypo}\label{stm:lumpa}
    The precise assumptions are as follows:
    \begin{enumerate}[label=(\roman*)]
        \item $(X_0,(\,\cdot\,,\,\cdot\,)_0)$ is a Hilbert space and $\mathfrak L:X_0\to (X_0)'$ is a linear invertible operator corresponding to a bounded, symmetric and coercive bilinear map $B_\mathfrak L:X_0\times X_0\to\real$, meaning that $\mathfrak L u=B_\mathfrak L(u,\,\cdot\,)$; we refer to Section~\ref{sec:dual} for the details; 
        \item the space $X_0$ is densely embedded in $(\dot{L}^2(\Omega),(\,\cdot\,,\,\cdot\,)_{L^2})$ via the map $i:X_0\hookrightarrow \dot{L}^2(\Omega)$;
	\item as in Section~\ref{sec:dual}, we then consider $X_0$ to be endowed with the inner product
        $$
        (u,v)_{0,\mathfrak L}\coloneqq B_\mathfrak L(u,v);
        $$
        \item we define the Hilbert space $X$ as the direct sum $X\coloneqq \real\oplus X_0$, with the inner product
        $$
        (m_1+u_1,m_2+u_2)_X \coloneqq |\Omega|m_1m_2+(u_1,u_2)_{0,\mathfrak L}\quad\mbox{for every }m_1,m_2\in\real\mbox{ and }u_1,u_2\in X_0;
        $$
	\item we extend $B_\mathfrak L$ to a bilinear map       $B_\mathfrak L:X\times X\to\real$ by setting
        $$
        B_\mathfrak L(m_1+u_1,m_2+u_2)=B_\mathfrak L(u_1,u_2);
        $$
        this is clearly a bounded and symmetric map. This also extends $\mathfrak L$ to the linear operator $\mathfrak L:X\to X'$ given by $\mathfrak L(m+u)=B_\mathfrak L(m+u,\,\cdot\,)$, which is no longer invertible, as indeed
	$$
	\textup{Ker}(\mathfrak L)=X_0^\perp=\real.
	$$
    \end{enumerate}
\end{hypo}
The embedding $i$ extends to an embedding $\hat{\imath}:X\hookrightarrow L^2(\Omega)$ as $m+u\mapsto m +i(u)$. In light of \eqref{eq:umpa} and Hypothesis~\ref{stm:lumpa}\,(ii), the space $X$ is densely embedded in $L^2(\Omega)$, hence we have the triple
$X\hookrightarrow L^2(\Omega)\simeq(L^2(\Omega))'\hookrightarrow X'$.
Recalling Remark~\ref{rmk:triple}, we have
$$
\langle u_1+m_1,u_2+m_2\rangle_{X'}=(\hat{\imath}(u_1+m_1),\hat{\imath}(u_2+m_2))_{L^2}=(i(u_1),i(u_2))_{L^2}+|\Omega|m_1m_2.
$$

With a slight ambiguity, from now on we refer to the factor $\real$ in the decomposition of $X$ as being the constant functions, and to the factor $X_0$ as the zero mass functions.

\begin{remark}\label{rmk:cuneo}
If we choose $\Phi(r)=\frac{1}{2}|r|^2$, the energy functional in \eqref{eq:asti} actually coincides with $\frac{1}{2}\|\,\cdot\,\|^2_{K_\Omega,2}$ and its first variation, as in Lemma~\ref{lem:flareon}, is
\begin{equation*}
\langle \mathfrak{L}w,\psi \rangle_{X'} = B_{\mathfrak{L}}(w,\psi) =  \iint_{\Omega^2} (w(x)-w(y))(\psi(x)-\psi(y)) K_\Omega(x,y) \dd x\dd y.
\end{equation*}
This is a symmetric bilinear form on the Hilbert space $X = \cW^{K_\Omega,2}(\Omega)$, which is bounded, as indeed
$$
|B_{\mathfrak{L}}(w,\psi)| \le 2\|w\|_{K_\Omega,2}\|\psi\|_{K_\Omega,2}.
$$
Moreover, the mass $\mathfrak m : X \to \real$ induces a splitting of $X$ as $X = X_0 \oplus \real = \widehat\cW^{K_\Omega,2}_0(\Omega) \oplus \real$, and $B_{\mathfrak{L}}$ is coercive on the Hilbert space $X_0$ by the fractional Poincar\'e inequality \eqref{eq:eandatavia}, since
$$
\|w\|_X^2 = \|w\|_{L^2(\Omega)}^2 + \|w\|_{K_\Omega,2}^2 \leq (\Lambda\, C+1)\|w\|_{K_\Omega,2}^2 = (\Lambda\, C+1) B_{\mathfrak{L}}(w,w).
$$
Since $X_0$ is densely embedded in $\dot L^2(\Omega)$ by Remark~\ref{rmk:porridge}, and we clearly have
$$
    B_{\mathfrak L}(m_1+w_1,m_2+w_2)=B_{\mathfrak L}(w_1,w_2),
$$
we are then in the setting of Hypothesis~\ref{stm:lumpa}.
\end{remark}

The general Cahn-Hilliard-type system that we consider is then:
\begin{align}
\partial_t u + \mathfrak L w = 0 &\qquad \text{in } \Omega \times (0,T), \label{eq:CH1regL}\\
w = \mathfrak I_\Omega u + F'(u) &\qquad \text{in } \Omega \times (0,T), \label{eq:CH2regL} \\
u(x,0) = u_0(x) &\qquad \text{in } \Omega. \label{eq:iniCHregL}
\end{align}
Here above $\mathfrak L : X \to X'$ satisfies Hypothesis~\ref{stm:lumpa}, while $\mathfrak I_\Omega$ is as in \eqref{eq:pruriginoso}, with the corresponding functional framework of Section~\ref{sec:chem}.
On the other hand, $F'$ formally denotes the derivative of a potential function $F$ satisfying Hypothesis~\ref{stm:hyp}.

\begin{definition}[Solution to the regional Cahn-Hilliard system]\label{def:weaksolregL}
Let $T>0$ be fixed. We say that $(u,w,\zeta)$ is a weak solution to the Cahn-Hilliard system \eqref{eq:CH1regL}--\eqref{eq:CH2regL} associated with the initial datum $u_0 \in \cW^{K_\Omega,q}(\Omega)$ if
$$
\begin{aligned}
u &\in L^\infty(0,T;\cW^{K_\Omega,q}(\Omega)) \cap W^{1,2}(0,T;X'), \\
w &\in L^2(0,T;X),\\
\zeta &\in L^2(0,T;L^2(\Omega)), \qquad \zeta\in\gamma(u) \quad \text{a.\,e. in } \Omega \times (0,T),
\end{aligned}
$$
and $(u,w,\zeta)$ satisfies the following weak formulation of \eqref{eq:CH1regL}--\eqref{eq:CH2regL}:
\begin{align}
\partial_t u + \mathfrak L w = 0 &\qquad \text{in } X', \label{eq:CHw1regL}\\
w = \mathfrak I_\Omega u + \zeta + \pi(u) &\qquad \text{in } (\cW^{K_\Omega,q}(\Omega))', \label{eq:CHw2regL}
\end{align}
almost everywhere in $(0,T)$, with $u(0) = u_0$ almost everywhere in $\Omega$.
\end{definition}

By Hypothesis~\ref{stm:lumpa} and Remark~\ref{rmk:porridge}, we have the triples
$$
X\hookrightarrow L^2(\Omega)\simeq (L^2(\Omega))'\hookrightarrow X'
$$
and
$$
\cW^{K_\Omega,q}(\Omega)\hookrightarrow L^2(\Omega)\simeq (L^2(\Omega))'\hookrightarrow (\cW^{K_\Omega,q}(\Omega))'.
$$

\begin{theorem}[Existence and uniqueness]\label{stm:existenceregL}
Let Hypothesis~\ref{stm:hyp} be satisfied. Assume
$$
u_0 \in \cW^{K_\Omega,q}(\Omega), \quad \Gamma(u_0) \in L^1(\Omega)\quad\mbox{and}\quad\mathfrak m(u_0)\in \textup{Int}\,D(\gamma).
$$
Then the Cahn-Hilliard system \eqref{eq:CH1regL}--\eqref{eq:iniCHregL} admits a solution $(u,w,\zeta)$, in the sense of Definition~\ref{def:weaksolregL} for every $T>0$. Moreover,
\begin{equation}\label{eq:starmie}
\mathfrak m(u(t))=\mathfrak m(u_0)\quad\mbox{for every }t\geq 0,
\end{equation}
and, defining the energy
$$
\mathscr E_\Omega(u) \coloneqq \iint_{\Omega^2}\Phi(u(x)-u(y))K_\Omega(x,y)\dd x\dd y + \io F(u(x))\dd x,
$$
the following energy estimate holds:
\begin{equation}\label{eq:en-estregL}
\frac{1}{2}\int_0^t \|w(\tau)-\mathfrak m(w(\tau))\|_{X_0}^2 \dd \tau + \mathscr E_\Omega(u(t)) \le \mathscr E_\Omega(u_0)\qquad\mbox{for every }t\geq0.
\end{equation}
If $\phi$ satisfies condition \eqref{eq:nameless1}, then $u$ is unique, and
\begin{equation}\label{eq:est-ensregL}
\mathscr E_\Omega(u(t)) \le \mathscr E_\Omega(u(\tau)) \qquad \mbox{for } \ 0 \le \tau \le t.
\end{equation}
If we also assume $\Gamma$ to be differentiable, then also $w$ and $\zeta$ are unique.
On the other hand, if there is an element of $\textup{Int}\,D(\gamma)$ at which $\gamma$ is multivalued, then there exists an initial datum for which the solution is not unique.
\end{theorem}

\begin{remark}\label{rmk:stj}
Some comments are in order:
    \begin{enumerate}
    \item the initial datum $u_0$ can actually be chosen so that $\mathfrak m(u_0)\in D(\gamma)\setminus\textup{Int}\,D(\gamma)$. However, if we also require $\Gamma(u_0)\in L^1(\Omega)$, then we must have $u_0\equiv m$ and it is trivial to find $w$ and $\zeta$ such that the triple $(m,w,\zeta)$ is a solution.
    \item Thanks to Lemma~\ref{lem:subdiff}\,(i) and Lemma~\ref{lem:viabrombeis}, we know that $\Gamma$ is differentiable at $r_0\in\textup{Int}\, D(\Gamma)$ if and only if $\gamma(r_0)$ is a singleton.
    \item Concerning the uniqueness of $w$ and $\zeta$, the following more general claims hold true. If $u$ is unique, then also the $X_0$-component of $w$ is unique. On the other hand the $\real$-component, i.\,e., the mass $\mathfrak{m}(\hat{\imath}(w))$, might not be unique.
    However, if
    $$
    u(x,t)\in \textup{Int}\,D(\Gamma) \quad\mbox{for almost every }x\in\Omega\mbox{ and }t\geq0,
    $$
    and $\Gamma$ is differentiable in $\textup{Int}\,D(\Gamma)$, then also $\zeta$, hence $\mathfrak{m}(\hat{\imath}(w))$ and $w$, are unique.
    \end{enumerate}
\end{remark}

\subsubsection{Proof of Theorem~\ref{stm:existenceregL}}\label{sec:ProofThmRegExi}
We split the proof into several steps. We begin by establishing the conservation of mass \eqref{eq:starmie}, thus justifying the formal computation in Remark~\ref{rmk:gut}. This motivates us to seek the solution $u$ within the space of functions of $\cW^{K_\Omega,q}(\Omega)$ having fixed mass, equal to that of the initial datum $u_0$. We are then led to a solution to
\begin{equation}\label{eq:hangedsponge}
w=\mathfrak I_\Omega u+\zeta+\pi(u)\qquad\mbox{in }(\widehat{\cW}^{K_\Omega,q}_0(\Omega))'.
\end{equation}
Moreover, since the mass of $u$ is constant in time, we know that
$$
\{\mbox{constant functions}\}\subset\textup{Ker}(\partial_t u),
$$
hence solving equation \eqref{eq:CHw1regL} is equivalent to solving
\begin{equation}\label{eq:espeon}
\partial_tu+\mathfrak L w=0
\qquad\mbox{in }(X_0)'.
\end{equation}
We also point out that, if $w$ solves \eqref{eq:espeon} for almost every $t\in(0,T)$, then also any function of the form $w+\vartheta$, with $\vartheta$ depending only on $t$, still solves the above equation.
We thus restrict ourselves to considering functions $w\in X_0$.
Since
$$
\{\mbox{constant functions}\}\subset\textup{Ker}(\mathfrak I_\Omega u),
$$
in order to pass from equation \eqref{eq:hangedsponge} to \eqref{eq:CHw2regL}, it is enough to consider $w+\mathfrak m(\zeta+\pi(u))$ in place of $w$. 
Moreover, as observed above, the function $w+\mathfrak m(\zeta+\pi(u))$ also solves \eqref{eq:espeon}, hence \eqref{eq:CHw1regL}.

Notice also that by iteratively gluing together solutions, as in \eqref{eq:superattak}, we obtain a global in time solution, once we have proved the finite-time existence. Hence, in what follows, we restrict ourselves to a fixed arbitrary time interval $(0,T)$.

\textbf{Step 1: conservation of mass.}
Let $(u,w,\zeta)$ be a solution in the sense of Definition~\ref{def:weaksolregL}.
Then we have
$$
\int_\Omega u(x,\tau)\dd x- \int_\Omega u(x,0)\dd x=\int_0^\tau\langle \partial_t u,1\rangle_{X'}\dd t=-\int_0^\tau\langle\mathfrak L w,1\rangle_{X'}\dd t=-\int_0^\tau B_\mathfrak L (w,1)\dd t=0,
$$
which proves \eqref{eq:starmie}.

\textbf{Step 2: time-discretization.}
We begin by fixing $T>0$ and we consider $N\in\nat$ and $\tau=\tau_N \coloneqq T/N$.
Then, we define the regularized energy
$$
\mathscr E_\Omega^\lambda(u)\coloneqq\mathfrak F_\Omega(u)+\io\big(\Gamma_\lambda(u)+\Pi(u)\big)\dd x
$$
for $\lambda\in(0,1)$, with $\mathfrak F_\Omega$ as defined in \eqref{eq:asti}, and we prove the analogues of Lemma~\ref{lemma:neutralquokkaq} and Lemma~\ref{lemma:emmaq}.

\begin{lemma}\label{lemma:neutralquokkaregL}
The functional $\mathscr E^\lambda$ is lower semicontinuous with respect to the $L^2$-convergence. Moreover
	\begin{equation}\label{eq:gargoyleregL}
		\mathscr E^\lambda_\Omega(u) \geq \frac{1}{2\Lambda q}\|u\|_{K_\Omega,q}^q-c_\natural|\mathfrak m(u)|^q-\beta_0|\Omega|\qquad\mbox{for every }u\in\cW^{K_\Omega,q}(\Omega),
	\end{equation}
	for $\lambda\in(0,\bar\lambda)$, 
	where $c_\natural,\beta_0>0$ and $\bar\lambda\in(0,1)$ depend only on $d,\Lambda,s,\varrho,q,\Omega,a_1,a_2$ and $C_\pi$.
\end{lemma}

\begin{proof}
The semicontinuity follows by arguing as in the proof of \eqref{eq:pandasaviour}.
The estimate \eqref{eq:gargoyleregL} follows by using the fractional Poincar\'e inequality \eqref{eq:eevee} and the fact that $q\geq2$, as in the proof of \eqref{eq:gargoyleq}.
\end{proof}

Lemma~\ref{lemma:neutralquokkaregL} provides one of the two fundamental ingredients of the Direct Method of the Calculus of Variations, namely the lower semicontinuity of the functional. In light of \eqref{eq:gargoyleregL} and the fractional Poincar\'e inequality \eqref{eq:eandatavia}, in order to ensure the compactness of a minimizing sequence, it is enough to fix the mass of the functions we work with. We then obtain the following result.

\begin{lemma}\label{lemma:emmaregL}
Given $\lambda\in(0,\bar\lambda)$, $\tau>0$ and any $g\in L^2(\Omega)$, we set $m\coloneqq\mathfrak m(g)$ and define the energy $\mathscr E_{\Omega,g}^\lambda:\widehat{\cW}^{K_\Omega,q}_m(\Omega)\to\real$ by setting
	$$
	\mathscr E_{\Omega,g}^\lambda(u) \coloneqq \frac{1}{2\tau}\|u-g\|^2_{(X_0)'} + \mathscr E_\Omega^\lambda(u).
	$$
	Then, there exists at least one function $u_\star\in\widehat{\cW}^{K_\Omega,q}_m(\Omega)$ such that
	$$
	\mathscr E_{\Omega,g}^\lambda(u_\star) = \inf_{u\in\widehat{\cW}^{K_\Omega,q}_m(\Omega)}\mathscr E_{\Omega,g}^\lambda(u).
	$$
 As a consequence
 \begin{equation}\label{eq:vaporeon}
 \mathfrak L^{-1}\bigg(\frac{u_\star-g}{\tau}\bigg)+\mathfrak I_\Omega u_\star+\gamma_\lambda(u_\star)+\pi(u_\star)=0\quad\mbox{in }\big(\widehat{\cW}^{K_\Omega,q}_0(\Omega)\big)'.
 \end{equation}
\end{lemma}

\begin{proof}
By \eqref{eq:gargoyleregL} the infimum is finite. Thus, if we consider a minimizing sequence, i.\,e., $\{v_k\}\subset\widehat{\cW}^{K_\Omega,q}_m(\Omega)$ such that
$$
\lim_{k\to\infty}\mathscr E^\lambda_{\Omega,g}(v_k) = \inf_{u\in\widehat{\cW}^{K_\Omega,q}_m(\Omega)}\mathscr E_{\Omega,g}^\lambda(u),
$$
then, by \eqref{eq:gargoyleregL} and the fractional Poincar\'e inequality \eqref{eq:eevee}, we have
$$
\begin{aligned}
&\|v_k\|_{K_\Omega,q}^q\leq 2\Lambda q\bigg(\beta_0|\Omega|+c_\natural|m|^q+\inf_{u\in\widehat{\cW}^{K_\Omega,q}_m(\Omega)}\mathscr E_{\Omega,g}^\lambda(u)+1\bigg),\\
&
\|v_k\|_{L^q(\Omega)}^q\leq C\big(\|v_k\|_{K_\Omega,q}^q+|m|^q\big)
\end{aligned}
$$
for every $k$ big enough. This implies that $\{v_k\}$ is bounded in $\cW^{K_\Omega,q}(\Omega)$ which is compactly embedded in $L^2(\Omega)$, hence there exists $u_\star\in \widehat{\cW}^{K_\Omega,q}_m(\Omega)$ such that
$$
v_k\to u_\star\qquad\mbox{strongly in }L^2(\Omega),
$$
up to a subsequence, that we do not relabel. By the $L^2$-convergence and the lower semicontinuity of $\mathscr E^\lambda_{\Omega}$ proved in Lemma~\ref{lemma:neutralquokkaregL},
$$
\mathscr E_{\Omega,g}^\lambda(u_\star)\leq\liminf_{k\to\infty}\mathscr E^\lambda_{\Omega,g}(v_k)=
\inf_{u\in \widehat{\cW}^{K_\Omega,q}_m(\Omega)}\mathscr E_{\Omega,g}^\lambda(u),
$$
proving that $u_\star$ is indeed a minimum.
To conclude, notice that if $v\in\widehat{\cW}^{K_\Omega,q}_0(\Omega)$, then $u_\star+\varepsilon v\in \widehat{\cW}^{K_\Omega,q}_m(\Omega)$ for every $\varepsilon\in\real$. Therefore, by the minimality of $u_\star$ and recalling Lemma~\ref{lem:flareon} and that $\mathfrak L^{-1}$ is the Fr\'echet derivative of $\|\,\cdot\,\|_{(X_0)'}^2/2$, we obtain \eqref{eq:vaporeon}.
\end{proof}

As anticipated, we fix $m\coloneqq \mathfrak m(u_0)$ and we carry out the argument within the space $\widehat{\cW}^{K_\Omega,q}_m(\Omega)$.

\medskip

We first apply Lemma~\ref{lemma:emmaregL} to pick a minimizer $u^\lambda_1\in\widehat{\cW}^{K_\Omega,q}_m(\Omega)$ of the functional
$$
\mathscr E_1^\lambda(u) \coloneqq \frac{\tau}{2}\left\|\frac{u-u_0}{\tau}\right\|^2_{(X_0)'} + \mathscr E^\lambda_\Omega(u)=\mathscr E_{\Omega,u_0}^\lambda(u).
$$
Then we proceed iteratively, choosing $u_n^\lambda\in\widehat{\cW}^{K,q}_m(\Omega)$ to be a minimizer of
$$
\mathscr E_n^\lambda(u) \coloneqq \frac{\tau}{2}\left\|\frac{u-u_{n-1}}{\tau}\right\|^2_{(X_0)'} + \mathscr E^\lambda_\Omega(u).
$$
We also define $\omega_n^\lambda\in X_0$ as
$$
\omega_n^\lambda\coloneqq(-\mathfrak L)^{-1}\bigg(\frac{u_n^\lambda-u_{n-1}^\lambda}{\tau}\bigg).
$$
Recalling \eqref{eq:vaporeon}, we find that the pair $(u_n^\lambda,\omega_n^\lambda)$ satisfies
\begin{align*}
\frac{u^\lambda_n - u^\lambda_{n-1}}{\tau} + \mathfrak L \, \omega^\lambda_n = 0 &\qquad \text{in } (X_0)',\\
\omega^\lambda_n = \mathfrak I_\Omega u^\lambda_n + \gamma_\lambda(u^\lambda_n) + \pi(u^\lambda_n) &\qquad \text{in } (\widehat{\cW}^{K_\Omega,q}_0(\Omega))',
\end{align*}
for every $n=1,\ldots,N$. Moreover, we obtain the counterpart of the energy estimate \eqref{eq:underclasshero}, that is
\begin{equation}\label{eq:somesayregL}
\frac{\tau}{2}\sum_{n=1}^{\bar{n}} \left\|\frac{u_n^\lambda - u^\lambda_{n-1}}{\tau}\right\|_{(X_0)'}^2 + \mathscr E^\lambda_\Omega(u_{\bar{n}}^\lambda)\le \mathscr E_\Omega(u_0),
\end{equation}
and we recall that, by \eqref{eq:lorenzace}, we have
$$
\|\omega_n^\lambda\|_{X_0}=\left\|\frac{u_n^\lambda - u^\lambda_{n-1}}{\tau}\right\|_{(X_0)'}.
$$

\textbf{Step 3: convergence as $\tau\to0$.}
We then define the interpolating functions as in \eqref{eq:inter}. Notice that
$$
\mathfrak m(\bar{u}^\lambda_\tau(t))=\mathfrak m(\hat{u}^\lambda_\tau(t))=m
$$
for every $\lambda,\tau,t$. Moreover, $\partial_t \hat{u}^\lambda_\tau\in\cW^{K_\Omega,q}(\Omega)$ with $\mathfrak m(\partial_t \hat{u}^\lambda_\tau)=0$, hence, if we interpret it as an element of $X'$ via $\cW^{K_\Omega,q}(\Omega)\hookrightarrow L^2(\Omega)\simeq(L^2(\Omega))'\hookrightarrow X'$, we have
$$
\{\mbox{constant functions}\}\subset\textup{Ker}(\partial_t \hat{u}^\lambda_\tau).
$$
As a consequence
$$
\|\partial_t \hat{u}^\lambda_\tau\|_{X'}=\|\partial_t \hat{u}^\lambda_\tau\|_{(X_0)'}.
$$
Therefore, from \eqref{eq:gargoyleregL} and \eqref{eq:somesayregL} we deduce the analogues of \eqref{eq:esti1_tau} and \eqref{eq:esti2_tau},
\begin{equation}\label{eq:esti1_tauregL}
\begin{aligned}
&\int_0^T \|\bar{\omega}^\lambda_\tau(t)\|_{X_0}^2 \dd t + \sup_{t \in [0,T]} \|\bar{u}^\lambda_\tau(t)\|_{K_\Omega,q}^q
\le C\big(\mathscr E_\Omega(u_0)+|\mathfrak m(u_0)|^q+1\big),\\
&
\int_0^T \|\partial_t\hat{u}^\lambda_\tau(t)\|_{X'}^2 \dd t \le C\big(\mathscr E_\Omega(u_0)+|\mathfrak m(u_0)|^q+1\big),\\
&
\sup_{t \in [0,T]} \|\hat{u}^\lambda_\tau(t)\|_{K_\Omega,q}^q \le C\big(\mathscr E_\Omega(u_0)+|\mathfrak m(u_0)|^q+1\big),
\end{aligned}
\end{equation}
where $C$ is a constant not depending on $\tau$ nor $\lambda$.

These estimates and the fractional Poincar\'e inequality \eqref{eq:eevee} lead to the convergence of the interpolants as $\tau\to0$, as in \eqref{eq:conv} and \eqref{eq:conv_ALS}, via compactness, to a pair $(u^\lambda,\omega^\lambda)$ such that
$$
u^\lambda \in L^\infty(0,T;\cW^{K_\Omega,q}(\Omega))\cap W^{1,2}(0,T;X'), \quad \omega^\lambda \in L^2(0,T;X_0).
$$
More precisely,
\begin{equation*}
\begin{array}{rcll}
\bar{\omega}^\lambda_\tau &\rightharpoonup& \omega^\lambda &\quad \text{weakly in } L^2(0,T;X_0), \\[2pt]
\hat{u}^\lambda_\tau &\to& u^\lambda &\quad \text{strongly in } C([0,T];L^2(\Omega)), \\[2pt]
 \bar{u}^\lambda_\tau &\to& u^\lambda &\quad \text{strongly in } L^\infty(0,T;L^2(\Omega)), \\[2pt]
\partial_t\hat{u}^\lambda_\tau &\rightharpoonup& \partial_t u^\lambda &\quad \text{weakly in } L^2(0,T;X').
\end{array}
\end{equation*}
Observe also that $\mathfrak m(u^\lambda(t))=m$ for every $\lambda$ and $t$.
The convergences only hold true for a specific subsequence $\{\tau_{N_h}\}$.

For the convergence of $\mathfrak I_\Omega \bar{u}^\lambda_\tau$ we repeat the argument leading to \eqref{eq:slowpoke}, thus obtaining 
$$
\int_0^T\vartheta(t)\langle\mathfrak I_\Omega \bar{u}_{\tau_{N_h}}^\lambda(t),v\rangle_{K_\Omega,q}\dd t
\xrightarrow{h\to\infty}\int_0^T\vartheta(t)\langle\mathfrak I_\Omega u^\lambda(t),v\rangle_{K_\Omega,q}\dd t,
$$
for every $\vartheta\in L^q(0,T)$ and $v\in\cW^{K_\Omega,q}(\Omega)$.
Hence, by passing to the limit $h\to\infty$, the limit functions $(u^\lambda,\omega^\lambda)$ satisfy the system
\begin{equation}\label{eq:drowzee}\begin{aligned}
&\int_0^T \langle  \partial_t u^\lambda(t)
+\mathfrak L\, \omega^\lambda(t),v_1\rangle_{(X_0)'} \,\vartheta(t) \dd t = 0,\\
&
\int_0^T \left(\langle\mathfrak I_\Omega u^\lambda(t),v_2\rangle_{K_\Omega,q} + \io \big(\gamma_\lambda(u^\lambda(x,t))+\pi(u^\lambda(x,t))-\omega^\lambda(x,t)\big)v_2(x) \dd x\right)\vartheta(t) \dd t=0,
\end{aligned}\end{equation}
for every $\vartheta\in L^q(0,T)$, $v_1\in X_0$ and $v_2\in\widehat{\cW}^{K_\Omega,q}_0(\Omega)$.

\textbf{Step 4: convergence as $\lambda\to0$.}
We define the function $w^\lambda\in L^2(0,T;X)$ as
$$
w^\lambda(x,t)\coloneqq\omega^\lambda(x,t)+\mathfrak m\big(\gamma_\lambda(u^\lambda(t))+\pi(u^\lambda(t))\big),
 $$
and we observe that by \eqref{eq:drowzee} the pair $(u^\lambda, w^\lambda)$ satisfies
 \begin{align}
\partial_t u^\lambda + \mathfrak L w^\lambda = 0 &\qquad \text{in } X',\label{eq:CHw1regLl} \\
w^\lambda = \mathfrak I_\Omega u^\lambda + \gamma_\lambda(u^\lambda) + \pi(u^\lambda) &\qquad \text{in } (\cW^{K_\Omega,q}(\Omega))',\label{eq:CHw2regLl}
\end{align}
almost everywhere in $(0,T)$.

Notice that estimates \eqref{eq:esti1_tauregL} are preserved when passing to the limit $\tau\to0$, yielding
\begin{equation*}
\begin{aligned}
&\int_0^T \|\omega^\lambda(t)\|_{X_0}^2 \dd t 
\le C\big(\mathscr E_\Omega(u_0)+|\mathfrak m(u_0)|^q+1\big),\\
&
\int_0^T \|\partial_t u^\lambda(t)\|_{X'}^2 \dd t \le C\big(\mathscr E_\Omega(u_0)+|\mathfrak m(u_0)|^q+1\big),\\
&
\sup_{t \in [0,T]} \|u^\lambda(t)\|_{K_\Omega,q}^q \le C\big(\mathscr E_\Omega(u_0)+|\mathfrak m(u_0)|^q+1\big).
\end{aligned}
\end{equation*}
Recalling Hypothesis~\ref{stm:lumpa}\,(iv), we observe that
\begin{equation}\label{eq:saintJimmy}
\|w^\lambda(t)\|_X^2=|\Omega|\big|\mathfrak m\big(\gamma_\lambda(u^\lambda(t))+\pi(u^\lambda(t))\big)\big|^2+\|\omega^\lambda(t)\|_{X_0}^2.
\end{equation}
Since $\pi$ is Lipschitz, by Jensen's inequality we have
$$
\begin{aligned}
|\mathfrak m\big(\pi(u^\lambda(t))\big)|^2&\leq |\Omega|^{-2}\|\pi(u^\lambda(t))\|^2_{L^2(\Omega)}\leq |\Omega|^{-2}C_\pi^2\|u^\lambda(t)\|^2_{L^2(\Omega)}\\
&
\leq|\Omega|^{-\frac{q+2}{q}}C_\pi^2\|u^\lambda(t)\|_{L^q(\Omega)}^2.
\end{aligned}
$$
Recalling the fractional Poincar\'e inequality \eqref{eq:eandatavia} and the fact that $q\geq 2$, we can estimate
$$
\begin{aligned}
\|u^\lambda(t)\|_{L^q(\Omega)}^2&\leq\Big(C\Lambda\|u^\lambda(t)\|_{K_\Omega,q}+|\Omega|^\frac{1}{q}|\mathfrak m(u^\lambda(t))|\Big)^2\leq C\Big(\|u^\lambda(t)\|_{K_\Omega,q}^2+|\mathfrak m(u_0)|^2\Big)\\
&
\leq
C\Big(\big(\mathscr E_\Omega(u_0)+|\mathfrak m(u_0)|^q+1\big)^\frac{2}{q}+|\mathfrak m(u_0)|^2\Big)\\
&
\leq
C\Big(\mathscr E_\Omega(u_0)+|\mathfrak m(u_0)|^q+1\Big),
\end{aligned}
$$
hence
$$
|\mathfrak m\big(\pi(u^\lambda(t))\big)|^2
\leq
C\Big(\mathscr E_\Omega(u_0)+|\mathfrak m(u_0)|^q+1\Big)
$$
uniformly in $\lambda$ and $t$.

Coming back to \eqref{eq:saintJimmy}, we are left to estimate the mass of $\gamma_\lambda(u^\lambda)$.
To this aim, we test equation \eqref{eq:CHw2regLl} with $u^\lambda(t)-m\in\widehat{\cW}^{K_\Omega,q}_0(\Omega)$. Recalling that $\phi(r)r\geq0$, we have
$$
\langle\mathfrak I_\Omega u^\lambda(t),u^\lambda(t)-m\rangle_{K_\Omega,q}
=\langle\mathfrak I_\Omega u^\lambda(t),u^\lambda(t)\rangle_{K_\Omega,q}\geq0,
$$
hence we obtain
\begin{equation}\label{eq:trapattoni}
\begin{aligned}
\int_\Omega \gamma_\lambda(u^\lambda(x,t))\big(u^\lambda(x,t)-m\big)\dd x&\leq \int_\Omega \big(\omega^\lambda(x,t)-\pi(u^\lambda(x,t))\big)\big(u^\lambda(x,t)-m\big)\dd x\\
&
\leq\int_\Omega\left(|\omega^\lambda(x,t)|^2+C_\pi^2|u^\lambda(x,t)|^2+\frac{1}{2}|u^\lambda(x,t)-m|^2\right)\dd x.
\end{aligned}
\end{equation}
On the other hand, we can estimate
\begin{equation}\label{eq:P2}
\int_\Omega \gamma_\lambda(u^\lambda(x,t))\big(u^\lambda(x,t)-m\big)\dd x
\geq\delta\int_\Omega|\gamma_\lambda(u^\lambda(x,t))|\dd x-C(|\Omega|,m,\delta),
\end{equation}
where $\delta>0$ is such that $[m-\delta,m+\delta]\subset \textup{Int }D(\gamma)$. This classical estimate is a consequence of the monotonicity of $\gamma$ and the properties of the Moreau-Yosida regularizations and can be found, e.\,g., in \cite[Section~5]{kenmochi1995subdifferential} (together with the explicit expression of the constant $C$).

Combining \eqref{eq:trapattoni} and \eqref{eq:P2}, integrating in time, and exploiting the uniform estimates obtained above, we get
$$\begin{aligned}
\int_0^T&\|\gamma_\lambda(u^\lambda(t))\|_{L^1(\Omega)}\dd t\\
&\leq \frac{C}{\delta}\int_0^T \left(\|\omega^\lambda(t)\|^2_{X_0}+\|u^\lambda(t)\|^2_{L^2(\Omega)}+\|u^\lambda(t)-m\|^2_{L^2(\Omega)}+C(|\Omega|,m,\delta)\right)\dd t\leq C,
\end{aligned}$$
uniformly in $\lambda$.
Then, as a consequence of \eqref{eq:saintJimmy}, we have
$$
\int_0^T \|w^\lambda(t)\|_X^2\dd t\leq C,
$$
uniformly in $\lambda$. We can now proceed to estimate the norm of $\gamma_\lambda(u^\lambda)$ in $L^2(0,T;L^2(\Omega))$. For this, as in previous sections, we test equation \eqref{eq:CHw2regLl} with $\gamma_\lambda(u^\lambda)\in\cW^{K_\Omega,q}(\Omega)$ and we integrate in time. Recalling that $\gamma_\lambda$ is increasing and $\phi(r)r\geq0$, we have
$$
\langle\mathfrak I_\Omega u^\lambda(t),\gamma_\lambda(u^\lambda(t))\rangle_{K_\Omega,q}\geq 0,
$$
hence we obtain
$$
\int_0^T\int_\Omega|\gamma_\lambda(u^\lambda(x,t))|^2\dd x\dd t\leq \int_0^T\int_\Omega\big(w^\lambda(x,t)-\pi(u^\lambda(x,t))\big)\gamma_\lambda(u^\lambda(x,t))\dd x\dd t.
$$
Thus, by Young's inequality
$$
\frac{1}{2}
\int_0^T\|\gamma_\lambda(u^\lambda(t))\|^2_{L^2(\Omega)}\dd t
\leq\int_0^T \Big(\|w^\lambda(t)\|_{L^2(\Omega)}^2+C_\pi^2\|u^\lambda(t)\|_{L^2(\Omega)}^2\Big)\dd t\leq C,
$$
uniformly in $\lambda$.

Together with the fractional Poincar\'e inequality \eqref{eq:eevee} and the fact that $\mathfrak m(u^\lambda(t))=m$ for every $\lambda$ and $t$, these estimates are enough to ensure the passage to the limit of $u^\lambda$, $\gamma_\lambda(u^\lambda)$ and $w^\lambda$, for a certain subsequence $\lambda_k\to0$, obtaining
\begin{equation*}
\begin{array}{rcll}
	w^{\lambda_k} &\rightharpoonup& w &\quad \text{weakly in } L^2(0,T;X), \\[2pt]
	u^{\lambda_k} &\to& u &\quad \text{strongly in } C([0,T];L^2(\Omega)), \\[2pt]
	\partial_t u^{\lambda_k} &\rightharpoonup& \partial_t u &\quad \text{weakly in } L^2(0,T;X'),\\[2pt]
    \gamma_{\lambda_k}(u^{\lambda_k}) &\rightharpoonup& \zeta &\quad \text{weakly in } L^2(0,T;L^2(\Omega)).
\end{array}
\end{equation*}
We have also the analogue of \eqref{eq:repetto}, that is
\begin{equation*}
\phi(u^{\lambda_k}(x,t)-u^{\lambda_k}(y,t))K_\Omega(x,y)^\frac{1}{q'}\rightharpoonup \phi(u(x,t)-u(y,t))K_\Omega(x,y)^\frac{1}{q'}\quad\mbox{weakly in }L^{q'}(\Omega^2\times(0,T)).
\end{equation*}
We can thus pass to the limit $\lambda_k\to0$ in \eqref{eq:drowzee} and use the fundamental lemma of the calculus of variations for Bochner spaces, to obtain \eqref{eq:CHw1regL}--\eqref{eq:CHw2regL} for almost every $t\in(0,T)$, as wanted.
Moreover we point out that, arguing as in the end of Step~7 of Section~\ref{sec:exince}, ensures that $\zeta(x,t)\in\gamma(u(x,t))$ for almost every $(x,t)\in\Omega\times(0,T)$.

Observing that, with a slight abuse of notation,
$$
w(t)-\mathfrak m(w(t))=-\mathfrak L^{-1}\frac{\dd}{\dd t}u(t),
$$
by the discrete energy estimate \eqref{eq:somesayregL} and \eqref{eq:lorenzace}, we derive \eqref{eq:en-estregL} arguing as in Step~8 of Section~\ref{sec:exince}.

\textbf{Step 5: uniqueness.}
For the sake of clarity, in this step of the argument it is convenient to maintain the distinction between $w$ and $\hat{\imath}(w)$. Let $(u_1,w_1,\zeta_1)$ and $(u_2,w_2,\zeta_2)$ be two solutions in the sense of Definition~\ref{def:weaksolregL}, corresponding to the same initial datum $u_0$. Let $\wp_1:X\to\real$ and $\wp_2:X\to X_0$ be the projection maps, and notice that $\wp_1=\mathfrak{m}\circ\hat{\imath}$. We denote $$
\tilde{u}=u_1-u_2,\;\;\tilde{w}=\wp_2(w_1-w_2)\ \mbox{ and } \ m=\wp_1(w_1-w_2),
$$
hence $w_1-w_2=\tilde{w}+m$. The conservation of mass \eqref{eq:starmie} ensures that $\tilde{u}(t)\in\widehat{\cW}^{K_\Omega,q}_0(\Omega)$, and clearly $\tilde{w}(t)\in X_0$ by definition, for every $t\geq0$.
Then
$$
0=\partial_t \tilde{u}(t)+\mathfrak L(\tilde{w}(t)+m(t))=\partial_t \tilde{u}(t)+\mathfrak L\tilde{w}(t)\textrm{ in }X'$$
for almost every $t\in(0,T)$, and, since $\mathfrak L:X_0\to (X_0)'$ is invertible, we also have
\begin{equation}\label{eq:tameimpala}
\mathfrak L^{-1}\frac{\dd}{\dd t}\tilde{u}(t)+\tilde{w}(t)=0\mbox{ in }X_0.
\end{equation}
On the other hand
$$
\tilde{w}+m=\mathfrak I_\Omega u_1-\mathfrak I_\Omega u_2+\zeta_1-\zeta_2+\pi(u_1)-\pi(u_2)\mbox{ in }(\cW^{K_\Omega,q})',
$$
hence, if we restrict ourselves to test functions in $\widehat{\cW}_0^{K_\Omega,q}$, we obtain
$$
\tilde{w}=\mathfrak I_\Omega u_1-\mathfrak I_\Omega u_2+\zeta_1-\zeta_2+\pi(u_1)-\pi(u_2)\mbox{ in }(\widehat{\cW}_0^{K_\Omega,q})'.
$$
Condition \eqref{eq:nameless1} implies the following analogue of
\eqref{eq:Okinawa},
$$
\frac{1}{\Lambda}\|v_1-v_2\|_{K_\Omega,q}^q\leq\langle\mathfrak I_\Omega v_1- \mathfrak I_\Omega v_2,v_1-v_2\rangle_{K_\Omega,q}
$$
for every $v_1,v_2\in\cW^{K_\Omega,q}(\Omega)$.
Exploiting Ehrling's lemma for the embeddings $\widehat{\cW}_0^{K_\Omega,q} \hookrightarrow \dot{L}^2(\Omega) \hookrightarrow (X_0)'$ and arguing as in Step~1 of Section~\ref{sec:exince}, we obtain the following counterpart of
\eqref{eq:fantasmaformaggino}
$$
\frac{\dd}{\dd t}\|u_1-u_2\|^2_{(X_0)'}+\|u_1-u_2\|^q_{K_\Omega,q}\leq C \|u_1-u_2\|^2_{(X_0)'},
$$
for almost every $t\in(0,T)$.
This is enough to conclude the desired uniqueness of $u$, via Gr\"onwall's inequality.
By \eqref{eq:tameimpala} we have $\tilde{w}=0$, hence $\wp_2(w_1)=\wp_2(w_2)$, so that $w_1$ and $w_2$ are equal up to an additive constant (which may vary in $t$).

\medskip

If $u$ is unique, we can then obtain the monotonicity of the energy \eqref{eq:est-ensregL} by arguing as in Step~8 of Section~\ref{sec:exince}. Indeed, the possible non-uniqueness of $w$ and $\zeta$ does not affect the conclusion.

\medskip

On the other hand, concerning $w$ and $\zeta$ we have the following results.

If $\Gamma$ is differentiable, then
$$
\zeta_1(x,t)=\gamma(u(x,t))=\zeta_2(x,t)
$$
for almost every $x\in\Omega$ and $t\geq0$. As a consequence, considering $v\equiv1\in\cW^{K_\Omega,q}(\Omega)$ as test function in \eqref{eq:CHw2regL} we obtain
$$
|\Omega|\wp_1(w_1)=\int_\Omega \hat{\imath}(w_1)(x,t)\dd x=\int_\Omega \big(\zeta(x,t)+\pi(u(x,t))\big)\dd x =\int_\Omega \hat{\imath}(w_2)(x,t)\dd x=|\Omega|\wp_1(w_2)
$$
for every $t\geq0$, proving that also $w$ is unique.

Suppose instead that there exists $m_0\in\textup{Int}\,D(\gamma)$ at which $\gamma$ is multivalued and let $\xi_1\not=\xi_2\in\gamma(m_0)$. We choose $u_0\equiv m_0$ as initial datum. Then setting $u_1\equiv m_0\equiv u_2$, $\zeta_1\equiv\xi_1$, $\zeta_2\equiv \xi_2$, and $w_1\equiv\xi_1+\pi(m_0)$, $w_2\equiv\xi_2+\pi(m_0)$,
the triples $(u_1,w_1,\zeta_1)$ and $(u_2,w_2,\zeta_2)$ are two distinct solutions corresponding to the same initial datum $u_0$.


\subsection{Examples of regional kernels}\label{sec:kernels}

Some meaningful examples of regional kernels satisfying \eqref{eq:uhu} and either \eqref{eq:sticweak} or \eqref{eq:stic}, thus representing an admissible choice in \eqref{eq:pruriginoso}, are the following:
\begin{itemize}
\item[(K1)] \textit{Regional fractional $q$-Laplacian:} $\Omega\subset\real^d$ is a bounded open set and
$$
K_\Omega(x,y)=|x-y|^{-d-sq};
$$
\item[(K2)] \textit{Periodic fractional $q$-Laplacian:} $\Omega=(0,1)^d$ and
$$
K_\Omega(x,y)=\sum_{\nu\in\mathbb Z^d}|x-y-\nu|^{-d-sq};
$$
\item[(K3)] \textit{Neumann fractional Laplacian:} $\Omega\subset\real^d$ is a bounded open set with $C^2$ boundary and
$$
K_\Omega(x,y)=\frac{1}{|x-y|^{d+2s}}+\int_{\real^d\setminus\Omega}\frac{\dd z}{|x-z|^{d+2s}|y-z|^{d+2s}\int_\Omega \frac{\dd\eta}{|z-\eta|^{d+2s}}}.
$$
\item[(K4)] \textit{Neumann spectral fractional Laplacian:} $\Omega\subset\real^d$ is a bounded open set with $C^2$ boundary and
$$
K_\Omega(x,y)=\int_0^{+\infty}\frac{p_N^\Omega(t,x,y)}{t^{1+s}} \dd t,
$$
where $p_N^\Omega$ is the Neumann heat kernel in $\Omega$.
\end{itemize}
Here above $s\in(0,1)$ and $q\geq2$. As for \eqref{eq:uhu}, notice that the first three kernels are trivially greater or equal than $|x-y|^{-d-sq}$ (with $q=2$ in the third case), so that the inequality holds true for any $\varrho\geq\textup{diam}(\Omega)$. In what follows we prove that also the kernel in (K4) satisfies \eqref{eq:uhu} and that the kernels in (K1), (K3), and (K4) satisfy \eqref{eq:stic}, while the kernel in (K2) satisfies \eqref{eq:sticweak} but, in general, not \eqref{eq:stic}.

\subsubsection{The kernel in (K1)}
We simply have
$$
\begin{aligned}
    \int_\Omega\int_\Omega |x-y|^q|x-y|^{-d-sq}\dd x\dd y&\leq\int_\Omega \int_{B_{\textup{diam}(\Omega)}(x)}|x-y|^{-d+(1-s)q}\dd x\dd y\\
    &=\frac{\mathcal H^{d-1}(\partial B_1)|\Omega|\textup{diam}(\Omega)^{(1-s)q}}{(1-s)q}.
\end{aligned}
$$

\subsubsection{The kernel in (K2)}\label{sec:K2}
Let $\Omega'\subset\subset\Omega$ and $\delta=\textup{dist}(\Omega',\Omega)>0$. Then,
$$
\begin{aligned}
\int_{\Omega'}\int_\Omega  &|x-y|^q \sum_{\nu\in\mathbb Z^d}|x-y-\nu|^{-d-sq}\dd x\dd y\\
&
=\int_{\Omega'}\int_\Omega |x-y|^q|x-y|^{-d-sq}\dd x\dd y+\int_{\Omega'}\int_\Omega  |x-y|^q \sum_{\nu\in\mathbb Z^d\setminus\{0\}}|x-y-\nu|^{-d-sq}\dd x\dd y\\
&
\leq \frac{\mathcal H^{d-1}(\partial B_1)d^\frac{(1-s)q}{2}}{(1-s)q}+\textup{diam}(\Omega)^q\int_{\Omega'}\int_{\Omega^c}|x-y|^{-d-sq}\dd x\dd y\\
&
\leq\frac{\mathcal H^{d-1}(\partial B_1)d^\frac{(1-s)q}{2}}{(1-s)q}+d^\frac{q}{2}\int_{\Omega'}\int_{\real^d\setminus B_\delta(x)}|x-y|^{-d-sq}\dd x\dd y\\
&
\leq \frac{\mathcal H^{d-1}(\partial B_1)d^\frac{(1-s)q}{2}}{(1-s)q}+\frac{d^\frac{q}{2}\mathcal H^{d-1}(\partial B_1)}{sq\delta^{sq}},
\end{aligned}
$$
so that \eqref{eq:sticweak} holds true, by the symmetry of $K_\Omega$.
If however we consider $d=1$ and $sq>1$, then
$$
\begin{aligned}
\int_0^1\int_0^1& |x-y|^q\sum_{\nu\in\mathbb Z^d}|x-y-\nu|^{-d-sq}\dd x\dd y\geq \int_0^\frac{1}{4}\int_\frac{3}{4}^1 |x-y|^q |x-y+1|^{-d-sq}\dd x\dd y\\
&
\geq 2^{-q}\int_1^\frac{5}{4}\int_\frac{3}{4}^1 \frac{\dd x\dd y}{|x-y|^{1+sq}}=\frac{1}{2^q sq}\int_1^\frac{5}{4}\bigg(\frac{1}{(x-1)^{sq}}-\frac{1}{(x-3/4)^{sq}}\bigg)\dd x=+\infty,
\end{aligned}
$$
proving that \eqref{eq:stic} is not satisfied.

\subsubsection{The kernel in (K3)}\label{sec:K3}
We begin by observing that, since $\Omega$ is a bounded open set with $C^2$ boundary, there exists a constant $C\geq1$ such that
\begin{equation}\label{eq:tartaruga}
C^{-1}K_\Omega(x,y)\leq \frac{1+\log^-\Big(\frac{d_{x,y}}{|x-y|}\Big)}{|x-y|^{d+2s}}\leq C K_\Omega(x,y)\qquad\mbox{for every }x,y\in\Omega.
\end{equation}
Here above we have employed the following notations: $\log^- t=\max\{0,-\log t\}$ is the negative part of the logarithm, and $d_{x,y}=\min\{d(x),d(y)\}$, where $d(\,\cdot\,)=\textup{dist}(\,\cdot\,,\partial\Omega)$ is the distance function from the boundary of $\Omega$.
Estimate \eqref{eq:tartaruga} can be found in this form in \cite[Proposition~2.1]{audrito2022neumann}, which is a refinement of estimates (11) and (12) in the seminal paper \cite{abatangelo2020remark}.

In light of \eqref{eq:tartaruga} we have
\begin{equation}\label{eq:vivaldi}
\begin{aligned}
\int_\Omega\int_\Omega |x-y|^2 K_\Omega(x,y)\dd x\dd y&\leq C \int_\Omega\int_\Omega \frac{1+\log^-\Big(\frac{d_{x,y}}{|x-y|}\Big)}{|x-y|^{d-2+2s}}\dd x\dd y\\
&
\leq C\Bigg(1+\int_\Omega\int_\Omega \frac{\log^-\Big(\frac{d_{x,y}}{\textup{diam}(\Omega)}\Big)}{|x-y|^{d-2+2s}}\dd x\dd y\Bigg).
\end{aligned}
\end{equation}
In order to prove that $K_\Omega$ satisfies \eqref{eq:stic}, we are thus reduced to show that the last term in the above inequality is finite. For this, it is enough to assume that $\Omega$ is a bounded open set with Lipschitz boundary. Then, by \cite{doktor1976approximation} we know that there exists $\delta_0>0$, which w.l.o.g. we can assume to be smaller than $1$, such that each of the open sets
$$
\Omega_\delta\coloneqq\{x\in\Omega\,:\,\textup{dist}(x,\partial\Omega)>\delta\}\quad\textrm{with }\delta\in[0,\delta_0]
$$
has Lipschitz boundary, and
\begin{equation}\label{eq:Doktor}
\sup_{\delta\in[0,\delta_0]}\mathcal H^{d-1}(\partial\Omega_\delta)<+\infty.
\end{equation}
In the following estimates we will make use of the coarea formula for the distance function $d$, keeping in mind that its Lipschitz constant is not bigger than $1$ and $|\nabla d|=1$ at every differentiability point (see, e.\,g., \cite[Section~4]{ambrosio2000geometric}).
We have
$$
\begin{aligned}
\int_\Omega&\int_\Omega \frac{\log^-\Big(\frac{d_{x,y}}{\textup{diam}(\Omega)}\Big)}{|x-y|^{d-2+2s}}\dd x\dd y
= \int_{\Omega_{\delta_0}}\int_{\Omega_{\delta_0}} \frac{\log^-\Big(\frac{d_{x,y}}{\textup{diam}(\Omega)}\Big)}{|x-y|^{d-2+2s}}\dd x\dd y\\
&
+2\int_{\Omega_{\delta_0}}\int_{\Omega\setminus\Omega_{\delta_0}} \frac{\log^-\Big(\frac{d_{x,y}}{\textup{diam}(\Omega)}\Big)}{|x-y|^{d-2+2s}}\dd x\dd y
+\int_{\Omega\setminus\Omega_{\delta_0}}\int_{\Omega\setminus\Omega_{\delta_0}} \frac{\log^-\Big(\frac{d_{x,y}}{\textup{diam}(\Omega)}\Big)}{|x-y|^{d-2+2s}}\dd x\dd y\eqqcolon I_1+I_2+I_3,
\end{aligned}
$$
and
$$
I_1\leq \log^-\bigg(\frac{\delta_0}{\textup{diam}(\Omega)}\bigg)
\int_{\Omega_{\delta_0}}\int_{\Omega_{\delta_0}} \frac{\dd x\dd y}{|x-y|^{d-2+2s}}<+\infty.
$$
As for $I_2$ we estimate
$$
\begin{aligned}
I_2&=2\int_{\Omega_{\delta_0}}\Bigg(\int_{\Omega\setminus\Omega_{\delta_0}} \frac{\log^-\Big(\frac{d(y)}{\textup{diam}(\Omega)}\Big)}{|x-y|^{d-2+2s}}\dd y\Bigg)\dd x\\
&
\leq 2 \int_{\Omega_{\delta_0}}\Bigg(\int_{\Omega\setminus\Omega_{\delta_0}} \frac{\log^-d(y)}{|x-y|^{d-2+2s}}\dd y\Bigg)\dd x
+2|\log(\textup{diam}(\Omega))|\int_{\Omega_{\delta_0}}\int_{\Omega\setminus\Omega_{\delta_0}} \frac{\dd x\dd y}{|x-y|^{d-2+2s}}\\
&
=2\int_0^{\delta_0}\Bigg(\int_{\partial\Omega_\delta}\log^- d(y)\Bigg(\int_{\Omega_{\delta_0}}\frac{\dd x}{|x-y|^{d-2+2s}}\Bigg)\dd \mathcal H^{d-1}_y\Bigg)\dd\delta+C\\
&
\leq 2\int_0^{\delta_0}\log^- \delta\Bigg(\int_{\partial\Omega_\delta}\Bigg(\int_{B_{\textup{diam}(\Omega)}(y)}\frac{\dd x}{|x-y|^{d-2+2s}}\Bigg)\dd \mathcal H^{d-1}_y\Bigg)\dd\delta+C\\
&
=\frac{\mathcal H^{d-1}(\partial B_1)}{1-s}\textup{diam}(\Omega)^{2(1-s)}\int_0^{\delta_0}\log^- \delta\Bigg(\int_{\partial\Omega_\delta}\dd\mathcal H^{d-1}_y\Bigg)\dd\delta+C<+\infty,
\end{aligned}
$$
where in the last step we have used the integrability of the logarithm near the origin and \eqref{eq:Doktor}.
To conclude, we can similarly bound $I_3$ as follows
$$
\begin{aligned}
 I_3 &\leq \iint_{(\Omega\setminus\Omega_{\delta_0})^2} \frac{\log^-d_{x,y}}{|x-y|^{d-2+2s}}\dd x\dd y
+|\log(\textup{diam}(\Omega))|\iint_{(\Omega\setminus\Omega_{\delta_0})^2} \frac{\dd x\dd y}{|x-y|^{d-2+2s}}\\
&
=2\iint_{(\Omega\setminus\Omega_{\delta_0})^2\cap\{d(x)\leq d(y)\}} \frac{\log^-d(x)}{|x-y|^{d-2+2s}}\dd x\dd y
+C\\
&
=2\int_{\Omega\setminus\Omega_{\delta_0}}\Bigg(\int_0^{d(y)}\Bigg(\int_{\partial\Omega_\delta}\frac{\log^-d(x)}{|x-y|^{d-2+2s}}\dd \mathcal H^{d-1}_x\Bigg)\dd\delta\Bigg)\dd y+C\\
&
=2\int_{\Omega\setminus\Omega_{\delta_0}}\Bigg(\int_0^{\delta_0}\log^-\delta\Bigg(\int_{\partial\Omega_\delta}\frac{\dd \mathcal H^{d-1}_x}{|x-y|^{d-2+2s}}\Bigg)\dd\delta\Bigg)\dd y+C\\
&
\leq 2\int_0^{\delta_0}\log^- \delta\Bigg(\int_{\partial\Omega_\delta}\Bigg(\int_{B_{\textup{diam}(\Omega)}(x)}\frac{\dd y}{|x-y|^{d-2+2s}}\Bigg)\dd \mathcal H^{d-1}_x\Bigg)\dd\delta+C<+\infty.
\end{aligned}
$$

Recalling \eqref{eq:vivaldi}, by the above estimates for $I_1,I_2$ and $I_3$, we conclude that $K_\Omega$ satisfies \eqref{eq:stic} as wanted.

\subsubsection{The kernel in (K4)}
It is enough to observe that there exists $\Lambda \ge 1$ such that
$$
\frac{1}{\Lambda} \le K_{\Omega}(x,y)|x-y|^{d+2s} \le \Lambda \quad \mbox{for every } x,y \in \Omega, \ x \neq y.
$$
This estimate is obtained by the usual bounds for the heat kernel, see, e.\,g., \cite[formula~(P.4)]{abatangelo2019getting}, multiplied by $t^{-1-s}$ and integrated over $t$ between $0$ and $+\infty$.


\subsection{Fractional Cahn-Hilliard system with other boundary conditions}

The generalized fractional Cahn-Hilliard system \eqref{eq:CH1q}--\eqref{eq:iniCHq} can be coupled with other boundary conditions than the Dirichlet ones in \eqref{eq:bouCH1q}--\eqref{eq:bouCH2q}.
For example, we can consider either periodic boundary conditions or the Neumann boundary conditions (for appropriate choices of $\mathfrak L$ and $\mathfrak I$).

Interestingly, the operators $\mathfrak I=(-\Delta)^s_q$ together with periodic boundary conditions and $\mathfrak I=(-\Delta)^s$ together with Neumann boundary conditions can be rewritten as regional operators allowing weak representations of the form \eqref{eq:pruriginoso}, which satisfy the assumptions of Section~\ref{sec:sqdef}. The specific expression of the kernels can be found in Section~\ref{sec:kernels}.
In the upcoming sections we solve the corresponding problems, with the choice $\mathfrak L=(-\Delta)^\sigma$, proving the existence of weak solutions through Theorem~\ref{stm:existenceregL}. Other choices are also admissible, as observed in Section~\ref{sec:motherproblem}.

\subsubsection{Periodic boundary conditions}\label{sec:PeriodicBdSec}
In this section we consider $\Omega = (0,1)^d$ and periodic boundary conditions. More precisely, we solve the following problem
\begin{align}
\partial_t u + (-\Delta)^\sigma w = 0 &\qquad \text{in } \Omega\times(0,T), \label{eq:lampone}\\
w = (-\Delta)^s_q u + \zeta + \pi(u) &\qquad \text{in } \Omega\times(0,T), \\
u(x,0) = u_0(x) &\qquad \text{in } \real^d,\\
u(x+\nu,t) = u(x,t), \ w(x+\nu,t) = w(x,t) &\qquad \mbox{for } (x,t) \in \Omega\times(0,T) \mbox{ and } \nu \in \inte^d. \label{eq:mirtillo}
\end{align}
It does indeed make sense to consider $\inte^d$-periodic boundary conditions also for $w$, since the fractional $q$-Laplacian of a $\inte^d$-periodic function is itself $\inte^d$-periodic. Indeed, if $u \in C^2(\real^d)$ is $\inte^d$-periodic,
$$
\begin{aligned}
(-\Delta)^s_q u(x+\nu) &= \PV \int_{\real^d} \frac{|u(x+\nu)-u(y)|^{q-2}(u(x+\nu)-u(y))}{|x+\nu-y|^{d+sq}} \dd y \\
&= \PV \int_{\real^d} \frac{|u(x)-u(y)|^{q-2}(u(x)-u(y))}{|x-(y-\nu)|^{d+sq}} \dd y \\
&= \PV \int_{\real^d} \frac{|u(x)-u(z+\nu)|^{q-2}(u(x)-u(z+\nu))}{|x-z|^{d+sq}} \dd z \\
&= \PV \int_{\real^d} \frac{|u(x)-u(z)|^{q-2}(u(x)-u(z))}{|x-z|^{d+sq}} \dd z = (-\Delta)^s_q u(x),
\end{aligned}
$$
for every $x\in\real^d$ and $\nu\in\inte^d$.

Notice also that, if $u \in C^2(\real^d)$ is $\inte^d$-periodic, then
\begin{equation}\label{eq:yoyoma}
\begin{aligned}
(-\Delta)^s_q u(x) &= \PV \int_{\real^d} \frac{|u(x)-u(y)|^{q-2}(u(x)-u(y))}{|x-y|^{d+sq}} \dd y \\
&= \PV \int_{(0,1)^d} \frac{|u(x)-u(y)|^{q-2}(u(x)-u(y))}{|x-y|^{d+sq}} \dd y \\
&\qquad\qquad\qquad + \sum_{\nu\in\inte^d\setminus\{0\}} \int_{(0,1)^d+\nu} \frac{|u(x)-u(y)|^{q-2}(u(x)-u(y))}{|x-y|^{d+sq}} \dd y \\
&= \PV \int_{(0,1)^d} \frac{|u(x)-u(y)|^{q-2}(u(x)-u(y))}{|x-y|^{d+sq}} \dd y \\
&\qquad\qquad\qquad + \sum_{\nu\in\inte^d\setminus\{0\}} \int_{(0,1)^d} \frac{|u(x)-u(\xi+\nu)|^{q-2}(u(x)-u(\xi+\nu))}{|x-\xi-\nu|^{d+sq}} \dd \xi \\
&= \PV \int_{(0,1)^d} |u(x)-u(y)|^{q-2}(u(x)-u(y)) \sum_{\nu\in\inte^d} \frac{1}{|x-y-\nu|^{d+sq}} \dd y.
\end{aligned}
\end{equation}
Recall that
$$
K_\Omega(x,y) \coloneqq \sum_{\nu\in\inte^d} \frac{1}{|x-y-\nu|^{d+sq}}
$$
is the kernel in (K2) defined in Section~\ref{sec:kernels}.
Therefore, for a $\inte^d$-periodic test function $\psi \in C^2(\real^d)$ we formally have
\begin{equation}\label{eq:cellosuit}
\io (-\Delta)^s_q u(x) \psi(x) \dd x = \frac12 \iint_{\Omega^2} |u(x)-u(y)|^{q-2}(u(x)-u(y))(\psi(x)-\psi(y)) K_\Omega(x,y) \dd x\dd y.
\end{equation}
Hence, if we set $\phi(r) \coloneqq \frac{1}{2}|r|^{q-2}r$, we obtain
\begin{equation}\label{eq:john}
\io (-\Delta)^s_q u(x) \psi(x) \dd x = \langle \mathfrak I_\Omega,\psi \rangle_{K_\Omega,q}
\end{equation}
in the sense of \eqref{eq:pruriginoso}. Since $K_\Omega$ satisfies assumptions \eqref{eq:uhu} and \eqref{eq:sticweak} (see Section~\ref{sec:K2}), we are in the functional framework of Section~\ref{sec:chem}.
This setting, which is restricted to $\Omega$, is compatible with the original problem in $\real^d$ with $\inte^d$-periodic conditions, as stated in the following lemma.

\begin{lemma}\label{lemma:gigi}
    If $u : \real^d \to \real$ is $\inte^d$-periodic and belongs to $W^{s,q}_{\textup{loc}}(\real^d)$, then $u|_\Omega \in \cW^{K_\Omega,q}(\Omega)$.

    On the other hand, if $u : \Omega \to \real$ belongs to $\cW^{K_\Omega,q}(\Omega)$, then its $\inte^d$-periodic extension $\tilde{u}$ belongs to $W^{s,q}_{\textup{loc}}(\real^d)$.
\end{lemma}

\begin{proof}
     The main ingredient of the proof is the fact that we can write $\real^d = \bigcup_{\nu\in\inte^d} \big((0,1)^d+\nu\big)$, up to a set of measure $0$. We can thus consider the relevant cubes for the argument and change variables by exploiting the periodicity. We begin with the first implication. For this, we write
    $$
    \begin{aligned}
        \iint_{\Omega^2}& |u(x)-u(y)|^q \sum_{\nu\in\inte^d} \frac{1}{|x-y-\nu|^{d+sq}} \dd x\dd y\\
        &
        =\iint_{\Omega^2} |u(x)-u(y)|^q \sum_{\nu\in\{-1,0,1\}^d} \frac{1}{|x-y-\nu|^{d+sq}} \dd x\dd y\\
        &\qquad\qquad\qquad\qquad
        +\iint_{\Omega^2} |u(x)-u(y)|^q \sum_{\nu\in\inte^d\setminus\{-1,0,1\}^d} \frac{1}{|x-y-\nu|^{d+sq}} \dd x\dd y\\
        &
        \leq \iint_{(-1,2)^{2d}}  \frac{|u(x)-u(y)|^q}{|x-y|^{d+sq}} \dd x\dd y+2^q\int_\Omega |u(x)|^q\Bigg(\int_{\real^d\setminus(-1,2)^d}\frac{\dd y}{|x-y|^{d+sq}}\Bigg)\dd x<+\infty.
    \end{aligned}
    $$
    For the other implication it is enough to observe that, by $\inte^d$-periodicity, $u\in W_\textup{loc}^{s,q}(\real^d)$ if and only if $u\in W^{s,q}((-1,2)^d)$, and
    $$
    \iint_{(-1,2)^{2d}}  \frac{|u(x)-u(y)|^q}{|x-y|^{d+sq}} \dd x\dd y\leq C(d) \iint_{\Omega^2} |u(x)-u(y)|^q \sum_{\nu\in\{-2,-1,0,1,2\}^d} \frac{1}{|x-y-\nu|^{d+sq}} \dd x\dd y.
    $$
\end{proof}

We also point out that if $u : \real^d \to \real$ is $\inte^d$-periodic, then
\begin{equation}\label{eq:eros}
u \in L^q((0,1)^d) \ \implies \int_{\real^d} \frac{|u(y)|^q}{1+|y|^{d+sq}} \dd y < +\infty \ \implies \int_{\real^d} \frac{|u(y)|^{q-1}}{1+|y|^{d+sq}} \dd y < +\infty.
\end{equation}
The first implication can be proved by exploiting the periodicity of u and appropriately changing variables in the integrals and sums, whereas the second one follows by Holder’s inequality and the integrability of the weight $(1+|x|^{d+sq})^{-1}$.

Recalling \eqref{eq:yoyoma} and \eqref{eq:cellosuit}, a weak formulation of $(-\Delta)^\sigma w$ for a $\inte^d$-periodic $w$ is
\begin{equation}\label{eq:thanatos}
\langle \mathfrak{L}w,\psi \rangle_{K_\Omega,2} = B_{\mathfrak{L}}(w,\psi) =  \frac{1}{2}\iint_{\Omega^2} (w(x)-w(y))(\psi(x)-\psi(y)) K_\Omega(x,y) \dd x\dd y,
\end{equation}
where we are tacitly assuming $q=2$ and $s=\sigma$ in the $K_\Omega$ considered above, with a slight abuse of notation.
As observed in Remark~\ref{rmk:cuneo}, we are in the setting of Hypothesis~\ref{stm:lumpa}.

Thus, given a $\inte^d$-periodic datum $u_0\in W^{s,q}_\textup{loc}(\real^d)$ such that
$$
\Gamma(u_0)\in L^1(\Omega)\quad\mbox{and}\quad\mathfrak m(u_0)\in\textup{Int}\,D(\gamma),
$$
we can apply Theorem~\ref{stm:existenceregL} to deduce the existence of a solution to the following system
\begin{align*}
\partial_t u + \mathfrak L w = 0 &\qquad \text{in } (\cW^{K_\Omega,2}(\Omega))',\\
w = \mathfrak I_\Omega u + \zeta + \pi(u) &\qquad \text{in } (\cW^{K_\Omega,q}(\Omega))',
\end{align*}
almost everywhere in $(0,+\infty)$, with $u(0) = u_0|_\Omega$ almost everywhere in $\Omega$. Moreover, the mass of $u$ is conserved and we have the uniform energy estimate
$$
\frac{1}{4}\int_0^t \|w(\tau)-\mathfrak m(w(\tau))\|_{K_\Omega,2}^2\dd \tau+\frac{1}{2q}\|u(t)\|_{K_\Omega,q}^q+\int_\Omega F(u(x,t))\dd x
\leq \frac{1}{2q}\|u_0\|_{K_\Omega,q}^q+\int_\Omega F(u_0(x))\dd x
$$
for every $t\geq 0$.

The above is a weak formulation of \eqref{eq:lampone}--\eqref{eq:mirtillo}. Indeed, if $\tilde{u}$ is the $\inte^d$-periodic extension of $u$, then $\tilde{u}\in W^{s,q}_\textup{loc}(\real^d)$ by Lemma~\ref{lemma:gigi}. Hence, since \eqref{eq:eros} holds true, if we also assume $u$ to be locally $C^2$ in $\Omega$, then $(-\Delta)^s_q \tilde{u}$ is well-defined everywhere in $\Omega$ and we can exploit identity \eqref{eq:john}. We can clearly argue in the same way for $\tilde{w}$.

\subsubsection{Fractional Neumann boundary conditions}
In this section the open set $\Omega$ is not necessarily connected (unless otherwise stated).

As we already mentioned in the introduction, in \cite{dipierro2017nonlocal} the authors proposed the following notion of nonlocal normal derivative
$$
\mathcal N_s u(x)\coloneqq (1-s)\int_\Omega\frac{u(x)-u(y)}{|x-y|^{d+2s}}\dd y\qquad\mbox{for }x\in\real^d\setminus\Omega,
$$
which allows a natural formulation of Neumann boundary conditions for the fractional Laplacian. Indeed, as shown in \cite{dipierro2017nonlocal} the following fractional counterpart of the classical integration by parts formulae hold true:
$$
\int_\Omega(-\Delta)^s u\dd x=-\int_{\real^d\setminus\Omega} \mathcal N_s u\dd x
$$
and
$$
\frac{1-s}{2}\iint_{Q(\Omega)}\frac{(u(x)-u(y))(v(x)-v(y))}{|x-y|^{d+2s}}\dd x\dd y=\int_\Omega v(-\Delta)^s u\dd x+\int_{\real^d\setminus\Omega}v\, \mathcal N_s u\dd x.
$$
Moreover, the corresponding Neumann problem
\begin{equation}\label{eq:nocciolina}
\begin{aligned}
    (-\Delta)^s u=f&\qquad\mbox{in }\Omega\\
    \mathcal N_s u=g&\qquad\mbox{in }\real^d\setminus\Omega
\end{aligned}
\end{equation}
has a variational nature, as indeed free critical points of the functional
$$
\frac{1-s}{4}\iint_{Q(\Omega)}\frac{|u(x)-u(y)|^2}{|x-y|^{d+2s}}\dd x\dd y-\int_\Omega fu\dd x-\int_{\real^d\setminus\Omega}gu\dd x
$$
correspond to weak solutions of \eqref{eq:nocciolina}.
As a consequence we can address Metatheorem~\ref{stm:Metathm3}, thus considering the following fractional Cahn-Hilliard system
\begin{align}
\partial_t u + (-\Delta)^\sigma w = 0 &\qquad \text{in } \Omega\times(0,T),\label{eq:liquido} \\
w = (-\Delta)^s u + \zeta + \pi(u) &\qquad \text{in } \Omega\times(0,T),\label{eq:solido} \\
\mathcal N_\sigma w = 0 &\qquad \text{in } (\real^{d}\setminus\Omega)\times(0,T),\label{eq:gassoso}\\
\mathcal N_s u = 0 &\qquad \text{in } (\real^{d}\setminus\Omega)\times(0,T),\\
u(x,0) = u_0(x) &\qquad \text{in } \Omega.\label{eq:narcotic}
\end{align}
The above observations yield a natural variational framework for the weak formulation of this problem.

As a first step we define the seminorm
$$
\|u\|^2_{\cH^s_\Omega}\coloneqq \frac{1-s}{2}\iint_{Q(\Omega)}\frac{|u(x)-u(y)|^2}{|x-y|^{d+2s}}\dd x\dd y,
$$
and the functional space
$$
\cH^s_\Omega\coloneqq\left\{u:\real^d\to\real\,:\,u\in L^2(\Omega)\mbox{ and }\|u\|_{\cH^s_\Omega}<+\infty\right\},
$$
endowed with the norm $\|\,\cdot\,\|_{L^2(\Omega)}+\|\,\cdot\,\|_{\cH^s_\Omega}$. As proved in \cite[Proposition~3.1]{dipierro2017nonlocal}, this is a Hilbert space.
Notice that
$$
u\in\cH^s_\Omega\quad\implies\quad\int_{\real^d}\frac{|u(x)|^2}{1+|x|^{d+2s}}\dd x < +\infty,
$$
hence $u\in L^2_\textup{loc}(\real^d)$. For a proof, see, e.\,g., the beginning of the proof of \cite[Proposition~7.1\,(i)]{caffarelli2010nonlocal} or \cite[Lemma~D.1.3]{lombardini2018minimization}.

To stress the difference between the functional spaces employed here and those for the operator \eqref{eq:w-slapl}, as a weak version of the fractional Laplacian we also employ the notation
$$
\langle\mathfrak D^s u,v\rangle_{(\cH^s_\Omega)'}\coloneqq \frac{1-s}{2}\iint_{Q(\Omega)}\frac{(u(x)-u(y))(v(x)-v(y))}{|x-y|^{d+2s}}\dd x\dd y
$$
for every $u,v\in\cH^s_\Omega$, which is the first variation of $\frac{1}{2}\|\,\cdot\,\|_{\cH^s_\Omega}^2$ as indeed
$$
\frac{\dd}{\dd\varepsilon}\Big|_{\varepsilon=0}\frac{1}{2}\|u+\varepsilon v\|_{\cH^s_\Omega}^2=\langle\mathfrak D^s u,v\rangle_{(\cH^s_\Omega)'}.
$$
This can be proved via the same proof of Lemma~\ref{lem:nonhouncomodino}. Hence, we can give the following definition of weak solution to system \eqref{eq:liquido}--\eqref{eq:narcotic}.

\begin{definition}\label{def:weaksolNeu}
Let $T>0$ be fixed. We say that $(u,w,\zeta)$ is a weak solution to the Cahn-Hilliard system \eqref{eq:liquido}--\eqref{eq:narcotic} associated with the initial datum $u_0:\real^d\to\real$ such that $u_0\in H^s(\Omega)$ if
$$
\begin{aligned}
u &\in L^\infty(0,T;\cH^s_\Omega) \cap W^{1,2}(0,T;(\cH^\sigma_\Omega)'), \\
w &\in L^2(0,T;\cH^\sigma_\Omega),\\
\zeta &\in L^2(0,T;L^2(\Omega)), \qquad \zeta\in\gamma(u) \quad \text{a.\,e. in } \Omega \times (0,T),
\end{aligned}
$$
and $(u,w,\zeta)$ satisfies the following weak formulation of \eqref{eq:liquido}--\eqref{eq:solido}:
\begin{align}
\partial_t u + \mathfrak D^\sigma w = 0 &\qquad \text{in } (\cH^\sigma_\Omega)', \label{eq:CH1sigomega}\\
w = \mathfrak D^s u + \zeta + \pi(u) &\qquad \text{in } (\cH^s_\Omega)',\label{eq:CH2somega}
\end{align}
almost everywhere in $(0,T)$, with $u(0) = u_0$ almost everywhere in $\Omega$.
\end{definition}
We recall that for every $x\in\real^d\setminus\overline{\Omega}$ we have the following equivalence
\begin{equation}\label{eq:sherpa}
\mathcal N_s v(x)=0\quad\Longleftrightarrow\quad v(x)=\frac{\displaystyle\int_\Omega\frac{v(y)}{|x-y|^{d+2s}}\dd y}{\displaystyle\int_\Omega\frac{\dd y}{|x-y|^{d+2s}}}.
\end{equation}

\begin{remark}
    It is possible to define a natural notion of nonlocal normal derivative also for the $(s,q)$-Laplacian, with $q\not=2$, by still relying on an integration by parts formula, as observed in \cite[Theorem~6.3]{barrios2020neumann}. As in the case $q=2$, this underpins a variational framework for studying problems of the kind \eqref{eq:nocciolina}, see, e.\,g., \cite[Section~6.3]{barrios2020neumann} and \cite{mugnai2019neumann}. As a consequence, we can give an analogous of Definition~\ref{def:weaksolNeu}, with the $(s,q)$-Laplacian in \eqref{eq:CH2somega}. On the other hand, the equivalence \eqref{eq:sherpa} relies on the linearity of the numerator of the integrand defining $\mathcal N_s$, which holds only for the case $q=2$. Since \eqref{eq:sherpa} is crucial in our solving system \eqref{eq:liquido}--\eqref{eq:narcotic}, as it allows us to reformulate it as a regional problem, we limit ourselves to the case $q=2$.
\end{remark}

\begin{remark}\label{rmk:fammilaspesa}
    We gather some observations clarifying Definition~\ref{def:weaksolNeu}.
\begin{enumerate}
    \item Concerning the interpretation of $u \in L^\infty(0,T;\cH^s_\Omega) \cap W^{1,2}(0,T;(\cH^\sigma_\Omega)')$ and of equation \eqref{eq:CH1sigomega}, we regard a function $v\in\cH^s_\Omega$ as an element of $(\cH^\sigma_\Omega)'$ via
    $$
\varphi\in\cH^\sigma_\Omega\mapsto\int_\Omega v\varphi\dd x.
    $$
    This is rigorously achieved as follows. We consider the restriction map $\cH^s_\Omega\to H^s(\Omega)$, which induces the embedding $(H^s(\Omega))'\hookrightarrow (\cH^s_\Omega)'$ since it is surjective by \cite[Theorem~5.4]{di2012hitchhiker}, and the triple $H^s(\Omega)\hookrightarrow L^2(\Omega)\simeq(L^2(\Omega))'\hookrightarrow (H^s(\Omega))'$. Similarly for $\sigma$. Combining these maps, we have in particular $\cH^s_\Omega\to (L^2(\Omega))'\hookrightarrow (\cH^\sigma_\Omega)'$.

    We do the same for equation \eqref{eq:CH2somega}.
    The precise interpretation of the system is thus
    \begin{align*}
\partial_t (u|_\Omega) + \mathfrak D^\sigma w = 0 &\qquad \text{in } (\cH^\sigma_\Omega)',\\
w|_\Omega = \mathfrak D^s u + \zeta + \pi(u|_\Omega) &\qquad \text{in } (\cH^s_\Omega)'.
\end{align*}
\item The fractional Neumann boundary conditions
\begin{align*}
    \mathcal N_\sigma w=0 &\qquad \text{a.\,e. in }\real^d\setminus\Omega,\\
    \mathcal N_s u=0 &\qquad \text{a.\,e. in }\real^d\setminus\Omega,
\end{align*}
hold true whenever both equations \eqref{eq:CH1sigomega}--\eqref{eq:CH2somega} are satisfied. We begin by observing that, if $v\in C_c^1(\real^d\setminus\overline{\Omega})$, then, since $\textup{dist}(\Omega,\textup{supp}\,v)>0$ and $\textup{supp}\,v$ is compact, by symmetry we have
\begin{align*}
    \iint_{Q(\Omega)}\frac{(w(x,t)-w(y,t))(v(x)-v(y))}{|x-y|^{d+2\sigma}}\dd x\dd y&
=-2\int_\Omega\Bigg(\int_{\textup{supp}\,v}\frac{w(x,t)-w(y,t)}{|x-y|^{d+2\sigma}}v(y)\dd y\Bigg)\dd x\\
&
=\frac{2}{1-\sigma}\int_{\real^d}v(y)\mathcal N_\sigma w(y,t)\dd y.
\end{align*}
Thus, if \eqref{eq:CH1sigomega} is satisfied in $t\in(0,T)$, by point (1) we obtain
$$
0=\langle\partial_t u(t),v\rangle_{(\cH^\sigma_\Omega)'}=-\langle\mathfrak D^\sigma w(t),v\rangle_{(\cH^\sigma_\Omega)'}=-\int_{\real^d}v(y)\mathcal N_\sigma w(y,t)\dd y,
$$
for every $v\in C_c^1(\real^d\setminus\overline{\Omega})$.
This implies that
$$
\mathcal N_\sigma w(y,t)=0\qquad\textrm{for a.\,e. }y\in\real^d\setminus\Omega,
$$
as claimed. Similarly for \eqref{eq:CH2somega} and $\mathcal N_s u$.

\item Concerning the initial datum $u_0$, since we are working in $\real^d$ it has to be defined in the whole space. Nevertheless, the meaningful part is only $u_0|_\Omega$. Indeed, by point (2) and \eqref{eq:sherpa} we know that
\begin{equation}\label{eq:spleen}
u(x,t)=\frac{\displaystyle\int_\Omega\frac{u(y,t)}{|x-y|^{d+2s}}\dd y}{\displaystyle\int_\Omega\frac{\dd y}{|x-y|^{d+2s}}}\ \mbox{ for a.\,e. }x\in\real^d\setminus\Omega\mbox{ and a.\,e. }t\in(0,T).
\end{equation}
Moreover, by point (1) we have the embeddings $H^s(\Omega)\hookrightarrow L^2(\Omega)\hookrightarrow (\cH^\sigma_\Omega)'$, so that, by the Aubin-Lions-Simon compactness lemma
$$
u\in L^\infty(0,T;H^s(\Omega))\cap W^{1,2}(0,T;(\cH^\sigma_\Omega)')\subset C([0,T];L^2(\Omega)).
$$
As a consequence, the initial condition $u(0)=u_0$ almost everywhere in $\Omega$ makes sense, and $u(t)\to u_0$ strongly in $L^2(\Omega)$ as $t\to0$.
By \eqref{eq:spleen}, this also implies that
$$
u(x,t)\xrightarrow{t\to0}\frac{\displaystyle\int_\Omega\frac{u_0(y)}{|x-y|^{d+2s}}\dd y}{\displaystyle\int_\Omega\frac{\dd y}{|x-y|^{d+2s}}}\quad\mbox{strongly in }L^2_\textup{loc}(\real^d\setminus\overline{\Omega}).
$$
We thus have two options for the choice of the initial condition: either we impose $\mathcal N_s u_0=0$ almost everywhere in $\real^d\setminus\Omega$, and then we have $u(0)=u_0$ almost everywhere in $\real^d$, or the initial condition can be satisfied only in $\Omega$.
\end{enumerate}
\end{remark}

As a second step, thanks to the equivalence \eqref{eq:sherpa}, in \cite[Corollary~1.2]{abatangelo2020remark} the author observed that if $u:\real^d\to\real$ is such that $u\in C^2_{\textup{loc}}(\Omega)\cap L^\infty(\Omega)$ and $\mathcal N_s u=0$ almost everywhere in $\real^d\setminus\Omega$, then the fractional Laplacian of $u$ in $\Omega$ can be rewritten as a regional operator,
$$
(-\Delta)^s u(x)=\PV\int_\Omega (u(x)-u(y))K_\Omega^s(x,y)\dd y\quad\mbox{for every }x\in\Omega,
$$
where $K_\Omega^s$ denotes the kernel in (K3) defined in Section~\ref{sec:kernels}, multiplied by $(1-s)$. Since $K_\Omega^s$ is symmetric in $x$ and $y$, we formally have the following natural weak formulations
\begin{equation}\label{eq:saudade}
\langle \mathfrak D^s u,v\rangle_{(\cH^s_\Omega)'} = \int_\Omega v(-\Delta)^s u\dd x=\frac{1}{2}\iint_{\Omega^2} (u(x)-u(y))(v(x)-v(y))K_\Omega^s(x,y)\dd x\dd y=\langle\mathfrak I_\Omega u,v\rangle_{K_\Omega^s,2}.
\end{equation}
If $\Omega\subset\real^d$ is a bounded open set with $C^2$ boundary, the operator, as written in the form of the last equality, with $\phi(r)=\frac{1}{2}r$, fits into the functional setting of Section~\ref{sec:chem}, thanks to Section~\ref{sec:K3}.
We thus have a second weak version of system \eqref{eq:liquido}--\eqref{eq:narcotic}, given by Definition~\ref{def:weaksolregL}, with $X=\cW^{K_\Omega^\sigma,2}(\Omega)$,
\begin{equation}\label{eq:biella}
\langle\mathfrak L u,v\rangle_{X'}=B_\mathfrak L(u,v)=\frac{1}{2}\iint_{\Omega^2} (u(x)-u(y))(v(x)-v(y))K_\Omega^\sigma(x,y)\dd x\dd y,
\end{equation}
and $\mathfrak I_\Omega$ as in \eqref{eq:saudade} (see Remark~\ref{rmk:cuneo} for the details).

These two weak formulations are actually equivalent, as stated in the following.

\begin{theorem}\label{stm:Teo}
Let $\Omega\subset\real^d$ be a bounded open set with $C^2$ boundary.
    If $(u,w,\zeta)$ is a weak solution in the sense of Definition~\ref{def:weaksolNeu}, then $(u|_\Omega,w|_\Omega,\zeta)$ is a weak solution in the sense of Definition~\ref{def:weaksolregL}. Conversely, if $(u,w,\zeta)$ is a weak solution in the sense of Definition~\ref{def:weaksolregL}, and we extend $u$ and $w$ by setting
    \begin{equation}\label{eq:Neumannext}
\tilde{u}(x,t)=\frac{\displaystyle\int_\Omega\frac{u(y,t)}{|x-y|^{d+2s}}\dd y}{\displaystyle\int_\Omega\frac{\dd y}{|x-y|^{d+2s}}}\quad\mbox{and}\quad
\tilde{w}(x,t)=\frac{\displaystyle\int_\Omega\frac{w(y,t)}{|x-y|^{d+2\sigma}}\dd y}{\displaystyle\int_\Omega\frac{\dd y}{|x-y|^{d+2\sigma}}}
    \end{equation}
    for every $(x,t)\in(\real^d\setminus\overline{\Omega})\times(0,T)$, then $(\tilde{u},\tilde{w},\zeta)$ is a weak solution in the sense of Definition~\ref{def:weaksolNeu}.
\end{theorem}

In order to prove Theorem~\ref{stm:Teo} we need two ingredients: the following compatibility result between the two weak functional settings, which is \cite[Lemma~A.2]{audrito2022neumann}, and an approximation result, whose proof we give in Appendix~\ref{sec:approximation}.

\begin{lemma}[{\cite[Lemma~A.2]{audrito2022neumann}}]\label{lem:Xavier}
Let $v:\real^d\to\real$ be such that $v\in L^2(\Omega)$ and $\mathcal N_s v=0$ almost everywhere in $\real^d\setminus\Omega$. Then $v\in\cH^s_\Omega$ if and only if $v|_\Omega\in\cW^{K^s_\Omega,2}(\Omega)$, and
$$
\|v\|_{\cH^s_\Omega}^2=\frac{1}{2}\|v\|_{K_\Omega^s,2}^2.
$$
Moreover, let $v,\varphi\in\cH^s_\Omega$ be such that $\varphi|_\Omega\in\cW^{K^s_\Omega,2}(\Omega)$ and $\mathcal N_s v=0$ almost everywhere in $\real^d\setminus\Omega$. Then
$$
\langle\mathfrak D^s v,\varphi\rangle_{(\cH^s_\Omega)'}=\langle\mathfrak I_\Omega v,\varphi\rangle_{K^s_\Omega,2}.
$$
\end{lemma}

\begin{lemma}\label{lem:approximation}
    Let $\Omega\subset\real^d$ be a bounded open set with $C^2$ boundary. Given $\varphi\in\cH^s_\Omega$, there exists a sequence $\{\varphi_k\}\subset C^\infty(\real^d)\cap L^\infty(\real^d)$ such that
    \begin{equation}\label{eq:ennui}
    \begin{aligned}
        & \varphi_k\to \varphi\quad\textrm{strongly in }L^2_\textup{loc}(\real^d)\mbox{ and a.\,e. in }\real^d,\\
        &
        \sup_k \|\varphi_k\|_{\cH^s_\Omega}<+\infty,\\
        &
        \langle\mathfrak D^s v,\varphi_k\rangle_{(\cH^s_\Omega)}\to
        \langle\mathfrak D^s v,\varphi\rangle_{(\cH^s_\Omega)}\quad\mbox{for every }v\in\cH^s_\Omega.
    \end{aligned}
    \end{equation}
\end{lemma}

The statement of the lemma is probably far from being optimal. Indeed, it is reasonable to expect that the $C^2$-regularity requirement may be weakened to Lipschitz, and that one can actually obtain the strong convergence in the $\cH^s_\Omega$-norm. However, Lemma~\ref{lem:approximation} in its present form is enough for our purposes and its proof is rather straightforward.

\begin{proof}[Proof of Theorem~\ref{stm:Teo}]
\textbf{Definition}~\ref{def:weaksolNeu} $\implies$  \textbf{Definition}~\ref{def:weaksolregL}. Let $(u,w,\zeta)$ be a weak solution in the sense of Definition~\ref{def:weaksolNeu}.
By Remark~\ref{rmk:fammilaspesa}\,(2) we know that
\begin{align*}
    \mathcal N_\sigma w(x,t)=0 &\qquad \text{for a.\,e. }x\in \real^d\setminus\Omega\mbox{ and }t\in(0,T),\\
    \mathcal N_s u(x,t)=0 &\qquad \text{for a.\,e. }x\in \real^d\setminus\Omega\mbox{ and }t\in(0,T).
\end{align*}
Hence, by Lemma~\ref{lem:Xavier} we have $u|_\Omega\in L^\infty(0,T;\cW^{K^s_\Omega,2}(\Omega))$ and $w|_\Omega\in L^2(0,T;\cW^{K^\sigma_\Omega,2}(\Omega))$. Moreover, again by Lemma~\ref{lem:Xavier}, we have the isomorphism
$$
\cW^{K^\sigma_\Omega,2}(\Omega)\xrightarrow{\sim}\big\{v\in\cH^\sigma_\Omega\,:\,\mathcal N_\sigma v=0\mbox{ a.\,e. in }\real^d\setminus\Omega\big\},
$$
given by
$$
v\longmapsto\tilde{v}\qquad\mbox{ s.\,t. }\tilde{v}|_\Omega=v\quad\mbox{and}\quad\tilde{v}(x)=\frac{\displaystyle\int_\Omega\frac{v(y)}{|x-y|^{d+2\sigma}}\dd y}{\displaystyle\int_\Omega\frac{\dd y}{|x-y|^{d+2\sigma}}}\mbox{ for }x\in\real^d\setminus\overline{\Omega}.
$$
Thus we have also $u|_\Omega\in W^{1,2}(0,T;(\cW^{K^\sigma_\Omega,2}(\Omega))')$ and
$$
\partial_t(u|_\Omega)=-\langle\mathfrak D^\sigma w,\tilde{v}\rangle_{(\cH^\sigma_\Omega)}=-\langle\mathfrak L w|_\Omega,v\rangle_{K^\sigma_\Omega,2}\quad\mbox{for every }v\in \cW^{K^\sigma_\Omega,2}(\Omega),
$$
by Lemma~\ref{lem:Xavier}, for almost every $t\in(0,T)$. Similarly for equation \eqref{eq:CH2somega}.

\textbf{Definition}~\ref{def:weaksolregL} $\implies$  \textbf{Definition}~\ref{def:weaksolNeu}. Let $(u,w,\zeta)$ be a weak solution in the sense of Definition~\ref{def:weaksolregL} and let $\tilde{u}$ and $\tilde{w}$ be extended as in \eqref{eq:Neumannext}. By Lemma~\ref{lem:Xavier} we have $\tilde{u}\in L^\infty(0,T;\cH^s_\Omega)$ and $\tilde{w}\in L^2(0,T;\cH^\sigma_\Omega)$. Now let $v\in\cH^\sigma_\Omega$ and let $\{v_k\}$ be an approximating sequence in the sense of Lemma~\ref{lem:approximation}. Since $v_k|_\Omega\in C^{0,1}(\overline{\Omega})$, by the first inclusion in \eqref{eq:lapras} and Section~\ref{sec:K3}, we know that $v_k|_\Omega\in\cW^{K_\Omega^\sigma,2}(\Omega)$. Then, given $\varphi\in C^\infty_c(0,T)$, by the definition of weak derivative in Bochner spaces, by equation \eqref{eq:CHw1regL} and by Lemma~\ref{lem:Xavier}, we obtain
$$
\begin{aligned}
    \int_0^T\varphi'(t)\int_\Omega u(x,t)v_k|_\Omega(x)\dd x\dd t&=-\int_0^T\varphi(t)\langle\partial_t u(t),v_k|_\Omega\rangle_{K^\sigma_\Omega,2}\dd t =\int_0^T\varphi(t)\langle\mathfrak L w(t),v_k|_\Omega\rangle_{K^\sigma_\Omega,2}\dd t\\
    &
    =\int_0^T\varphi(t)\langle\mathfrak D^\sigma \tilde{w}(t),v_k\rangle_{(\cH^\sigma_\Omega)'}\dd t.
\end{aligned}
$$
By Lemma~\ref{lem:approximation}, we have
$$
\lim_{k\to\infty}\int_0^T\varphi'(t)\int_\Omega u(x,t)v_k|_\Omega(x)\dd x\dd t=\int_0^T\varphi'(t)\int_\Omega u(x,t)v|_\Omega(x)\dd x\dd t,
$$
and
$$
\lim_{k\to\infty}\varphi(t)\langle\mathfrak D^\sigma\tilde{w}(t),v_k\rangle_{(\cH^\sigma_\Omega)'}=\varphi(t)\langle\mathfrak D^\sigma\tilde{w}(t),v\rangle_{(\cH^\sigma_\Omega)'}
$$
for every $t\in(0,T)$. Moreover
$$
|\varphi(t)\langle\mathfrak D^\sigma\tilde{w}(t),v_k\rangle_{(\cH^\sigma_\Omega)'}|\leq \|\varphi\|_{L^\infty(0,T)}\|\tilde{w}(t)\|_{\cH^\sigma_\Omega}\sup_k\|v_k\|_{\cH^\sigma_\Omega}.
$$
Lebesgue's dominated convergence theorem then yields
$$
\lim_{k\to\infty}\int_0^T\varphi(t)\langle\mathfrak D^\sigma \tilde{w}(t),v_k\rangle_{(\cH^\sigma_\Omega)'}\dd t=\int_0^T\varphi(t)\langle\mathfrak D^\sigma \tilde{w}(t),v\rangle_{(\cH^\sigma_\Omega)'}\dd t,
$$
hence
$$
\int_0^T\varphi'(t)\int_\Omega \tilde{u}|_\Omega(x,t)v|_\Omega(x)\dd x\dd t=\int_0^T\varphi(t)\langle\mathfrak D^\sigma \tilde{w}(t),v\rangle_{(\cH^\sigma_\Omega)'}\dd t.
$$
Since this holds for every $v\in\cH^\sigma_\Omega$ and $\varphi\in C^\infty_c(0,T)$, by the definition of weak derivative in Bochner spaces we have that $\tilde{u}|_\Omega\in W^{1,2}(0,T;(\cH^\sigma_\Omega)')$ and $\partial_t(\tilde{u}|_\Omega)=-\mathfrak D^\sigma\tilde{w}$. The fundamental lemma of the calculus of variations for Bochner spaces thus implies \eqref{eq:CH1sigomega} almost everywhere in $(0,T)$.
A similar argument ensures the validity of \eqref{eq:CH2somega}, concluding the proof.
\end{proof}

As a consequence, we are able to prove the existence of weak solutions as detailed here below.

\begin{coro}[Existence and uniqueness]\label{cor:NeuExi}
Let $\Omega\subset\real^d$ be an open set with $C^2$ boundary (not necessarily connected) and let Hypothesis~\ref{stm:hyp} be satisfied. Assume
$$
u_0 \in \cW^{K_\Omega^s,2}(\Omega), \quad \Gamma(u_0) \in L^1(\Omega)\quad\mbox{and}\quad\mathfrak m(u_0|_\Omega)\in \textup{Int}\,D(\gamma).
$$
Then, the Cahn-Hilliard system \eqref{eq:liquido}--\eqref{eq:narcotic} admits a solution $(u,w,\zeta)$, in the sense of Definition~\ref{def:weaksolNeu} for every $T>0$. Moreover,
\begin{equation*}
\mathfrak m(u|_\Omega(t))=\mathfrak m(u_0|_\Omega)\quad\mbox{for every }t\geq 0,
\end{equation*}
and, defining the energy
$$
\mathscr E(u(t)) \coloneqq \frac{1}{2}\|u(t)\|_{\cH^s_\Omega}^2 + \io F(u(x,t))\dd x,
$$
the following energy estimate holds:
\begin{equation*}
\frac{1}{4}\int_0^t \|w(\tau)-\mathfrak m(w|_\Omega(\tau))\|_{\cH^\sigma_\Omega}^2 \dd \tau +  \mathscr E(u(t)) \le \mathscr E(\tilde{u_0})\qquad\mbox{for every }t\geq 0,
\end{equation*}
where the function $\tilde{u_0}$ is defined as follows
$$
\tilde{u_0}|_\Omega=u_0\quad\mbox{and}\quad\tilde{u_0}(x)=\frac{\displaystyle\int_\Omega\frac{u_0(y)}{|x-y|^{d+2\sigma}}\dd y}{\displaystyle\int_\Omega\frac{\dd y}{|x-y|^{d+2\sigma}}}\mbox{ for }x\in\real^d\setminus\overline{\Omega}.
$$

Moreover, $u$ is unique and
\begin{equation*}
\mathscr E(u(t)) \le \mathscr E(u(\tau)) \qquad \mbox{for } \ 0 \le \tau \le t.
\end{equation*}
If we also assume $\Gamma$ to be differentiable, then also $w$ and $\zeta$ are unique.
On the other hand, if there is an element of $\textup{Int}\,D(\gamma)$ at which $\gamma$ is multivalued, then there exists an initial datum for which the solution is not unique.
\end{coro}

The proof consists only in applying Theorem~\ref{stm:existenceregL} for $\mathfrak L$ as in \eqref{eq:biella} and $\mathfrak I_\Omega$ as in \eqref{eq:saudade}, and then exploiting Theorem~\ref{stm:Teo} and Lemma~\ref{lem:Xavier} to translate everything in the language of Definition~\ref{def:weaksolNeu}. This motivates in particular the assumption of $C^2$-regularity for the boundary of $\Omega$. The connectedness assumption, usually required in order to ensure the validity of the Poincar\'e inequality, can be dropped in view of Proposition~\ref{stm:poinconnected} (which also gives the necessary coercivity of the bilinear functional $B_\mathfrak L$ defined in \eqref{eq:biella}). The assumptions on the initial datum further explain point (3) of Remark~\ref{rmk:fammilaspesa} as they only involve the restriction of $u_0$ to $\Omega$. In particular, concerning the assumption $u_0 \in \cW^{K_\Omega^s,2}(\Omega)$, in light of \eqref{eq:lapras} it is enough to require $u_0|_\Omega\in C^{0,1}(\overline{\Omega})$, regardless of its behavior outside of $\Omega$.

\subsubsection{Other examples}\label{sec:motherproblem}
Some natural operators $\mathfrak L$ satisfying Hypothesis~\ref{stm:lumpa}, which can hence be considered in equation \eqref{eq:CHw1regL}, are the following ones:
\begin{itemize}
    \item $\Omega=(0,1)^d$ and $\mathfrak L=-\Delta$, with periodic boundary conditions;
    \item $\Omega\subset\real^d$ a bounded and connected open set with Lipschitz boundary, and $\mathfrak L=-\Delta$, with Neumann boundary conditions;
    \item $\Omega=(0,1)^d$ and $\mathfrak L=(-\Delta)^\sigma$, with periodic boundary conditions, in the form \eqref{eq:thanatos};
    \item $\Omega\subset\real^d$ a bounded open set with $C^2$ boundary, and $\mathfrak L=(-\Delta)^\sigma$, with fractional Neumann boundary conditions, in the form \eqref{eq:biella}; an alternative notion of weak solution can be given in the spirit of Definition~\ref{def:weaksolNeu}, as in \eqref{eq:CH1sigomega}. One can actually show that these two formulations are equivalent, by arguing as in the proof of Theorem~\ref{stm:Teo}.
\end{itemize}
Each of the resulting equations can be coupled with
equation \eqref{eq:CHw2regL}, where we consider one of the following:
\begin{itemize}
    \item $\Omega=(0,1)^d$ and $\mathfrak I_\Omega=(-\Delta)^s_q$, with periodic boundary conditions, as in \eqref{eq:cellosuit}--\eqref{eq:john};
    \item $\Omega\subset\real^d$ a bounded open set with $C^2$ boundary, and $\mathfrak I_\Omega=(-\Delta)^s$, with fractional Neumann boundary conditions, in the form \eqref{eq:saudade}; as in the last point above, an alternative notion of weak solution can be given in the spirit of Definition~\ref{def:weaksolNeu}, as in \eqref{eq:CH2somega}. Once again, one can actually show that these two formulations are equivalent, by arguing as in the proof of Theorem~\ref{stm:Teo}.
\end{itemize}

\subsection{Spectral fractional Laplacian}\label{sec:geniodelletartarughe}
A different definition of fractional Laplacian can be given via spectral theory. For example, let $\psi_k^D$ be an eigenfunction corresponding to the $k$th eigenvalue of the Dirichlet Laplacian, namely
$$
\left\lbrace
\begin{aligned}
    -\Delta\psi_k^D &= \lambda_k^D\psi_k^D \ \mbox{in } \Omega, \\
    \psi_k^D &\in H^1_0(\Omega),
\end{aligned}
\right.
$$
with $0<\lambda_0^D<\lambda_1^D\leq\lambda_2^D\leq\dots$, normalized to make $\{\psi_k^D\}_k$ a normal basis of $L^2(\Omega)$. We can then write
$$
u(x) = \sum_{k=0}^{+\infty} u_k^D\psi_k^D(x),
$$
for any $u \in L^2(\Omega)$. The corresponding spectral $s$-fractional Laplacian with Dirichlet boundary conditions is defined as
$$
(-\Delta)^s_{D,\Omega} u(x) \coloneqq \sum_{k=0}^{+\infty} (\lambda_k^D)^s u_k^D\psi_k^D(x).
$$
In general, this operator is different from the integrodifferential fractional Laplacian $(-\Delta)^s$, see, e.\,g., \cite{servadei2014spectrum}.
Interestingly, they coincide instead, up to a constant $c=c(s,d)>0$, for periodic functions, that is
$$
(-\Delta)^s_{D,\Omega} u(x)= c (-\Delta)^s u(x),
$$
if $u:\real^d\to\real$ is $\inte^d$-periodic (and regular enough), see, e.\,g., \cite[formula~(2.53) and Appendix~Q]{abatangelo2019getting}.
As observed in \eqref{eq:yoyoma}, for such functions, we also have
$$
(-\Delta)^s u(x)=\PV\int_{(0,1)^d}(u(x)-u(y))K_\Omega(x,y)\dd y,
$$
where $K_\Omega$ is the kernel in (K2) defined in Section~\ref{sec:kernels} (with $q=2$). In the case in which $\Omega=(0,1)^d$, the spectral fractional Laplacian $(-\Delta)^s_{D,\Omega}$ can thus be represented as a regional integrodifferential operator that fits into the functional setting of Section~\ref{sec:chem}.

In a general open set, the operator $(-\Delta)^s_{D,\Omega}$ can still be represented as a regional integrodifferential operator with kernel $J$, plus a weighted local term of order zero, see \cite[formula~(3)]{abatangelo2017nonhomogeneous}.

Another possibility in the definition of the spectral fractional Laplacian consists in considering Neumann boundary conditions instead of Dirichlet ones.
Let now $\Omega\subset\real^d$ be a bounded open set with $C^2$ boundary, and let $\psi_k^N$ be an eigenfunction corresponding to the $k$th eigenvalue of the Neumann Laplacian, namely
$$
\left\lbrace
\begin{aligned}
    -\Delta\psi_k^N &= \lambda_k^N\psi_k^N \ \mbox{in } \Omega, \\
    \frac{\partial\psi_k^N}{\partial_\nu} &=0 \ \mbox{on } \partial\Omega,\\
    \psi_k^N &\in H^1(\Omega),
\end{aligned}
\right.
$$
with $0=\lambda_0^N<\lambda_1^N\leq\lambda_2^N\leq\dots$, normalized to make $\{\psi_k^N\}_k$ a normal basis of $L^2(\Omega)$. Then we write
$$
u(x) = \sum_{k=0}^{+\infty} u_k^N\psi_k^N(x),
$$
for any $u \in L^2(\Omega)$, and the corresponding spectral $s$-fractional Laplacian with Neumann boundary conditions is defined as
$$
(-\Delta)^s_{N,\Omega} u(x) \coloneqq \sum_{k=0}^{+\infty} (\lambda_k^N)^s u_k^N\psi_k^N(x).
$$
As proved in \cite[Appendix~P]{abatangelo2019getting}, this can be written as the regional integrodifferential operator
$$
(-\Delta)^s_{N,\Omega} u(x)=\PV\int_{\Omega}(u(x)-u(y))K_\Omega(x,y)\dd y,
$$
where $K_\Omega$ is the kernel in (K4) defined in Section~\ref{sec:kernels}.
Also this operator thus fits into the functional setting of Section~\ref{sec:chem}.


\appendix

\addtocontents{toc}{\protect\setcounter{tocdepth}{1}}


\section{Fractional Poincar\'e inequalities}


\subsection{Proof of Proposition~\ref{prop:onthesofa}}\label{sec:smoothie}

    We rely on a careful covering argument to extend the elementary, well-known proof. We fix $\delta:=\varrho/4$ and we consider a (finite) sequence of open sets with Lipschitz boundary $\Omega_\ell\subset\Omega_{\ell-1}\subset\ldots\subset\Omega_2\subset\Omega\subset\Omega_1$, such that
	$$
	\mbox{diam}(\Omega_\ell)<2\delta\qquad\mbox{and}\qquad\partial\Omega_{i+1}
	\subset\big\{y\in\real^d\,:\,\mbox{dist}(y,\partial\Omega_i)<\delta\big\},\quad\mbox{for every }i=1,\ldots,\ell-1.
	$$
	By compactness, for every $i=1,\ldots,\ell$, we can find a finite number of points $x_1^i,\ldots,x_{k_i}^i\in\partial\Omega_i$ such that
	$$
	\Omega_i\setminus\Omega_{i+1}\subset\bigcup_{j=1}^{k_i} B_{2\delta}(x^i_j),\quad\mbox{for every }i=1,\ldots,\ell-1,
	\qquad\mbox{and}\qquad
	\Omega_\ell \subset\bigcup_{j=1}^{k_\ell} B_{2\delta}(x^\ell_j).
	$$
	Next, we define the sets $A_j^i$ by setting
	$$
	A^i_1:=B_{2\delta}(x^i_1)\cap\big(\Omega_i\setminus\Omega_{i+1}\big)\quad\mbox{and}\quad A^i_j:=B_{2\delta}(x^i_j)\cap\big(\Omega_i\setminus\Omega_{i+1}\big)\setminus\bigcup_{n=1}^{j-1} B_{2\delta}(x^i_n),\quad\mbox{for }j=2,\ldots,k_i,
	$$
	for every $i=1,\ldots,\ell-1$, and
	$$
	A^\ell_1:=B_{2\delta}(x^\ell_1)\cap\Omega_\ell\quad\mbox{and}\quad A^\ell_j:=B_{2\delta}(x^\ell_j)\cap\Omega_\ell\setminus\bigcup_{n=1}^{j-1} B_{2\delta}(x^\ell_n),\quad\mbox{for }j=2,\ldots,k_\ell.
	$$
	We can assume that $|A^i_j|>0$ for every $i,j$. Moreover, by definition, the sets $A^i_j$'s are pairwise disjoint, and
	$$
	\Omega\subset\Omega_1=\Omega_\ell\cup\bigcup_{i=1}^{\ell-1}\big(\Omega_i\setminus\Omega_{i+1}\big)=\bigcup_{i=1}^\ell\bigcup_{j=1}^{k_i} A_j^i.
	$$
	We also remark that, since each $\Omega_i$ is a bounded open set with Lipschitz boundary,
	$$
	c_\ast:=\min_{i,j}\big|B_{2\delta}(x^i_j)\setminus\Omega_i\big|>0.
	$$
	We are now ready to prove \eqref{eq:gen_Poincare}. We start by observing that we can assume the right hand side of \eqref{eq:gen_Poincare} to be finite, otherwise there is nothing to prove. As a first step, since $u=0$ almost everywhere in $\real^d\setminus\Omega$, hence in $\real^d\setminus\Omega_1$, we have
	$$
	|u(x)|^q=\frac{1}{|B_{2\delta}(x^1_j)\setminus\Omega_1|}\int_{B_{2\delta}(x^1_j)\setminus\Omega_1} |u(x)-u(y)|^q\dd y
	\leq\frac{1}{c_\ast}\int_{B_{2\delta}(x^1_j)\setminus\Omega_1} \frac{|u(x)-u(y)|^q}{|x-y|^{d+sq}}|x-y|^{d+sq}\dd y,
	$$
	for every $j=1,\ldots,k_1$, and every $x\in A_j^1$. Thus,
	\begin{align*}
	\|u\|_{L^q(A_j^1)}^q&\leq\frac{\varrho^{d+sq}}{c_\ast}\int_{A^1_j}\int_{B_{2\delta}(x^1_j)\setminus\Omega_1} \frac{|u(x)-u(y)|^q}{|x-y|^{d+sq}}\dd x\dd y\\
	&
	\leq
	\frac{\varrho^{d+sq}}{c_\ast}\int_{A^1_j}\int_{\real^d\cap\{|x-y|<\varrho\}}\frac{|u(x)-u(y)|^q}{|x-y|^{d+sq}}\dd x\dd y,
	\end{align*}
	which implies
	$$
	\|u\|_{L^q(\Omega\setminus\Omega_2)}^q=\|u\|_{L^q(\Omega_1\setminus\Omega_2)}^q=\sum_{j=1}^{k_1}\|u\|_{L^q(A_j^1)}^q\leq
	\frac{\varrho^{d+sq}}{c_\ast}\int_{\Omega_1\setminus\Omega_2}\int_{\real^d\cap\{|x-y|<\varrho\}}\frac{|u(x)-u(y)|^q}{|x-y|^{d+sq}}\dd x\dd y.
	$$
	Next we estimate
	\begin{align*}
		|u(x)|^q&=\frac{1}{|B_{2\delta}(x^2_j)\setminus\Omega_2|}\int_{B_{2\delta}(x^2_j)\setminus\Omega_2} |u(x)-u(y)+u(y)|^q\dd y\\
		&
		\leq \frac{2^{q-1}}{c_\ast}\int_{B_{2\delta}(x^2_j)\setminus\Omega_2} |u(x)-u(y)|^q\dd y+\frac{2^{q-1}}{c_\ast}\|u\|_{L^q(\real^d\setminus\Omega_2)}^q\\
		&
		\leq \frac{2^{q-1}\varrho^{d+sq}}{c_\ast}\int_{B_{2\delta}(x^2_j)\setminus\Omega_2} \frac{|u(x)-u(y)|^q}{|x-y|^{d+sq}}\dd y\\
		&\qquad\qquad
		+\frac{2^{q-1}\varrho^{d+sq}}{c_\ast^2}\int_{\Omega_1\setminus\Omega_2}\int_{\real^d\cap\{|x-y|<\varrho\}}\frac{|u(x)-u(y)|^q}{|x-y|^{d+sq}}\dd x\dd y,
	\end{align*}
	for every $j=1,\ldots,k_2$, and every $x\in A_j^2$. Therefore,
	\begin{align*}
		\|u\|_{L^q(\Omega_2\setminus\Omega_3)}^q&=\sum_{j=1}^{k_2}\|u\|_{L^q(A_j^2)}^q\\
		&
		\leq
		\frac{2^{q-1}\varrho^{d+sq}}{c_\ast}\int_{\Omega_2\setminus\Omega_3}\int_{\real^d\cap\{|x-y|<\varrho\}}\frac{|u(x)-u(y)|^q}{|x-y|^{d+sq}}\dd x\dd y\\
		&\qquad\qquad
		+\frac{|\Omega_2\setminus\Omega_3|2^{q-1}\varrho^{d+sq}}{c_\ast^2}\int_{\Omega_1\setminus\Omega_2}\int_{\real^d\cap\{|x-y|<\varrho\}}\frac{|u(x)-u(y)|^q}{|x-y|^{d+sq}}\dd x\dd y\\
		&
		\leq
		\frac{(|\Omega_2\setminus\Omega_3|+c_\ast)2^{q-1}\varrho^{d+sq}}{c_\ast^2}\int_{\Omega_1\setminus\Omega_3}\int_{\real^d\cap\{|x-y|<\varrho\}}\frac{|u(x)-u(y)|^q}{|x-y|^{d+sq}}\dd x\dd y.
	\end{align*}
	Proceeding in this way we obtain
	\begin{align*}
		\|u\|_{L^q(\Omega_i\setminus\Omega_{i+1})}^q&
		\leq
		\frac{2^{q-1}\varrho^{d+sq}}{c_\ast}\int_{\Omega_i\setminus\Omega_{i+1}}\int_{\real^d\cap\{|x-y|<\varrho\}}\frac{|u(x)-u(y)|^q}{|x-y|^{d+sq}}\dd x\dd y\\
		&
		\qquad\qquad+\frac{2^{q-1}|\Omega_i\setminus\Omega_{i+1}|}{c_\ast}\|u\|_{L^q(\real^d\setminus\Omega_i)}^q\\
		&
		\leq C \int_{\Omega_1\setminus\Omega_{i+1}}\int_{\real^d\cap\{|x-y|<\varrho\}}\frac{|u(x)-u(y)|^q}{|x-y|^{d+sq}}\dd x\dd y,
	\end{align*}
	for every $i=1,\ldots,\ell-1$, and, similarly,
	$$
	\|u\|_{L^q(\Omega_\ell)}^q\leq C \int_{\Omega_1}\int_{\real^d\cap\{|x-y|<\varrho\}}\frac{|u(x)-u(y)|^q}{|x-y|^{d+sq}}\dd x\dd y.
	$$
	Summing these inequalities concludes the proof.


\subsection{Proof of Proposition~\ref{stm:lorenzanonce}}\label{sec:bubbletea}

The proof follows the usual argument by contradiction, which can be found e.\,g. in \cite[Section~5.8.1]{evans2010partial} for the local case and \cite[Lemma~3.10]{dipierro2017nonlocal} for the nonlocal case with the fractional seminorm on the whole $\Omega^2$.

First of all, we observe that we need to prove \eqref{eq:eandatavia} only for those measurable functions $u:\Omega\to\real$ which have finite right-hand side. By Remark~\ref{rmk:luca} any such function belongs to $W^{s,q}(\Omega)$.
We then assume by contradiction that there exist functions $v_k\in W^{s,q}(\Omega)$ such that
\begin{equation}\label{eq:uhm}
\mathfrak m(v_k)=0,\quad\|v_k\|_{L^q(\Omega)}=1,
\end{equation}
and
\begin{equation}\label{eq:mah}
\iint_{\Omega^2\cap\{|x-y|<\varrho\}}\frac{|v_k(x)-v_k(y)|^q}{|x-y|^{d+sq}}\dd x\dd y<\frac{1}{k}.
\end{equation}
As a consequence,
$$
[v_k]_{W^{s,q}(\Omega)}^q\leq\frac{1}{k}+\frac{2^q|\Omega|}{\varrho^{d+sq}}.
$$
Therefore, by the compact embedding of $W^{s,q}(\Omega)$ in $L^q(\Omega)$, there exists a function $v\in W^{s,q}(\Omega)$ such that
$$
v_k\to v\quad\mbox{strongly in }L^q(\Omega)\mbox{ and a.\,e. in }\Omega,
$$
up to a subsequence that we do not relabel.
By \eqref{eq:uhm} and \eqref{eq:mah} we have
$$
\mathfrak m(v)=0,\quad\|v\|_{L^q(\Omega)}=1\quad\mbox{and}\quad\iint_{\Omega^2\cap\{|x-y|<\varrho\}}\frac{|v(x)-v(y)|^q}{|x-y|^{d+sq}}\dd x\dd y=0.
$$
In particular, given any open set $B\subset\Omega$ with $\textup{diam}(B)<\varrho$, we have
\begin{equation}\label{eq:erste}
	[v]_{W^{s,q}(B)}\leq \iint_{\Omega^2\cap\{|x-y|<\varrho\}}\frac{|v(x)-v(y)|^q}{|x-y|^{d+sq}}\dd x\dd y= 0.
\end{equation}
Hence $v$ takes a (possibly different) constant value on each ball $B\subset\Omega$ having radius smaller than $\varrho/2$. Since $\Omega$ is connected, it is then standard to prove that $v$ is constant on the whole $\Omega$.

A possible proof is the following. We set $\delta(y)\coloneqq \min\{\mbox{dist}(y,\partial\Omega),\varrho/2\}$, for every $y\in\Omega$.
Then, \eqref{eq:erste} implies that for every $x\in\Omega$ there exists a constant $c(x)\in\real$ such that
\begin{equation*}
v(y)=c(x)\quad\mbox{for almost every }y\in B_{\delta(x)}(x).
\end{equation*}
Notice that
\begin{equation}\label{eq:hilfe}
c(y)=c(x)\quad\mbox{for every }y\in B_{\delta(x)}(x).
\end{equation}
We then pick a point $x_0\in\Omega$ and we define $c_0\coloneqq c(x_0)$ and
$$
\Omega_0\coloneqq\{x\in\Omega\,:\,c(x)=c_0\}.
$$
As a consequence of \eqref{eq:hilfe}, both $\Omega_0$ and its complement are open sets. Since $\Omega$ is connected, $\Omega_0=\Omega$ and $v$ is constant on the whole $\Omega$. However, this gives a contradiction with the fact that $\mathfrak m(v)=0$ and $\|v\|_{L^q(\Omega)}=1$, concluding the proof.


\subsection{Generalizations for disconnected domains}\label{sec:melange}

Differently from the classical case of the $L^q$-norm of the gradient, in the fractional framework, given the nonlocal nature of the $W^{s,q}$-seminorm, it is still possible to prove some versions of the Poincar\'e inequality when $\Omega$ is not connected.
It is enough to require some compatibility between the geometry of the domain $\Omega$ and the interaction radius $\varrho$.

\begin{proposition}\label{stm:jigglypuff}
Let $s\in(0,1)$, $q\geq1$ and $\Omega\subset\real^d$ be a bounded open set with Lipschitz boundary. If $\Omega_1,\dots,\Omega_j$ are the connected components of $\Omega$ and $\varrho>0$ is such that
 \begin{equation}\label{eq:amo}
 \textup{dist}(\Omega_i,\Omega_{i+1})<\varrho \quad\mbox{for every }i=1,\dots,j-1,
 \end{equation}
then there exists a constant $C=C(d,\Omega,s,q,\varrho)>0$ such that
\begin{equation*}
\|u-\mathfrak m(u)\|_{L^q(\Omega)}^q\leq C\iint_{\Omega^2\cap\{|x-y|<\varrho\}}\frac{|u(x)-u(y)|^q}{|x-y|^{d+sq}}\dd x\dd y,
\end{equation*}
for every measurable $u:\Omega\to\real$.
\end{proposition}

\begin{proof}
By arguing by contradiction as in the proof of Proposition~\ref{stm:lorenzanonce}, we end up with a function $v\in W^{s,q}(\Omega)$ such that
$$
\mathfrak m(v)=0,\quad\|v\|_{L^q(\Omega)}=1\quad\mbox{and}\quad\iint_{\Omega^2\cap\{|x-y|<\varrho\}}\frac{|v(x)-v(y)|^q}{|x-y|^{d+sq}}\dd x\dd y=0,
$$
which implies that $v$ is constant on the connected components,
$$
v\equiv c_i\quad\mbox{in }\Omega_i.
$$
Since $\textup{dist}(\Omega_i,\Omega_{i+1})<\varrho$, we can find a point $x_i\in\real^d$ such that
$$
|\Omega_i\cap B_{\varrho/2}(x_i)|>0
\quad\mbox{and}\quad |\Omega_{i+1}\cap B_{\varrho/2}(x_i)|>0.
$$
Then,
$$
\int_{\Omega_i\cap B_{\varrho/2}(x_i)}\int_{\Omega_{i+1}\cap B_{\varrho/2}(x_i)}\frac{|c_i-c_{i+1}|^q}{|x-y|^{d+sq}}\dd x\dd y
\leq \iint_{\Omega^2\cap\{|x-y|<\varrho\}}\frac{|v(x)-v(y)|^q}{|x-y|^{d+sq}}\dd x\dd y=0,
$$
hence $c_{i+1}=c_i$. Therefore $c_j=c_{j-1}=\dots=c_1$ and $v$ is constant on the whole $\Omega$.
This gives a contradiction with the fact that $\mathfrak m(v)=0$ and $\|v\|_{L^q(\Omega)}=1$, concluding the proof.
\end{proof}
We can interpret Proposition~\ref{stm:lorenzanonce} as the particular case of Proposition~\ref{stm:jigglypuff} in which $\Omega$ is connected and $\varrho$ is arbitrary.

We stress that a condition like \eqref{eq:amo} or the connectedness of $\Omega$ is necessary to ensure the validity of the fractional Poincar\'e inequality, as shown by the following example.

\begin{example}
    We consider $\Omega\coloneqq B_1(-4e_1)\cup B_1(4e_1)$ and $\varrho=2$. Then, we define
    $$
    v_k\coloneqq k\chi_{B_1(-4e_1)}-k\chi_{B_1(4e_1)}.
    $$
    Then,
    $$
    \mathfrak m(v_k)=0,\quad\|v_k\|_{L^q(\Omega)}^q=2k^q|B_1|
    $$
    and
    $$
    \iint_{\Omega^2\cap\{|x-y|<2\}}\frac{|v_k(x)-v_k(y)|^q}{|x-y|^{d+sq}}\dd x\dd y=[k]_{W^{s,q}(B_1(-4e_1))}^q+[-k]_{W^{s,q}(B_1(4e_1))}^q=0.
    $$
    This proves that \eqref{eq:eandatavia}
    cannot hold true in this case.
\end{example}

On the other hand, the Lipschitz regularity of the boundary of $\Omega$ was required only to ensure the validity of the compact embedding of $W^{s,q}(\Omega)$ in $L^q(\Omega)$, which was used in the above proof by contradiction. In the case in which $\varrho\geq\textup{diam}(\Omega)$ this hypothesis can be dropped. Indeed, the argument leading to formula (8.3) in \cite{di2012hitchhiker} provides a direct proof of the following result.
\begin{proposition}\label{stm:poinconnected}
Let $s\in(0,1),\,q\geq1$ and $\Omega\subset\real^d$ be a bounded open set.
Then,
\begin{equation*}
\|u-\mathfrak m(u)\|_{L^q(\Omega)}^q\leq \frac{\textup{diam}(\Omega)^{d+sq}}{|\Omega|}\iint_{\Omega^2}\frac{|u(x)-u(y)|^q}{|x-y|^{d+sq}}\dd x\dd y,
\end{equation*}
for every measurable $u:\Omega\to\real$.
\end{proposition}

Notice that in the statement of the above proposition $\Omega$ is not assumed to be connected.


\section{Further properties of the fractional spaces}


\subsection{A characterization of the integrability of the kernels}

We prove that the integrability of the kernel is equivalent to requiring Lipschitz functions to have finite energy.

\begin{lemma}\label{lem:AloeNatale}
    Let $K : \real^d \times \real^d \to [0,+\infty]$ be Borel-measurable. Then condition \eqref{eq:Kintegrability} holds true if and only if
    $$
    \iint_{Q(\Omega)}|u(x)-u(y)|^q K(x,y)\dd x\dd y < +\infty
    $$
    for every $u : \real^d \to \real$ which is uniformly Lipschitz and bounded.
\end{lemma}

\begin{proof}
    Since $\Omega$ is bounded, there exists $R\geq 1$ such that $\Omega\subset B_R$. We consider the functions $u_i : \real^d \to \real$ defined as $u_i\coloneqq\min\{x_i,R+1\}$, which are uniformly Lipschitz and bounded, for every $i=1,\dots,d$. We claim that
    $$
    \min\{1,|x-y|^q\}\leq d\max_{1\leq i\leq d}|u_i(x)-u_i(y)|^q
    $$
    for every $(x,y)\in Q(\Omega)$. Indeed, if $x,y\in B_{R+1}$, then
    $$
    \min\{1,|x-y|^q\}\leq |x-y|^q\leq d\max_{1\leq i\leq d}|x_i-y_i|^q=d\max_{1\leq i\leq d}|u_i(x)-u_i(y)|^q.
    $$
    On the other hand, if $x\in\Omega$ and $y\in \real^d\setminus B_{R+1}$, then
    $$
    \min\{1,|x-y|^q\}=1\leq d\,|R-R-1|^q\leq d\max_{1\leq i\leq d}|u_i(x)-u_i(y)|^q,
    $$
    and similarly for $x\in \real^d\setminus B_{R+1}$ and $y\in \Omega$. Hence, recalling also Remark~\ref{rmk:Kqcompact},
    $$
    \begin{aligned}
    \iint_{Q(\Omega)}\min &\{1,|x-y|^q\} K(x,y)\dd x\dd y\leq d\max_{1\leq i\leq d}\iint_{Q(\Omega)}|u_i(x)-u_i(y)|^q K(x,y)\dd x\dd y\\
    &
    \leq d\max_{1\leq i\leq d}\Big([u_i]_{C^{0,1}(\real^d)}^q+2^q\|u_i\|^q_{L^\infty(\real^d)}\Big)\iint_{Q(\Omega)}\min\{1,|x-y|^q\} K(x,y)\dd x\dd y.
    \end{aligned}
    $$
    Therefore, requiring the $\|\,\cdot\,\|_{K,q}$-seminorm to be finite for every function $u : \real^d \to \real$ which is uniformly Lipschitz and bounded is indeed equivalent to condition \eqref{eq:Kintegrability}. Actually, this is equivalent to requiring the functions $u_i$ to have finite $\|\,\cdot\,\|_{K,q}$-seminorm.
\end{proof}


\subsection{An approximation result}\label{sec:approximation}

We divide the proof of Lemma~\ref{lem:approximation} in several steps.

\textbf{Step 1.} Without loss of generality we assume that $\varphi\in\cH^s_\Omega\cap L^\infty(\real^d)$. Indeed, if we define the truncated functions $\varphi_k\coloneqq \max\{-k,\min\{\varphi,k\}\}$, then
    $$
    |\varphi_k|^2\leq|\varphi|^2\quad\mbox{and}\quad\frac{|\varphi_k(x)-\varphi_k(y)|^2}{|x-y|^{d+2s}}\leq \frac{|\varphi(x)-\varphi(y)|^2}{|x-y|^{d+2s}},
    $$
    hence
    $$
    \lim_{k\to\infty}\big(\|\varphi_k-\varphi\|_{L^2(B_R)}+\|\varphi_k-\varphi\|_{\cH^s_\Omega}\big)=0\quad\mbox{for every }R>0.
    $$

\textbf{Step 2.} We observe that, if $\varphi\in\cH^s_\mathcal O$ for some open set $\mathcal O\subset\real^d$ such that $\Omega\subset\subset\mathcal O$, and $\{\eta_\varepsilon\}$ is a standard sequence of mollifiers, then
$$
    \lim_{\varepsilon\to0}\big(\|\varphi\ast\eta_\varepsilon -\varphi\|_{L^2(B_R)}+\|\varphi\ast\eta_\varepsilon-\varphi\|_{\cH^s_\Omega}\big)=0\quad\mbox{for every }R>0.
    $$
    For the proof, see, e.\,g., the proof of \cite[Lemma~11]{fiscella2015density} or \cite[Lemma~3.2\,(i)]{lombardini2018approximation}.

    \textbf{Step 3.} We construct a sequence of bounded open sets $\mathcal O_\delta\subset\real^d$ with $\Omega\subset\subset\mathcal O_\delta$ and a sequence of functions $\varphi_\delta\in\cH^s_{\mathcal O_\delta}$, such that
    $$
    \begin{aligned}
        & \varphi_\delta\to \varphi\quad\mbox{a.\,e. in }\real^d,\\
        \sup_{\delta\in(0,\delta_0)} \|\varphi_\delta\|_{L^\infty(\real^d)}\leq \|\varphi&\|_{L^\infty(\real^d)}\quad\mbox{and}\quad\sup_{\delta\in(0,\delta_0)} \|\varphi_\delta\|_{\cH^s_{\mathcal O_\delta}}\leq C\|\varphi\|_{\cH^s_\Omega},
    \end{aligned}
    $$
    for some $\delta_0,C>0$.
    
 Since $\partial\Omega$ is of class $C^2$, we know (see, e.\,g., \cite[Lemma~14.16]{gilbarg2001elliptic} and \cite[Theorem~2]{ambrosio2000geometric}) that there exists $r_0>0$ such that $\bar{d}_\Omega\in C^2(N_{3r_0}(\partial\Omega))$, where
 $$
 \bar{d}_\Omega(x)\coloneqq\textup{dist}(x,\Omega)-\textup{dist}(x,\real^d\setminus\Omega)\quad\mbox{and}\quad N_\varrho(\partial\Omega)\coloneqq\{y\in\real^d\,:\,\textup{dist}(y,\partial\Omega)<\varrho\}.
 $$
We consider a cut-off function $\vartheta\in C^\infty_c(\real^d)$ such that $0\leq\vartheta\leq 1$, $\vartheta\equiv 1$ in $N_{r_0}(\partial\Omega)$ and $\vartheta\equiv 0$ in $\real^d\setminus N_{2r_0}(\partial\Omega)$, and we define
$\Theta_\delta:\real^d\to\real^d$ as
$$
\Theta_\delta(x)=x+\delta\vartheta(x)\nabla\bar{d}_\Omega(x),
$$
    for $\delta>0$. For $\delta$ small enough this is a $C^1$-diffeomorphism of $\real^d$, such that
    $$
    \textup{supp}(\textup{Id}-\Theta_\delta)\subset N_{2r_0}(\partial\Omega)\quad\mbox{and}\quad \Omega\subset\subset N_\delta(\Omega)\subset\Theta_\delta(\Omega).
    $$
The second statement follows from the fact that
$$
\nabla\bar{d}_\Omega(x)=\nabla\bar{d}_\Omega(\wp(x))=\nu_\Omega(\wp(x)),
$$
    where $\wp(x)\in\partial\Omega$ is the point of minimal distance between $x\in N_{2r_0}(\partial\Omega)$ and $\partial\Omega$, and $\nu_\Omega$ is the outer unit normal (see, e.\,g., \cite[Remark~3\,(1) and Theorem~2\,(i)]{ambrosio2000geometric}).

We then define the function $\varphi_\delta\coloneqq\varphi\circ\Theta_\delta^{-1}$ which belongs to $L^\infty(\real^d)$, since we clearly have $\|\varphi_\delta\|_{L^\infty(\real^d)}\leq \|\varphi\|_{L^\infty(\real^d)}$. Moreover, $\Theta_\delta^{-1}(x)\to x$ uniformly in $x\in\real^d$ as $\delta\to0$, hence $\varphi_\delta\to\varphi$ almost everywhere in $\real^d$. Furthermore
$$
\begin{aligned}
\|\varphi_\delta\|^2_{\cH^s_{N_\delta(\Omega)}}&\leq \|\varphi_\delta\|_{\cH^s_{\Theta_\delta(\Omega)}}^2=\frac{1-s}{2}\iint_{Q(\Theta_\delta(\Omega))}\frac{|(\varphi\circ\Theta_\delta^{-1})(x)-(\varphi\circ\Theta_\delta^{-1})(y)|^2}{|x-y|^{d+2s}}\dd x\dd y\\
&
=\frac{1-s}{2}\iint_{Q(\Omega)}\frac{|\varphi(\xi)-\varphi(\zeta)|^2}{|\Theta_\delta(\xi)-\Theta_\delta(\zeta)|^{d+2s}}|\det D\Theta_\delta(\xi)||\det D\Theta_\delta(\zeta)|\dd \xi\dd \zeta\\
&
\leq \frac{1-s}{2}\big(\textup{Lip}(\Theta_\delta^{-1})\big)^{d+2s}\iint_{Q(\Omega)}\frac{|\varphi(\xi)-\varphi(\zeta)|^2}{|\xi-\zeta|^{d+2s}}|\det D\Theta_\delta(\xi)||\det D\Theta_\delta(\zeta)|\dd \xi\dd \zeta\\
&
\leq C\|\varphi\|_{\cH^s_\Omega}^2,
\end{aligned}
$$
with $C$ uniform in $\delta$ for $\delta$ small enough, proving the claim.

\textbf{Step 4.} By Step 3 we can find a sequence $\delta_k\searrow0$ for which the functions $\varphi_{\delta_k}=\varphi\circ\Theta_{\delta_k}^{-1}\in\cH^s_{N_{\delta_k}(\Omega)}$ are such that
    $$
    \|\varphi_{\delta_k}-\varphi\|_{L^2(B_k)}<\frac{1}{2k},\quad\sup_k\|\varphi_{\delta_k}\|_{L^\infty(\real^d)}\leq \|\varphi\|_{L^\infty(\real^d)}\quad\mbox{and}\quad
    \sup_k\|\varphi_{\delta_k}\|_{\cH^s_{N_{\delta_k}(\Omega)}}\leq C \|\varphi\|_{\cH^s_\Omega}.
    $$
    We can then apply Step 2 to each $\varphi_{\delta_k}$ and find a sequence $\varepsilon_k\searrow0$ such that
    $$
    \|\varphi_{\delta_k}\ast\eta_{\varepsilon_k}-\varphi_{\delta_k}\|_{L^2(B_k)}<\frac{1}{2k}
    \quad\mbox{and}\quad
    \|\varphi_{\delta_k}\ast\eta_{\varepsilon_k}-\varphi_{\delta_k}\|_{\cH^s_\Omega}<\frac{1}{k}.
    $$
    Moreover, by standard properties of regularization by mollification, we have $\varphi_{\delta_k}\ast\eta_{\varepsilon_k}\in C^\infty(\real^d)\cap L^\infty(\real^d)$, with
    $$
    \sup_k\|\varphi_{\delta_k}\ast\eta_{\varepsilon_k}\|_{L^\infty(\real^d)}\leq \sup_k\|\varphi_{\delta_k}\|_{L^\infty(\real^d)}\leq \|\varphi\|_{L^\infty(\real^d)}.
    $$
    Therefore, if we denote $\varphi_k\coloneqq\varphi_{\delta_k}\ast\eta_{\varepsilon_k}$, we obtain
    $$
    \|\varphi_k-\varphi\|_{L^2(B_k)}<\frac{1}{k}\quad\mbox{and}\quad \sup_k\|\varphi_k\|_{\cH^s_\Omega}\leq C \|\varphi\|_{\cH^s_\Omega}+1.
    $$
    The strong $L^2_\textup{loc}(\real^d)$ convergence implies that $\varphi_k\to\varphi$ almost everywhere in $\real^d$ (up to extracting a subsequence, that we do not relabel).

    We are left to prove the third statement of \eqref{eq:ennui}. For this, we observe that the function
    $$
    \xi_k(x,y)\coloneqq\frac{\varphi_k(x)-\varphi_k(y)}{|x-y|^{\frac{d}{2}+s}}
    $$
    is such that
    $$
    \bigg(\frac{1-s}{2}\bigg)^\frac{1}{2}\sup_k\|\xi_k\|_{L^2(Q(\Omega))}=\sup_k\|\varphi_k\|_{\cH^s_\Omega}\leq C \|\varphi\|_{\cH^s_\Omega}+1,
    $$
    hence $\xi_k\rightharpoonup\xi$ weakly in $L^2(Q(\Omega))$, for some $\xi\in L^2(Q(\Omega))$. The pointwise convergence of $\varphi_k$ ensures that
    $$
    \xi(x,y)=\frac{\varphi(x)-\varphi(y)}{|x-y|^{\frac{d}{2}+s}}.
    $$
    Thus
    $$
    \iint_{Q(\Omega)}\frac{(v(x)-v(y))(\varphi_k(x)-\varphi_k(y))}{|x-y|^{d+2s}}\dd x\dd y\to \iint_{Q(\Omega)}\frac{(v(x)-v(y))(\varphi(x)-\varphi(y))}{|x-y|^{d+2s}}\dd x\dd y,
    $$
    for every $v\in\cH^s_\Omega$, concluding the proof.


\section{Moreau-Yosida regularization}\label{sec:MYreg}

\begin{proof}[Proof of Lemma~\ref{lemma:scoponescientifico}]
In order to prove (1), we begin by observing that
$$
\Gamma_\lambda(r)\geq\inf_{z\in\real}\left\{-\Pi(z)-a_1|z|^p-a_2+\frac{1}{2\lambda}|z-r|^2\right\}.
$$
Next we show that the infimum is achieved in a neighborhood of $r$, i.\,e., there exists $\delta=\delta(r,\lambda)>0$ such that
\begin{equation}\label{eq:killer_koala}
\inf_{z\in\real}\left\{-\Pi(z)-a_1|z|^p-a_2+\frac{1}{2\lambda}|z-r|^2\right\}
=
\inf_{|z-r|<\delta}\left\{-\Pi(z)-a_1|z|^p-a_2+\frac{1}{2\lambda}|z-r|^2\right\}.
\end{equation}
Since
$$
C_\pi r^2>-\Pi(r)-a_1|r|^p-a_2\geq \inf_{z\in\real}\left\{-\Pi(z)-a_1|z|^p-a_2+\frac{1}{2\lambda}|z-r|^2\right\},
$$
in order to prove \eqref{eq:killer_koala}, it is enough to show that for $|z-r|\geq\delta$ we have
\begin{align*}
-\Pi(z)-a_1|z|^p-a_2+\frac{1}{2\lambda}|z-r|^2\geq C_\pi r^2.
\end{align*}
For this, we point out that
\begin{align*}
-\Pi(z)-a_1|z|^p-a_2&+\frac{1}{2\lambda}|z-r|^2- C_\pi r^2
\geq
-\Big(a_1+\frac{C_\pi}{2}\Big)z^2+\frac{1}{2\lambda}|z-r|^2- C_\pi r^2-a_1-a_2\\
&
\geq \Big(\frac{1}{2\lambda}-C_\pi-2a_1\Big)|z-r|^2-\big(4C_\pi+2a_1\big) r^2-a_1-a_2,
\end{align*}
where we have used Young's inequality in the last estimate, in the following form
$$
rz\leq\frac{1}{4}z^2+r^2.
$$
This proves \eqref{eq:killer_koala} with
$$
\delta\coloneqq\lambda^\frac{1}{2}\sqrt{\frac{2\big(4C_\pi+2a_1\big) r^2+2a_1+2a_2}{1-2\lambda \big(C_\pi+2a_1\big)}}
\leq 2\lambda^\frac{1}{2}\big(|r|2\sqrt{2C_\pi+a_1}+\sqrt{2a_1+2a_2}\big)\leq c\lambda^\frac{1}{2}\big(|r|+1\big),
$$
with $\lambda<3/(8 C_\pi+16a_1)$ and $c=c(a_1,a_2,C_\pi)>0$.
Therefore
\begin{align*}
\sup_{|z-r|<\delta}a_1|z|^p
\leq a_1\Big(|r|\big(c\lambda^\frac{1}{2}+1\big)+c\lambda^\frac{1}{2}\Big)^p
\leq a_1\big(|r|(c+1)+c\big)^p
\leq 2 a_1 (|r|^p(c+1)^p+c^p),
\end{align*}
hence
$$
\inf_{|z-r|<\delta}-a_1|z|^p
\geq
-a_3|r|^p-a_4\qquad\mbox{for every }r\in\real,
$$
with $a_3,a_4>0$ depending only on $a_1,a_2$ and $C_\pi$.
Thus
\begin{align*}
	\Gamma_\lambda(r)+\Pi(r)&\geq \inf_{|z-r|<\delta}\left\{-\Pi(z)-a_1|z|^p-a_2+\Pi(r)+\frac{1}{2\lambda}|z-r|^2\right\}\\
    &
    \geq -a_3|r|^p-a_4-a_2
	+
	\inf_{|z-r|<\delta}\left\{\Pi(r)-\Pi(z)\right\}\\
	&
	= -a_3|r|^p-a_4-a_2+\inf_{|z-r|<\delta}\int_z^r \pi(\zeta)\dd\zeta\\
    &
    \geq-a_3|r|^p-a_4-a_2-\delta C_\pi\big(|r|(c+1)+c\big)\\
    &
    \geq -a_3|r|^p-\beta -\alpha\lambda^\frac{1}{2} r^2,
\end{align*}
which is the desired estimate.

As for point (2), reasoning as above we are led to
$$
\Gamma_\lambda(r)+\Pi(r)\geq \inf_{|z-r|<\delta}\left\{-\Pi(z)-a_1 z^2-a_2+\Pi(r)+\frac{1}{2\lambda}|z-r|^2\right\},
$$
where
$$
\delta\coloneqq\lambda^\frac{1}{2}\sqrt{\frac{2\big(4C_\pi+2a_1\big) r^2+2a_2}{1-2\lambda \big(C_\pi+2a_1\big)}}
\leq c\lambda^\frac{1}{2}\big(|r|+1\big),
$$
with $c=c(a_1,a_2,C_\pi)>0$ and $\lambda>0$ small enough. Since
\begin{equation}\label{eq:linkpostali}
    \inf_{|z-r|<\delta}-a_1 z^2
    \geq
    -2a_1 r^2-4a_1c^2\lambda r^2-4a_1 c^2\lambda\qquad\mbox{for every }r\in\real,
\end{equation}
we conclude that
$$
\Gamma_\lambda(r)+\Pi(r)
\geq -2a_1r^2-\beta -\alpha\lambda^\frac{1}{2} r^2
$$
as claimed.
\end{proof}

\begin{proof}[Proof of Lemma~\ref{lemma:parzialmenteserena}]
We begin by recalling that the $\Gamma$-liminf is defined as
$$
\Big(\Gamma\textrm{-}\liminf_{\lambda\to0}\Gamma_\lambda\Big)(r_0)\coloneqq\inf\Big\{\liminf_{k\to\infty}\Gamma_{\lambda_k}(r_k)\,:\,\lambda_k\to0^+,\ r_k\to r_0\Big\}.
$$

\textbf{Inequality} $\leq$.
Let $\lambda_k\to0$ and $r_k\to r_0$. By the lower semicontinuity of $\Gamma$, for every $\varepsilon>0$ there exists $\delta>0$ such that
$$
\min\Big\{\Gamma(r_0)-\varepsilon,\frac{1}{\varepsilon}\Big\}\leq\Gamma(r)\qquad\mbox{for every }r\in\real\mbox{ such that }|r-r_0|<2\delta.
$$
Notice that
$$
\min\Big\{\Gamma(r_0)-\varepsilon,\frac{1}{\varepsilon}\Big\}
 \ \nearrow \ 
\Gamma(r_0)\qquad\mbox{as }\varepsilon\searrow0,
$$
finite or not.
Moreover, there exists $k_0$ such that $|r_k-r_0|<\delta$ for every $k>k_0$. If $|r-r_0|<2\delta$, then
$$
\Gamma(r)+\frac{1}{2\lambda_k}|r-r_k|^2\geq \min\Big\{\Gamma(r_0)-\varepsilon,\frac{1}{\varepsilon}\Big\},
$$
whereas, if $|r-r_0|\geq 2\delta$, then
$$
\Gamma(r)+\frac{1}{2\lambda_k}|r-r_k|^2\geq \frac{1}{2\lambda_k}\delta^2.
$$
Therefore, taking $k$ big enough,
$$
\Gamma_{\lambda_k}(r_k)=\inf_{r\in\real}\left\{\Gamma(r)+\frac{1}{2\lambda}|r-r_k|^2\right\}\geq \min\Big\{\Gamma(r_0)-\varepsilon,\frac{1}{\varepsilon}\Big\}.
$$
Since this holds for every $\varepsilon>0$, we have proved that
$$
\Gamma(r_0)\leq\Big(\Gamma\textrm{-}\liminf_{\lambda\to0}\Gamma_\lambda\Big)(r_0).
$$

\textbf{Inequality} $\geq$. 
Since $\Gamma_\lambda\leq\Gamma$, by considering $r_k\equiv r_0$, we have
$$
\Big(\Gamma\textrm{-}\liminf_{\lambda\to 0}\Gamma_{\lambda}\Big)(r_0)\leq\liminf_{k\to\infty}\Gamma_{\lambda_k}(r_k)\le\Gamma(r_0),
$$
for every sequence $\lambda_k\to 0$.
This concludes the proof of the lemma.
\end{proof}

\begin{proof}[Alternative proof of~\eqref{eq:barba}]
Let us define $\zeta^\lambda \coloneqq \gamma_\lambda(u^\lambda)$. Take any $\xi,\eta \in L^2(0,T;\cL^2_0)$ such that $\xi(x,t) \in \gamma(\eta(x,t))$ for almost every $(x,t) \in \Omega \times (0,T)$. We can rewrite
\begin{equation}\label{eq:lambda}
\begin{aligned}
&\int_0^T \io (\zeta^\lambda(x,t)-\eta(x,t))(u^\lambda(x,t)-\xi(x,t)) \dd x \dd t \\
&= \int_0^T \io(\zeta^\lambda(x,t)-\eta(x,t))(u^\lambda(x,t)-J_\lambda(u^\lambda(x,t))) \dd x \dd t \\
& \qquad + \int_0^T \io(\zeta^\lambda(x,t)-\eta(x,t))(J_\lambda(u^\lambda(x,t))-\xi(x,t)) \dd x \dd t.
\end{aligned}
\end{equation}
By the definition of Yosida approximation, the first term at the right-hand side of \eqref{eq:lambda} is such that
$$
\int_0^T \io(\zeta^\lambda-\eta)(u^\lambda-J_\lambda(u^\lambda)) \dd x \dd t = \lambda \int_0^T \io(\zeta^\lambda-\eta)\zeta^\lambda \dd x \dd t.
$$
Moreover, Lemma~\ref{lemma:yosida}\,(i) yields $\zeta^\lambda(x,t) \in \gamma(J_\lambda(u^\lambda(x,t)))$ for almost every $(x,t) \in \Omega \times (0,T)$, so that the second term at the right-hand side of \eqref{eq:lambda} satisfies
$$
\int_0^T \io(\zeta^\lambda-\eta)(J_\lambda(u^\lambda)-\xi) \dd x \dd t \ge 0
$$
by the monotonicity of $\gamma$. Therefore, from \eqref{eq:lambda} we deduce
$$
\int_0^T \io (\zeta^\lambda-\eta)(u^\lambda-\xi) \dd x \dd t \ge \lambda \int_0^T \io(\zeta^\lambda-\eta)\zeta^\lambda \dd x \dd t,
$$
that is,
\begin{equation}\label{eq:lambda_ineq}
\int_0^T \io (\zeta^\lambda-\eta)(u^\lambda-u) \dd x \dd t + \int_0^T \io (\zeta^\lambda-\eta)(u-\xi) \dd x \dd t \ge \lambda \int_0^T \io(\zeta^\lambda-\eta)\zeta^\lambda \dd x \dd t.
\end{equation}
By H\"older's inequality, the first term at the left-hand side is such that
\begin{align*}
&\left|\int_0^T \io (\zeta^\lambda(x,t)-\eta(x,t))(u^\lambda(x,t)-u(x,t)) \dd x \dd t\right| \\
&\le \left(\int_0^T \io |\zeta^\lambda(x,t)-\eta(x,t)|^2 \dd x \dd t\right)^{1/2} \left(\int_0^T \io |u^\lambda(x,t)-u(x,t)|^2 \dd x \dd t\right)^{1/2} \\
&\le CT^{1/2} \left(\sup_{t \in (0,T)} \|u^\lambda(t)-u(t)\|_{\cL^2_0}^2\right)^{1/2} \xrightarrow{\lambda \to 0} 0,
\end{align*}
whereas for the second one we directly have, by the weak convergence of $\zeta^\lambda$,
\begin{align*}
&\int_0^T \io (\zeta^\lambda(x,t)-\eta(x,t))(u(x,t)-\xi(x,t)) \dd x \dd t \\
&\xrightarrow{\lambda \to 0} \int_0^T \io (\zeta(x,t)-\eta(x,t))(u(x,t)-\xi(x,t)) \dd x \dd t.
\end{align*}
Finally, again by H\"older's inequality, the right-hand side of \eqref{eq:lambda_ineq} is such that
$$
\lambda\left|\int_0^T \io (\zeta^\lambda(x,t)-\eta(x,t))\zeta^\lambda(x,t) \dd x \dd t\right| \le \lambda C \xrightarrow{\lambda \to 0} 0.
$$
Hence, passing to the limit as $\lambda \to 0$ in \eqref{eq:lambda_ineq} gives
$$
\int_0^T \io (\zeta(x,t)-\eta(x,t))(u(x,t)-\xi(x,t)) \dd x \dd t \ge 0.
$$
The maximal monotonicity of $\gamma$ then yields $\zeta(x,t) \in \gamma(u(x,t))$ for almost every $(x,t) \in \Omega \times (0,T)$.
\end{proof}


\section*{Acknowledgments}
This research was funded by the Austrian Science Fund (FWF) projects \href{https://doi.org/10.55776/F65}{10.55776/F65}, \href{https://doi.org/10.55776/V662}{10.55776/V662}, \href{https://doi.org/10.55776/Y1292}{10.55776/Y1292}, and \href{https://doi.org/10.55776/P35359}{10.55776/P35359}. Support from the Austrian Federal Ministry of Education, Science and Research (BMBWF) through the OeAD-WTZ project CZ09/2023 is also gratefully acknowledged. Most of this work was done while CG was affiliated with TU Wien.


\bibliographystyle{siam}
\bibliography{biblioCH}

\begin{thebibliography}{10}

\bibitem{abatangelo2020remark}
{\sc N.~Abatangelo}, {\em A remark on nonlocal {Neumann} conditions for the
  fractional {Laplacian}}, Archiv der Mathematik, 114 (2020), pp.~699--708.

\bibitem{abatangelo2017nonhomogeneous}
{\sc N.~Abatangelo and L.~Dupaigne}, {\em Nonhomogeneous boundary conditions
  for the spectral fractional {Laplacian}}, Annales de l'Institut Henri
  Poincar{\'e} C, Analyse non lin{\'e}aire, 34 (2017), pp.~439--467.

\bibitem{abatangelo2018loss}
{\sc N.~Abatangelo, S.~Jarohs, and A.~Salda{\~n}a}, {\em On the loss of maximum
  principles for higher-order fractional {Laplacians}}, Proceedings of the
  American Mathematical Society, 146 (2018), pp.~4823--4835.

\bibitem{abatangelo2019getting}
{\sc N.~Abatangelo and E.~Valdinoci}, {\em Getting acquainted with the
  fractional {L}aplacian}, Contemporary research in elliptic PDEs and related
  topics,  (2019), pp.~1--105.

\bibitem{abels2015cahn}
{\sc H.~Abels, S.~Bosia, and M.~Grasselli}, {\em {Cahn-Hilliard} equation with
  nonlocal singular free energies}, Annali di Matematica Pura ed Applicata
  (1923-), 194 (2015), pp.~1071--1106.

\bibitem{akagi2016fractional}
{\sc G.~Akagi, G.~Schimperna, and A.~Segatti}, {\em Fractional {Cahn-Hilliard},
  {Allen-Cahn} and porous medium equations}, Journal of Differential Equations,
  261 (2016), pp.~2935--2985.

\bibitem{akagi2019convergence}
\leavevmode\vrule height 2pt depth -1.6pt width 23pt, {\em Convergence of
  solutions for the fractional {Cahn-Hilliard} system}, Journal of Functional
  Analysis, 276 (2019), pp.~2663--2715.

\bibitem{ambrosio2000geometric}
{\sc L.~Ambrosio}, {\em Geometric evolution problems, distance function and
  viscosity solutions}, in Calculus of Variations and Partial Differential
  Equations, Springer Berlin, 2000, pp.~5--93.

\bibitem{ambrosio2005gradient}
{\sc L.~Ambrosio, N.~Gigli, and G.~Savar{\'e}}, {\em Gradient flows: in Metric
  Spaces and in the Space of Probability Measures}, Springer Science \&
  Business Media, 2005.

\bibitem{audrito2022neumann}
{\sc A.~Audrito, J.-C. Felipe-Navarro, and X.~Ros-Oton}, {\em The {Neumann}
  problem for the fractional {Laplacian}: regularity up to the boundary},
  Annali della Scuola Normale Superiore di Pisa, Classe di Scienze, XXIV
  (2022), pp.~1155--1222.

\bibitem{barbu2010nonlinear}
{\sc V.~Barbu}, {\em Nonlinear Differential Equations of Monotone Types in
  Banach Spaces}, Springer Monographs in Mathematics, Springer New York, 2010.

\bibitem{baroncini2018continuity}
{\sc C.~Baroncini, J.~F. Bonder, and J.~F. Spedaletti}, {\em Continuity results
  with respect to domain perturbation for the fractional $p$-{Laplacian}},
  Applied Mathematics Letters, 75 (2018), pp.~59--67.

\bibitem{barrios2020neumann}
{\sc B.~Barrios, L.~Montoro, I.~Peral, and F.~Soria}, {\em Neumann conditions
  for the higher order $s$-fractional {Laplacian} $(-{\Delta})^s u$ with
  $s>1$}, Nonlinear Analysis, 193 (2020).

\bibitem{bates2005dirichlet}
{\sc P.~W. Bates and J.~Han}, {\em The {Dirichlet} boundary problem for a
  nonlocal {Cahn-Hilliard} equation}, Journal of Mathematical Analysis and
  Applications, 311 (2005), pp.~289--312.

\bibitem{bogdan2003censored}
{\sc K.~Bogdan, K.~Burdzy, and Z.-Q. Chen}, {\em Censored stable processes},
  Probability theory and related fields, 127 (2003), pp.~89--152.

\bibitem{bourgain2001another}
{\sc J.~Bourgain, H.~Brezis, and P.~Mironescu}, {\em Another look at {Sobolev}
  spaces}, in Optimal control and partial differential equations, IOS,
  Amsterdam, 2001, pp.~439--455.

\bibitem{brasco2015stability}
{\sc L.~Brasco, E.~Parini, and M.~Squassina}, {\em Stability of variational
  eigenvalues for the fractional $p$-{Laplacian}}, Discrete and Continuous
  Dynamical Systems, 36 (2016), pp.~1813--1845.

\bibitem{brezis1973operateurs}
{\sc H.~Brezis}, {\em Opérateurs Maximaux Monotones et semi-groupes de
  contractions dans les espaces de Hilbert}, North Holland/American Elsevier,
  1973.

\bibitem{caffarelli2010nonlocal}
{\sc L.~A. Caffarelli, J.-M. Roquejoffre, and O.~Savin}, {\em Nonlocal minimal
  surfaces}, Communications on Pure and Applied Mathematics, 9 (2010),
  pp.~1111--1144.

\bibitem{cahn1958free}
{\sc J.~W. Cahn and J.~E. Hilliard}, {\em Free energy of a nonuniform system.
  {I}. {Interfacial} free energy}, The Journal of Chemical Physics, 28 (1958),
  pp.~258--267.

\bibitem{colli2012global}
{\sc P.~Colli, S.~Frigeri, and M.~Grasselli}, {\em Global existence of weak
  solutions to a nonlocal {Cahn-Hilliard-Navier-Stokes} system}, Journal of
  Mathematical Analysis and Applications, 386 (2012), pp.~428--444.

\bibitem{colli2019well}
{\sc P.~Colli, G.~Gilardi, and J.~Sprekels}, {\em Well-posedness and regularity
  for a generalized fractional {Cahn-Hilliard} system}, Rendiconti Lincei --
  Matematica E Applicazioni, 30 (2019), pp.~437--479.

\bibitem{correa2018nonlocal}
{\sc E.~Correa and A.~De~Pablo}, {\em Nonlocal operators of order near zero},
  Journal of Mathematical Analysis and Applications, 461 (2018), pp.~837--867.

\bibitem{davoli2023local}
{\sc E.~Davoli, E.~Rocca, L.~Scarpa, and L.~Trussardi}, {\em Local asymptotics
  and optimal control for a viscous {Cahn-Hilliard-Reaction-Diffusion} model
  for tumor growth}, arXiv preprint
  \href{http://arxiv.org/abs/2311.10457}{arXiv:2311.10457},  (2023).

\bibitem{di2016local}
{\sc A.~Di~Castro, T.~Kuusi, and G.~Palatucci}, {\em Local behavior of
  fractional $p$-minimizers}, Annales de l'Institut Henri Poincar{\'e} C,
  Analyse non lin{\'e}aire, 33 (2016), pp.~1279--1299.

\bibitem{di2012hitchhiker}
{\sc E.~Di~Nezza, G.~Palatucci, and E.~Valdinoci}, {\em Hitchhiker's guide to
  the fractional {Sobolev} spaces}, Bulletin des Sciences Math{\'e}matiques,
  136 (2012), pp.~521--573.

\bibitem{dipierro2017nonlocal}
{\sc S.~Dipierro, X.~Ros-Oton, and E.~Valdinoci}, {\em Nonlocal problems with
  {Neumann} boundary conditions}, Revista Matem{\'a}tica Iberoamericana, 33
  (2017), pp.~377--416.

\bibitem{doktor1976approximation}
{\sc P.~Doktor}, {\em Approximation of domains with {Lipschitzian} boundary},
  {\v{C}}asopis pro p{\v{e}}stov{\'a}ni matematiky, 101 (1976), pp.~237--255.

\bibitem{drabek2011manifolds}
{\sc P.~Dr{\'a}bek, R.~F. Man{\'a}sevich, and P.~Tak{\'a}c}, {\em Manifolds of
  critical points in a quasilinear model for phase transitions}, Contemporary
  Mathematics, 540 (2011), pp.~95--134.

\bibitem{dyda2020regularity}
{\sc B.~Dyda and M.~Kassmann}, {\em Regularity estimates for elliptic nonlocal
  operators}, Analysis \& PDE, 13 (2020), pp.~317--370.

\bibitem{evans2010partial}
{\sc L.~C. Evans}, {\em Partial Differential Equations}, Graduate Studies in
  Mathematics, American Mathematical Society, 2010.

\bibitem{fall2022regional}
{\sc M.~M. Fall}, {\em Regional fractional laplacians: boundary regularity},
  Journal of Differential Equations, 320 (2022), pp.~598--658.

\bibitem{felsinger2015dirichlet}
{\sc M.~Felsinger, M.~Kassmann, and P.~Voigt}, {\em The {Dirichlet} problem for
  nonlocal operators}, Mathematische Zeitschrift, 279 (2015), pp.~779--809.

\bibitem{fiscella2015density}
{\sc A.~Fiscella, R.~Servadei, and E.~Valdinoci}, {\em Density properties for
  fractional {Sobolev} spaces}, Annales Academi\ae\ Scientiarum Fennic\ae\
  Mathematica, 40 (2015), pp.~235--253.

\bibitem{folino2022generalized}
{\sc R.~Folino, L.~F.~L. R{\'\i}os, and M.~Strani}, {\em On a generalized
  {Cahn-Hilliard} model with $p$-laplacian}, Advances in Differential
  Equations, 27 (2022), pp.~647--682.

\bibitem{gajewski2003nonlocal}
{\sc H.~Gajewski and K.~Zacharias}, {\em On a nonlocal phase separation model},
  Journal of Mathematical Analysis and Applications, 286 (2003), pp.~11--31.

\bibitem{gal2017strong}
{\sc C.~G. Gal}, {\em On the strong-to-strong interaction case for doubly
  nonlocal {Cahn-Hilliard} equations}, Discrete and Continuous Dynamical
  Systems, 37 (2017), pp.~131--167.

\bibitem{gal2018doubly}
\leavevmode\vrule height 2pt depth -1.6pt width 23pt, {\em Doubly nonlocal
  {Cahn–Hilliard} equations}, Annales de l'Institut Henri Poincar{\'e} C,
  Analyse non lin{\'e}aire, 35 (2018), pp.~357--392.

\bibitem{gal2023separation}
{\sc C.~G. Gal, A.~Giorgini, and M.~Grasselli}, {\em The separation property
  for $2d$ {Cahn-Hilliard} equations: local, nonlocal and fractional energy
  cases}, Discrete and Continuous Dynamical Systems, 43 (2023), pp.~2270--2304.

\bibitem{garofalo2017fractional}
{\sc N.~Garofalo}, {\em Fractional thoughts}, arXiv preprint
  \href{http://arxiv.org/abs/1712.03347}{arXiv:1712.03347},  (2017).

\bibitem{giacomin1997phase}
{\sc G.~Giacomin and J.~L. Lebowitz}, {\em Phase segregation dynamics in
  particle systems with long range interactions. {I}. {Macroscopic} limits},
  Journal of Statistical Physics, 87 (1997), pp.~37--61.

\bibitem{gilbarg2001elliptic}
{\sc D.~Gilbarg and N.~S. Trudinger}, {\em Elliptic Partial Differential
  Equations of Second Order}, Classics in Mathematics, Springer-Verlag, Berlin,
  2001.
\newblock Reprint of the 1998 edition.

\bibitem{heinonen2006nonlinear}
{\sc J.~Heinonen, T.~Kilpel{\"a}inen, and O.~Martio}, {\em Nonlinear Potential
  Theory of Degenerate Elliptic Equations}, Dover Publications, 2006.

\bibitem{kassmann2009priori}
{\sc M.~Kassmann}, {\em A priori estimates for integro-differential operators
  with measurable kernels}, Calculus of Variations and Partial Differential
  Equations, 34 (2009), pp.~1--21.

\bibitem{kenmochi1995subdifferential}
{\sc N.~Kenmochi, M.~Niezgodka, and I.~Pawlow}, {\em Subdifferential operator
  approach to the {Cahn-Hilliard} equation with constraint}, Journal of
  Differential Equations, 117 (1995), pp.~320--356.

\bibitem{korvenpaa2017fractional}
{\sc J.~Korvenp{\"a}{\"a}, T.~Kuusi, and G.~Palatucci}, {\em Fractional
  superharmonic functions and the {Perron} method for nonlinear
  integro-differential equations}, Mathematische Annalen, 369 (2017),
  pp.~1443--1489.

\bibitem{kuusi2015nonlocal}
{\sc T.~Kuusi, G.~Mingione, and Y.~Sire}, {\em Nonlocal equations with measure
  data}, Communications in Mathematical Physics, 337 (2015), pp.~1317--1368.

\bibitem{kwasnicki2017ten}
{\sc M.~Kwa{\'s}nicki}, {\em Ten equivalent definitions of the fractional
  {L}aplace operator}, Fractional Calculus and Applied Analysis, 20 (2017),
  pp.~7--51.

\bibitem{lisini2012cahn}
{\sc S.~Lisini, D.~Matthes, and G.~Savar{\'e}}, {\em {Cahn-Hilliard} and thin
  film equations with nonlinear mobility as gradient flows in
  weighted-{Wasserstein} metrics}, Journal of Differential Equations, 253
  (2012), pp.~814--850.

\bibitem{lombardini2018approximation}
{\sc L.~Lombardini}, {\em Approximation of sets of finite fractional perimeter
  by smooth sets and comparison of local and global $s$-minimal surfaces},
  Interfaces and Free Boundaries, 20 (2018), pp.~261--296.

\bibitem{lombardini2018minimization}
\leavevmode\vrule height 2pt depth -1.6pt width 23pt, {\em Minimization
  problems involving nonlocal functionals: nonlocal minimal surfaces and a free
  boundary problem}, PhD thesis, Universit\`a degli Studi di Milano and
  Universit\'e de Picardie Jules Verne, 2018.
\newblock Available at
  \href{http://arxiv.org/abs/1811.09746}{arXiv:1811.09746}.

\bibitem{miranville2017cahn}
{\sc A.~Miranville}, {\em The {Cahn-Hilliard} equation and some of its
  variants}, AIMS Mathematics, 2 (2017), pp.~479--544.

\bibitem{mugnai2019neumann}
{\sc D.~Mugnai and E.~Proietti~Lippi}, {\em Neumann fractional $p$-{Laplacian}:
  {Eigenvalues} and existence results}, Nonlinear Analysis, 188 (2019),
  pp.~455--474.

\bibitem{nec2008front}
{\sc Y.~Nec, A.~Nepomnyashchy, and A.~Golovin}, {\em Front-type solutions of
  fractional {Allen-Cahn} equation}, Physica D: Nonlinear Phenomena, 237
  (2008), pp.~3237--3251.

\bibitem{ponce2004estimate}
{\sc A.~C. Ponce}, {\em An estimate in the spirit of {Poincar{\'e}}'s
  inequality}, Journal of the European Mathematical Society, 6 (2004),
  pp.~1--15.

\bibitem{ponce2004new}
\leavevmode\vrule height 2pt depth -1.6pt width 23pt, {\em A new approach to
  {Sobolev} spaces and connections to ${\Gamma}$-convergence}, Calculus of
  Variations and Partial Differential Equations, 19 (2004), pp.~229--255.

\bibitem{saldana2018fractional}
{\sc A.~Salda{\~n}a}, {\em On fractional higher-order {Dirichlet} boundary
  value problems: between the {Laplacian} and the bilaplacian}, Mexican
  Mathematicians in the World,  (2018).

\bibitem{santambrogio2017euclidean}
{\sc F.~Santambrogio}, {\em $\{$Euclidean, metric, and Wasserstein$\}$ gradient
  flows: an overview}, Bulletin of Mathematical Sciences, 7 (2017),
  pp.~87--154.

\bibitem{servadei2014spectrum}
{\sc R.~Servadei and E.~Valdinoci}, {\em On the spectrum of two different
  fractional operators}, Proceedings of the Royal Society of Edinburgh: Section
  A Mathematics, 144 (2014), pp.~831--855.

\bibitem{simon1986compact}
{\sc J.~Simon}, {\em Compact sets in the space {$L^p(0,T;B)$}}, Annali di
  Matematica Pura ed Applicata, 146 (1986), pp.~65--96.

\bibitem{takavc2009stationary}
{\sc P.~Tak{\'a}{\v{c}}}, {\em Stationary radial solutions for a quasilinear
  {Cahn-Hilliard} model in space dimensions.}, Electronic Journal of
  Differential Equations, 2009 (2009), pp.~227--254.

\bibitem{vazquez2014recent}
{\sc J.-L. V{\'a}zquez}, {\em Recent progress in the theory of nonlinear
  diffusion with fractional {Laplacian} operators}, Discrete and Continuous
  Dynamical Systems - S, 7 (2014), pp.~857--885.

\bibitem{ziemer2017modern}
{\sc W.~P. Ziemer}, {\em Modern Real Analysis}, Graduate Texts in Mathematics,
  Springer International Publishing, 2017.

\end{thebibliography}

\end{document}